%% file: main.tex
\documentclass{article}
\usepackage[twoside,
paperwidth=210mm,
paperheight=297mm,
textheight=622pt,
textwidth=468pt,
centering,
headheight=50pt,
headsep=12pt,
footskip=18pt,
footnotesep=24pt plus 2pt minus 12pt,
columnsep=2pc]{geometry}

\usepackage[auth-sc]{authblk}

\usepackage{mathtools}

\allowdisplaybreaks

\usepackage{amssymb}
\usepackage[mathscr]{eucal}

\usepackage{mismath}

\usepackage{hyperref}

\usepackage{bm} % Bold math font.

\usepackage{graphicx}
\graphicspath{{figs/}}
\usepackage{pgf}
\usepackage{pgfplots}
\pgfplotsset{compat=1.14}

\usepackage{tikz}
\usetikzlibrary{math}
\usetikzlibrary{calc}

\usepackage{subcaption}

\usepackage{multirow}

\usepackage{makecell}
\setcellgapes{2pt}

\usepackage{booktabs}

\usepackage{siunitx}
\usepackage{algorithm}
\usepackage{algorithmic}

\usepackage[capitalise]{cleveref}

\newcommand{\dx}{\,\mathrm{d}}
\newcommand{\Real}{\mathbb{R}}

\DeclareMathOperator{\Dist}{dist}
\DeclareMathOperator{\Span}{span}

\DeclareMathOperator{\Diag}{diag}

\DeclarePairedDelimiter{\RoundBrackets}{(}{)}
\DeclarePairedDelimiter{\CurlyBrackets}{\{}{\}}
\DeclarePairedDelimiter{\SquareBrackets}{[}{]}

\usepackage{stmaryrd}
\DeclarePairedDelimiter{\DSquareBrackets}{\llbracket}{\rrbracket}
\DeclareSymbolFont{sfletters}{OML}{cmbrm}{m}{it}

\DeclareMathSymbol{\salpha}{\mathord}{sfletters}{"0B}
\DeclareMathSymbol{\sbeta}{\mathord}{sfletters}{"0C}
\DeclareMathSymbol{\sgamma}{\mathord}{sfletters}{"0D}
\DeclareMathSymbol{\sdelta}{\mathord}{sfletters}{"0E}
\DeclareMathSymbol{\sepsilon}{\mathord}{sfletters}{"0F}
\DeclareMathSymbol{\szeta}{\mathord}{sfletters}{"10}
\DeclareMathSymbol{\seta}{\mathord}{sfletters}{"11}
\DeclareMathSymbol{\stheta}{\mathord}{sfletters}{"12}
\DeclareMathSymbol{\siota}{\mathord}{sfletters}{"13}
\DeclareMathSymbol{\skappa}{\mathord}{sfletters}{"14}
\DeclareMathSymbol{\slambda}{\mathord}{sfletters}{"15}
\DeclareMathSymbol{\smu}{\mathord}{sfletters}{"16}
\DeclareMathSymbol{\snu}{\mathord}{sfletters}{"17}
\DeclareMathSymbol{\sxi}{\mathord}{sfletters}{"18}
\DeclareMathSymbol{\spi}{\mathord}{sfletters}{"19}
\DeclareMathSymbol{\srho}{\mathord}{sfletters}{"1A}
\DeclareMathSymbol{\ssigma}{\mathord}{sfletters}{"1B}
\DeclareMathSymbol{\stau}{\mathord}{sfletters}{"1C}
\DeclareMathSymbol{\supsilon}{\mathord}{sfletters}{"1D}
\DeclareMathSymbol{\sphi}{\mathord}{sfletters}{"1E}
\DeclareMathSymbol{\schi}{\mathord}{sfletters}{"1F}
\DeclareMathSymbol{\spsi}{\mathord}{sfletters}{"20}
\DeclareMathSymbol{\somega}{\mathord}{sfletters}{"21}
\DeclareMathSymbol{\svarepsilon}{\mathord}{sfletters}{"22}
\DeclareMathSymbol{\svartheta}{\mathord}{sfletters}{"23}
\DeclareMathSymbol{\svarpi}{\mathord}{sfletters}{"24}
\DeclareMathSymbol{\svarrho}{\mathord}{sfletters}{"25}
\DeclareMathSymbol{\svarsigma}{\mathord}{sfletters}{"26}
\DeclareMathSymbol{\svarphi}{\mathord}{sfletters}{"27}
\DeclareMathSymbol{\sDelta}{\mathord}{sfletters}{"01}
\DeclareMathSymbol{\sTheta}{\mathord}{sfletters}{"02}

\usepackage{amsthm}
\newtheorem{theorem}{Theorem}[section]
\newtheorem{lemma}[theorem]{Lemma}

\theoremstyle{definition}
\newtheorem{definition}{Definition}

\crefname{assumption}{assumption}{assumptions}
\Crefname{assumption}{Assumption}{Assumptions}

\crefname{problem}{problem}{problems}
\Crefname{problem}{Problem}{Problems}

\theoremstyle{remark}
\newtheorem*{remark}{Remark}

\usepackage{lineno}
\usepackage{easyReview}

\title{Learning a generalized multiscale prolongation operator}
\author[1]{Yucheng Liu}
\author[2]{Shubin Fu}
\author[1]{Yingjie Zhou}
\author[1]{Changqing Ye}
\author[1]{Eric~T.~Chung}
\affil[1]{Department of Mathematics, The Chinese University of Hong Kong, Shatin, Hong~Kong~SAR, China.}
\affil[2]{Eastern Institute for Advanced Study, Ningbo, China.}

\date{}
 
\begin{document}
\maketitle

\begin{abstract}
In this research, we address Darcy flow problems with random permeability using iterative solvers, enhanced by a two-grid preconditioner based on a generalized multiscale prolongation operator, which has been demonstrated to be stable for high contrast profiles.
To circumvent the need for repeatedly solving spectral problems with varying coefficients, we harness deep learning techniques to expedite the construction of the generalized multiscale prolongation operator.
Considering linear transformations on multiscale basis have no impact on the performance of the preconditioner, we devise a loss function by the coefficient-based distance between subspaces instead of the plain $l^2$-norm of the difference of the corresponding multiscale bases.
We discover that leveraging the inherent symmetry in the local spectral problem can effectively accelerate the neural network training process.
In scenarios where training data are limited, we utilize the Karhunen--Lo\`eve expansion to augment the dataset.
Extensive numerical experiments with various types of random coefficient models are exhibited, showing that the proposed method can significantly reduce the time required to generate the prolongation operator while maintaining the original efficiency of the two-grid preconditioner.
Notably, the neural network demonstrates strong generalization capabilities, as evidenced by its satisfactory performance on unseen random permeability fields.


\end{abstract}

\section{Introduction}
Simulating subsurface fluid flow in heterogeneous porous media is frequently encountered in real-world applications, including reservoir simulation, groundwater resource management, groundwater contamination predictions, and the exploration of oil and gas fields \cite{bear2010modeling, national1996rock}.
The properties of a heterogeneous porous medium can vary significantly across different spatial regions and exhibit complex structures, making it essential to carefully account for heterogeneity in simulations \cite{zimmerman1996hydraulic}.
These detailed models result in large linear systems, which in turn lead to computational challenges.

Model reduction techniques, such as upscaling based on homogenization \cite{3.7,3.19,2.50,3.19,2.6}, have been rapidly developed.
Upscaling is an averaging process that scales the static properties of a fine-grid model to equivalent characteristics at a coarse-grid level.
This process aims to ensure that the coarse and fine models share a similar formulation.
However, it often fails to capture all the small-scale details of the quantities of interest.
Another approach is multiscale model reduction \cite{chung2023multiscale}, which include the multiscale finite element method \cite{2.33}, the mixed multiscale finite element method \cite{2.15,2.1}, the multiscale finite volume method \cite{2.39, 2.27}, the multiscale mortar mixed finite element method \cite{2.4} and the generalized multiscale finite element method \cite{2.21, 2.17}.
The main idea behind these methods is to reduce the problem defined on fine grids to one on coarse grids.
The multiscale basis functions encapsulate the heterogeneity of the original geological medium and are defined by solving a series of carefully designed local problems on each coarse grid element.
While these multiscale methods have been successfully applied to various multiscale models, some of the effectiveness can diminish with increasing contrast in permeability and correlation length.
When highly precise fine-scale solutions are required, it becomes essential to design an efficient solver that can directly handle the original fine-scale information, rather than relying solely on multiscale methods.
Multigrid \cite{1} is one of the most efficient and commonly used algorithms for solving large-scale linear systems \cite{2}.
The smoother and the connection (prolongation operator) between different levels are therefore crucial components.
The mathematical theory shows that the optimal prolongation operator is defined on eigenvalue problems for fixed smoothing operation \cite{Xu2017}.
We focus on utilizing neural networks to learn such an optimal prolongation operator.

With the rapid advancement of deep learning, it has become a research hotspot for the numerical solution of differential equations in recent years. Numerous methods have been proposed to apply deep learning directly to solving PDEs.
A notable development in this area is the introduction of physics-informed neural networks (PINNs) \cite{5.3}, which have gained significant attention.
However, PINNs are limited in their applicability as they learn the solution for a single problem instance.
Another popular research area is operator learning \cite{5.7,5.8,5.9}, which demonstrates improved generalization abilities by focusing on the differential operators rather than one solution determined by the PDEs, initial value conditions and boundary conditions.
Compared to traditional algorithms for solving PDEs, neural network-based methods are highly demanding for training data and significant computational resource and have unresolved issues with numerical accuracy and stability.
For the interplay between multiscale methods and neural networks, we name a few: a deep learning strategy is proposed to determine multiscale basis functions and coarse-scale parameters for non-static coefficients in \cite{r2}; a new deep neural networks approach tailored for model reduction in multiscale porous media flow problems has been developed in \cite{r3}.
Attention has also been paid to combine learning technique with traditional linear solvers: geometric multigrid methods have been advanced through the application of evolutionary computation for optimizing their techniques \cite{6.30} and the deep learning approaches for tuning prolongation and restriction operators in two-grid algorithms \cite{6.21, 6.14}; graph based neural networks have been adapted to algebraic multigrid methods, enabling the performance on unstructured problems \cite{6.24}.

In this study, we explore the application of neural networks in simulating Darcy flow, focusing specifically on the prolongation operator related to generalized multiscale basis functions derived from the coefficient. We design our neural networks architecture with convolutional layers complemented by upsampling and downsampling layers, which significantly streamline computational and spatial complexities.
We have found that exploiting the inherent symmetry in local spectral problems (LSPs) can substantially expedite the training of our neural networks.
When faced with limited training data, the Karhunen--Lo\`eve expansion proves effective for data augmentation purposes.
Linear transformations on generalized multiscale basis do not affect the performance of the preconditioner; therefore, we propose using the coefficient-based distance between subspaces as our loss function instead of $l^{2}$-norm of the difference of the corresponding generalized multiscale basis and provide a theoretical proof to demonstrate its efficacy.
Ultimately, our approach enables the direct computation of the prolongation operator, thus obviating the need for the previously complex method of training individual neural networks for each basis function.
We underscore that the operators derived via our proposed data-driven deep learning approach exhibit universality, applicable across diverse boundary conditions, transcending the constraints of singular problem environments. 
Moreover, we test the generalization ability of neural networks which is sufficiently strong to be applicable to permeability profiles from unseen datasets and coarse elements of varying resolutions.

The rest of this paper is organized as follows. \Cref{prel} introduces the model, discretization, and the iterative solver.
\Cref{nn} details the deep learning method, defines our loss function by coefficient-based distance between subspaces, proves its efficacy and presents data augmentation methods including inherent symmetry in the LSPs and Karhunen--Lo\`eve expansion.
Experimental results are showcased in \cref{nume}.
Finally, we synthesize the key results and provided outlooks for future research in \cref{con}.

\section{Preliminaries}
\label{prel}
\subsection{Model problem and numerical discretization}
We consider the following Darcy flow model of a single-phase fluid in a potentially high-contrast and random medium:
\begin{equation}\label{eq:orgional_equation}
  \left\{
  \begin{aligned}
    \kappa^{-1}\bm{u}+\nabla p & =\bm{0}\quad \text{in}\quad \Omega, \\
    \nabla \cdot \bm{u}        & =f \quad \text{in}\quad \Omega,
  \end{aligned}
  \right.
\end{equation}
where $\Omega$ is the computational domain, $\kappa$ (scalar-valued random field) is the permeability of the porous medium, $\bm{u}$ is the Darcy velocity, $p$ is the pressure field and $f$ is the source term.

For specific numerical discretization, a boundary condition must be established.
However, it should be emphasized that our method can be applied for different boundary conditions.
Here, we present the numerical discretization under the no-flux boundary condition.
Consider the two-point flux approximation (TPFA) to discrete \cref{eq:orgional_equation}.
Let $\mathcal{T}^{h}$ be a quadrilaterals partition of the computational domain $\Omega$ with mesh size $h_{x}$ in x-direction and $h_{y}$ in y-direction, $N$ is the number of the element in mesh, $\mathcal{E}^{h}$ is the set of the internal edges of mesh.
By applying mixed finite element method (the lowest-order polynomial finite spaces $W_{h}$ and the lowest-order Raviart-Thomas finite element space $\bm{V}_{h}$ \cite{boffi2013mixed}) with velocity elimination technique \cite{ye2024robust} and the trapezoidal quadrature rule, we achieve
\begin{equation}\label{eq:apxfinal}
  \sum_{e\in\mathcal{E}^h}\kappa_{e} \DSquareBrackets{p_h}_{e}
  \DSquareBrackets{w_h}_{e}\frac{\abs{e}^{2}}{h_{x}h_{y}}= \int_{\Omega} f w_{h}
  \dx \bm{x}, \quad \forall w_{h} \in W_{h},
\end{equation}
where $\kappa_{e}^{-1}= \RoundBrackets{\kappa_{e, +}^{-1} + \kappa_{e, -}^{-1}}/{2}$, $\RoundBrackets{\kappa_{e,-}, \kappa_{e,+}}$ is a neighbor pair of an edge $e \in \mathcal{E}^{h}$, $\abs{e}$ is the length of the edge $e$ and $\DSquareBrackets{p_h}_{e}$, $\DSquareBrackets{w_h}_{e}$ are the jump of element-wisely constant functions $p_{h}$, $w_{h}$ across edge $e$ towards the positive direction.

We can express \cref{eq:apxfinal} as a linear system of equations:
\begin{equation}\label{eq:linearsystem}
  A P = F,
\end{equation}
where $A \in \Real^{N \times N}$, $P \in \Real^{N}$ and $F \in \Real^{N}$.

\subsection{Two-grid preconditioner}
We further elucidated the process of implementing the Preconditioned Conjugate Gradient (PCG) method \cite{YE2024116982} with the two-grid preconditioner $B^{-1}$ based on generalized multiscale basis to effectively solve the linear system defined in \cref{eq:linearsystem}.
In the subsequent discussion, we have utilized generalized multiscale space notation to highlight the importance of geometric information, and a notation with ``$c$'' superscript is associated with the coarse grid.
To grasp the key ingredients directly, we introduce several key notations and their specific roles in a two-grid preconditioner:
\begin{itemize}
  \item $R$ serves as the smoother in fine space $W_{h}$.
  \item $B^{-1}$ acts as a preconditioner of $A$.
  \item $P$ is a prolongation operator which links coarse space $W_{h}^{c}$ and fine space $W_{h}$.
  \item $B^{c} = \RoundBrackets{A^c}^{-1}= \RoundBrackets{P^\intercal A P}^{-1}$ functions as a solver in coarse space $W_{h}^{c}$.
\end{itemize}
The two-grid preconditioner can be divided into three main parts: pre-smoothing, coarse grid correction, and post-smoothing.
Each stage targets different characteristics of the problem to enhance the overall efficiency and accuracy of the solution.
The detailed representation of these phases is shown in \cref{alg:tgp}.

\begin{algorithm}
  [!ht]
  \caption{A two-grid preconditioner $B^{-1}$}
  \label{alg:tgp} Given right-hand side term $r$:
  \begin{algorithmic}[1]
    \STATE Pre-smoothing: $w \gets Rr$.
    \STATE Coarse grid correction: $w \gets w + P B^{c} P^{\intercal}\RoundBrackets{r - A w}$.
    \STATE Post-smoothing: $B^{-1}r \gets w + R \RoundBrackets{r - A w}$.
  \end{algorithmic}
\end{algorithm}

\subsection{Prolongation operators}
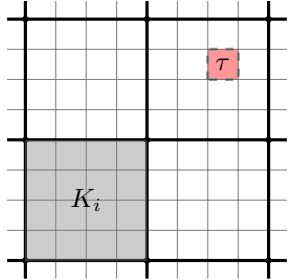
\begin{figure}[!ht]
  \centering
  \begin{tikzpicture}[scale=0.8]
    \draw[step=0.5, gray, thin] (-2.3, -2.3) grid (2.3, 2.3);
    \draw[step=2.0, black, very thick] (-2.3, -2.3) grid (2.3, 2.3);
    \foreach \x in {-1,...,1} \foreach \y in {-1,...,1}{ \fill (2.0 * \x, 2.0 * \y) circle (1.5pt); }
    \fill[gray, opacity=0.4] (-2.0, -2.0) rectangle (0.0, 0.0);
    \node at (-1.0, -1.0) {$K_{j}$};
    \draw[dashed, very thick, fill=red, opacity=0.4]
    (1.0, 1.0) rectangle (1.5, 1.5);
    \node at (1.25, 1.25) {$\tau$};
  \end{tikzpicture}
  \caption{An illustration of the two-scale mesh: $\tau$ is a fine element and $K_{j}$ is a coarse element which contains $16$ fine elements.}
  \label{fig:grid}
\end{figure}

Assume the domain $\Omega$ is divided into a union of disjoint subdomains denoted by $\CurlyBrackets{K_j}_{j=1}^{n}$ with each $K_{j}$ consisting of $N/n$ fine elements from $\mathcal{T}_{h}$.
Denote $\mathcal{T}_{H}$ as the coarse quadrilaterals partition of the computational domain $\Omega$.
The basis functions of finite space $W_{h}$ restricted to the coarse element $K_{j}$ form the snapshot space $W_{h}\RoundBrackets{K_j}$.
\Cref{fig:grid} is an illustration of the two-scale mesh where $\tau$ is a fine element and $K_{j}$ is a coarse element which contains $16$ fine elements.

We introduce the local generalized multiscale basis functions $\CurlyBrackets{\Phi_{j,k}}_{k=1}^{n^c}$ associated with the coarse element $K_{j}$, where $n^{c}$ is the number of basis functions on each coarse element.
In GMsFEMs \cite{2.21}, two bilinear forms corresponding to the $j$-th coarse element need to be defined to obtain effective low-dimensional representations from snapshot space $W_{h}\RoundBrackets{K_j}$.
Here we define
\begin{equation}\label{eq:two_bilinear}
  \begin{aligned}
    a_{j}(w_{h}', w_{h}'') & \coloneqq \sum_{e \in \mathcal{E}^h(K_j)}\kappa_{e} \left\llbracket w_{h}'\right\rrbracket_{e} \left\llbracket w_{h}''\right\rrbracket_{e} \frac{\lvert e \rvert^2}{h_x h_y}, \quad \forall w_{h}', w_{h}''\in W_{h}(K_{j}), \\
    s_{j}(w_{h}', w_{h}'') & \coloneqq \int_{K_j}\tilde{\kappa}w_{h}'w_{h}''\, d\bm{x}, \quad \forall w_{h}', w_{h}''\in W_{h}(K_{j}),
  \end{aligned}
\end{equation}
where $\mathcal{E}^{h}(K_{j})$ is the internal edge set in the coarse element $K_{j}$.
It is reasonable to simply set $\tilde{\kappa}=\kappa$, see \cite{ye2024robust}.
The local generalized multiscale basis functions $\CurlyBrackets{\Phi_{j,k}}_{k=1}^{n^c}$ is defined by solving the following spectral problem on coarse element $K_{j}$:
\begin{equation}\label{eq:spepb}
  \sum_{e\in \mathcal{E}^h(K_j)}\kappa_{e}\DSquareBrackets{\Phi_h}_{e}\DSquareBrackets{w_h}_{e}\frac{\abs{e}^{2}}{h_{x}h_{y}}= \lambda \int_{K_j} \tilde{\kappa}\Phi_{h} w_{h} \dx \bm{x}, \quad \forall w_{h} \in W_{h}(K_{j}).
\end{equation}
Solving the LSP \cref{eq:spepb}, eigenvectors $\CurlyBrackets{\Phi_{j,k}}_{k=1}^{n^c}$ corresponding to the $n^{c}$ smallest eigenvalues will be the local generalized multiscale basis functions.
The direct sum of the basis functions on each coarse element spans the global coarse space and its projection onto the fine space is represented as a prolongation operator $P$, where each column of $P$ corresponds to a generalized multiscale basis and denote the number of column of $P$ as $N^c = n \times n^c$.

\section{Methodology} \label{nn}
The coefficient $\kappa$ is drawn from a random field, indicating that for any $\kappa$ selected from a specific distribution, multiple spectral problems must be recalculated in order to derive the prolongation operator.
The process of solving spectral problems can be viewed as extracting several features of the same size from a coarse element, which is quite similar to panoptic segmentation \cite{kirillov2019panopticsegmentation}.
Therefore, we aim to propose a neural network architecture that takes a coarse element as input and directly outputs the corresponding part of prolongation operator.

\subsection{Neural network architectures}
\begin{figure}[!ht]
  \centering
  \includegraphics[width=\textwidth]{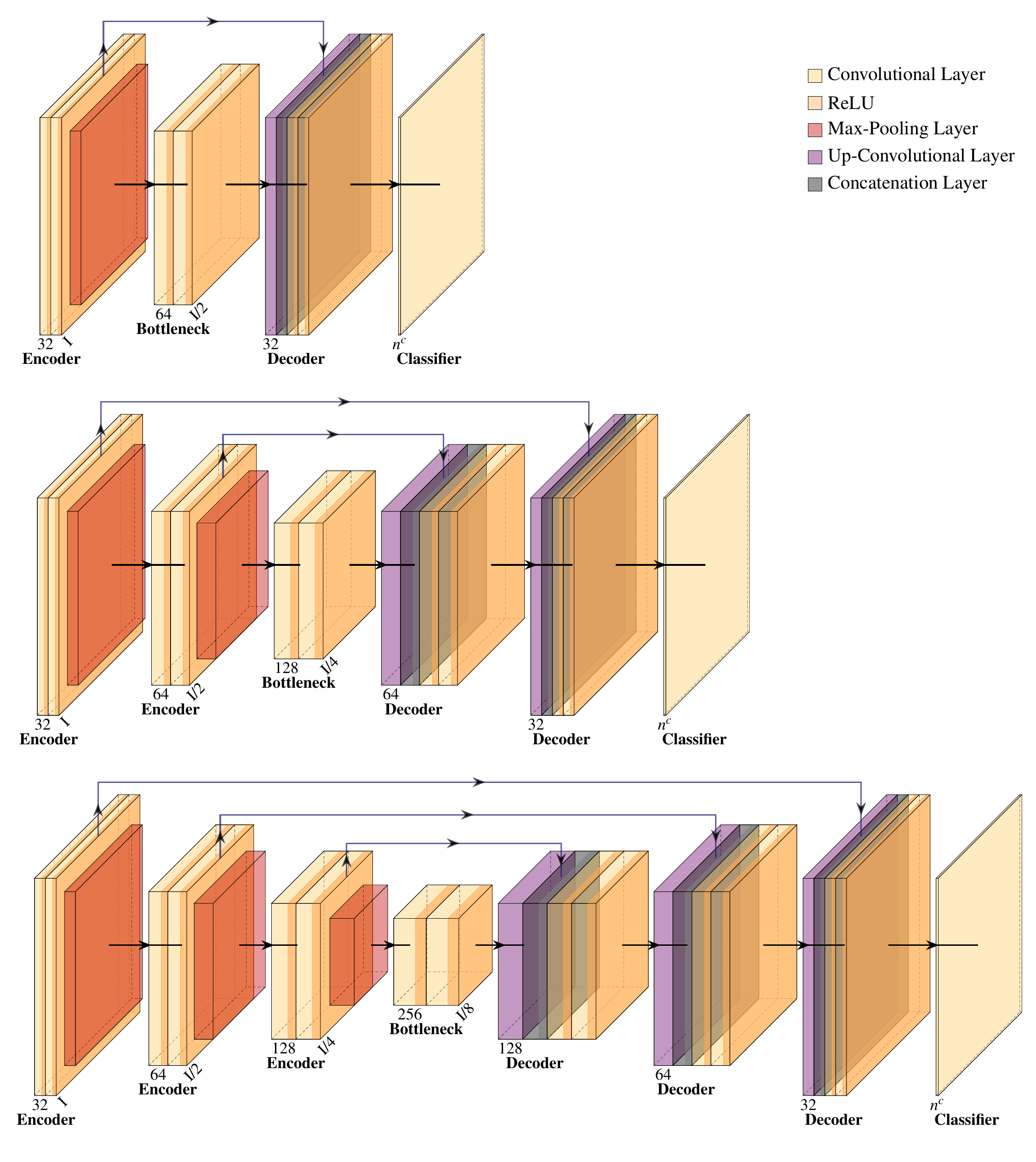}
  \caption{Two-level, three-level and four-level U-Net architectures.}
  \label{fig:Arch}
\end{figure}

Leveraging the efficiency of the U-Net \cite{Unet}, which excels in panoptic segmentation tasks, we propose multiple level neural networks shown in \cref{fig:Arch} specifically designed for spectral problems.
In the diagram, the labels along the x-axis represent the number of channels in the corresponding layer, while the labels along the z-axis indicate the dimensions of the image in that layer.

The multiple level neural networks contain consists of the encoder part (contracting path), a bottleneck, the decoder part (expansive path) and a classifier.
The specifics include:
\begin{itemize}
  \item \textbf{Convolutional Layers:} Two consecutive convolution processes using a kernel size of $3 \times 3$ and padding of $1$ are implemented, maintaining the output feature map size unchanged.

  \item \textbf{Pooling Layers:} After the convolutional operations, a max pooling layer with a pool size of $2 \times 2$ and a stride of $2$ is used to halve the spatial dimensions of the feature maps, providing a form of translation invariance.

  \item \textbf{Upconvolutional Layers:} Initially, a convolution operation with a kernel size of $3 \times 3$, stride of $2$, a padding of $1$, and output padding of $1$ is used to double the spatial dimensions of the feature maps.
\end{itemize}

The encoder and decoder are directly connected through skip connections.
These connections transfer feature maps from the encoder directly to corresponding layers in the decoder, assisting the decoder in better utilizing both high-level features and fine details transferred from the encoder, thereby more accurately reconstructing and localizing the input data.

\subsection{Loss function}
The prolongation operator $P$ is a block diagonal matrix.
Each block is assembled from generalized multiscale basis functions of a coarse element. Since we aim to learn each block of prolongation operator directly, a special distance between two subspaces of $\Real^{N/n}$ is proposed based on coefficient $\kappa$ on a coarse element.
In our neural networks, we take $\Dist_{\tilde{\kappa}}\RoundBrackets{\text{target}, \text{prediction}}$ as our loss function.
Specific definition of the distance is given below.

\begin{definition}
  \label{defdist} Denote
  \[
    P_{j} = \SquareBrackets*{p_{j, 1},\ p_{j, 2},\dots, p_{j, n^c}}, \quad
    X_{j} = \Span \CurlyBrackets*{p_{j, 1},\ p_{j, 2}, \dots, p_{j,n^c}},
  \]
  where $P_{j}$ is a block of prolongation operator based on the $j$-th coarse element, $X_{j}$ is a subspace of $\Real^{N/n}$ and $\CurlyBrackets{\psi_{j, 1},\ \psi_{j, 2}, \dots, \psi_{j, n^{c}}}$ is an orthonormal basis of subspace $X_{j}$ with respectively to the inner product $\RoundBrackets{\cdot, \cdot}_{\tilde{\kappa}}$.
  These basis functions are arranged column-wise to form the matrix denoted as $T_{j}$.
  For two different prolongation operators $P_{j}^{(1)}$ and $P_{j}^{(2)}$, define
  \[
    \Dist_{\tilde{\kappa}_j}\RoundBrackets{P_{j}^{(1)}, P_{j}^{(2)}}= \RoundBrackets*{n^{c} - \norm{\RoundBrackets{T_{j}^{(1)}}^{\intercal} \tilde{\kappa}_j T_{j}^{(2)}}_F^{2}}^{1/2},
  \]
  where $\norm{\cdot }_{F}$ is Frobenius norm and the $\tilde{\kappa}_{j}$ here is an ${N/n \times N/n}$ diagonal matrix which is from $\tilde{\kappa}$ restricted on coarse element $K_{j}$.
\end{definition}

To simplify these expressions, we will omit the index $j$ for the coarse element index in subsequent discussions, whenever it does not lead to ambiguity.

\begin{figure}[!ht]
  \centering
  \begin{tikzpicture}[scale=1, x={(1cm,0cm)}, y={(0.5cm,0.2cm)}, z={(0cm,1cm)}]
    \coordinate (O) at (0,0,0);
    \coordinate (A) at (2,1,1);
    \coordinate (B) at (2,1,0);

    \filldraw[fill=gray!20, draw=gray, opacity=0.9]
    (-1,-1.7,0) --
    (4,-1.7,0) --
    (4,3.2,0) --
    (-1,3.2,0) --
    cycle;

    \draw[->, thick, blue] (O) -- (A) node[midway, above] {$\psi^{(2)}_{k}$};
    \draw[->, thick, red] (O) -- (B) node[midway, below] {$\text{Projection}$};

    \draw[dashed, thick] (A) -- (B);

    \node at (A) [right] {$A_{k}$};
    \node at (B) [below] {$B_{k}$};
    \node at (O) [left] {$O$};
    \node at (3.2,1,0) [above] {Subspace $X^{(1)}$};
  \end{tikzpicture}
  \caption{An intuitive diagram of distance defined in the loss function.
  $\overrightarrow{OA_{k}}$ is an orthonormal basis of subspace $X^{(2)}$, $\overrightarrow{OB_{k}}$ is the projection of $\overrightarrow{OA_k}$ onto the subspace $X^{(1)}$ and $\overrightarrow{A_{k}B_{k}}$ is the residual of $\overrightarrow{OA_{k}}$ relative to the subspace $X^{(1)}$.}
  \label{fig:Dist}
\end{figure}
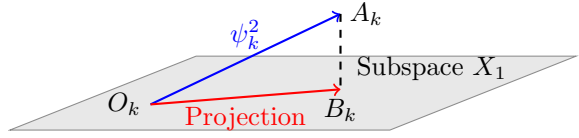

\begin{remark}
  An intuitive explanation of the above definition of distance in \cref{fig:Dist}.
  We first project all the basis vectors in an orthonormal basis of subspace $X^{(2)}$ onto subspace $X^{(1)}$ with respectively to $\RoundBrackets{\cdot, \cdot}_{\tilde{\kappa}}$.
  The square of the distance is equal to the sum of the squares of the lengths of residual $A_{k}B_{k}$.
\end{remark}

Next, we need to prove the distance in \cref{defdist} is well-defined.
\begin{theorem}
  The following properties hold for the distance in \cref{defdist}.
  \begin{itemize}
    \item The distance does not depend on the choice of basis.
    \item Non-negativity:
          $n^{c} - \norm{\RoundBrackets{T^{(1)}}^{\intercal} \tilde{\kappa} T^{(2)}}_{F}^{2} \geq 0$.

    \item Positive definiteness:
          $\Dist_{\tilde{\kappa}}\RoundBrackets{P^{(1)}, P^{(2)}}= 0$ if and only if $P^{(1)} = P^{(2)}$ in the subspace sense, which means $X^{(1)} = X^{(2)}$.

    \item Symmetry: $\Dist_{\tilde{\kappa}}\RoundBrackets{P^{(1)}, P^{(2)}}= \Dist_{\tilde{\kappa}}\RoundBrackets
            {P^{(2)}, P^{(1)}}$.

    \item Triangle inequality: $\Dist_{\tilde{\kappa}}\RoundBrackets{P^{(1)}, P^{(2)}}+ \Dist_{\tilde{\kappa}}\RoundBrackets
            {P^{(2)}, P^{(3)}}\geq \Dist_{\tilde{\kappa}}\RoundBrackets{P^{(1)}, P^{(3)}}$.
  \end{itemize}
\end{theorem}

\begin{proof}
  We can prove the properties one by one.
  \begin{itemize}
    \item Assume $\CurlyBrackets{\phi^{(i)}_1, \phi^{(i)}_2, \dots, \phi^{(i)}_{n^c}}$ is another orthonormal basis of subspace $X^{(i)}$ with respectively to $\RoundBrackets{\cdot, \cdot}_{\tilde{\kappa}}$ for $i \in {1, 2}$.
          Similarly, we define the orthogonal matrix $\hat{T}^{(i)}$.
          There is an orthogonal matrix $M_{i}$, such that $\hat{T}^{(i)} = T^{(i)} M^{(i)}$.
          Due to the orthogonal invariant of Frobenius norm, we can see that
          \begin{align*}
            \norm{(T^{(1)})^\intercal \tilde{\kappa} \hat{T}^{(2)}}_{F}^{2} & =  \norm{(T^{(1)})^\intercal \tilde{\kappa} T^{(2)} M^{(2)}}_{F}^{2} = \norm{(T^{(1)})^\intercal \tilde{\kappa} T^{(2)}}_{F}^{2}                   \\
                                                                            & = \norm{(M^{(1)})^\intercal (T^{(1)})^\intercal \tilde{\kappa} T^{(2)}}_{F}^{2} = \norm{(\hat{T}^{(1)})^\intercal \tilde{\kappa} T^{(2)}}_{F}^{2},
          \end{align*}
          which means that $\Dist_{\tilde{\kappa}}\RoundBrackets{P^{(1)}, P^{(2)}}$ does not depend on the choice of basis.

    \item Take $B_{kl}= \RoundBrackets{\psi^{(1)}_k, \psi^{(2)}_l}_{\tilde{\kappa}}$, we have
          \[
            \sum_{k=1}^{n^{c}}B_{kl}^{2} \leq \norm{ \psi^{(2)}_l}^{2}_{\tilde{\kappa}} =1 , \quad \sum_{k=1}^{n^{c}}\sum_{l=1}^{n^{c}}B_{kl}^{2} \leq n^{c}.
          \]
          Hence, it holds that
          \[
            n^{c} - \norm{ (T^{(1)}) ^\intercal \tilde{\kappa} T^{(2)}}_{F}^{2} = n^{c} - \norm{B}_{F}^{2} \geq 0.
          \]

    \item If $P^{(1)} = P^{(2)}$ in the subspace sense, we have
          \[
            \psi^{(2)}_{l} = \sum_{k=1}^{n^{c}}\RoundBrackets{\psi^{(1)}_k, \psi^{(2)}_l}_{\tilde{\kappa}}\psi^{(1)}_{k} , \quad \sum_{k=1}^{n^{c}}B_{kl}^{2}= \norm{\psi^{(2)}_l}^{2}_{\tilde{\kappa}}= 1, \quad \sum_{k=1}^{n^{c}}\sum_{l=1}^{n^{c}}B_{kl}^{2} = n^{c},
          \]
          which leads
          \[
            \Dist_{\tilde{\kappa}}\RoundBrackets{P^{(1)}, P^{(2)}}= \sqrt{n^{c} - \norm{(T^{(1)}) ^\intercal \tilde{\kappa} T^{(2)}}_{F}^{2}}= \sqrt{ n^{c} - \norm{B}_{F}^{2}}= 0.
          \]
          If
          \[
            \Dist_{\tilde{\kappa}}\RoundBrackets{P^{(1)}, P^{(2)}}= \sqrt{n^{c} - \norm{ (T^{(1)}) ^\intercal \tilde{\kappa}  T^{(2)}}_{F}^{2}}= \sqrt{\sum_{l=1}^{n^{c}}\RoundBrackets*{1 - \sum_{k=1}^{n^{c}} B_{kl}^2}}= 0.
          \]
          Due to for all $l \in \CurlyBrackets{1,2,3,\dots, n^c}$, $1 - \sum_{k=1}^{n^{c}}B_{kl}^{2} \geq 0$, we hence have $1 - \sum_{k=1}^{n^{c}}B_{kl}^{2} = 0$, which means that for all $l \in \CurlyBrackets{1,2,3,\dots, n^c}$, $\phi^{(2)}_{l} \in X^{(1)}$.
          Similarly, we can obtain for all $k = 1,2,3,\dots, n^{c}$, $\phi^{(1)}_{k} \in X^{(2)}$.
          Consequently, $P^{(1)} = P^{(2)}$ in the subspace sense.

    \item Note that $\tilde{\kappa}$ is a diagonal matrix, and we can derive
          \begin{align*}
            \Dist_{\tilde{\kappa}}\RoundBrackets{P^{(1)}, P^{(2)}} & = \sqrt{n^{c} - \norm{(T^{(1)}) ^\intercal \tilde{\kappa} T^{(2)}}_{F}^{2}}                                                           \\
                                                                   & = \sqrt{n^{c} - \norm{(T^{(2)}) ^\intercal \tilde{\kappa} T^{(1)}}_{F}^{2}}= \Dist_{\tilde{\kappa}} \RoundBrackets{P^{(2)}, P^{(1)}}.
          \end{align*}

    \item Denote $\xi^{(i)}_{k} = \sqrt{\tilde{\kappa}}\psi^{(i)}_{k}$ and
          $\CurlyBrackets{\psi^{(i)}_1, \psi^{(i)}_2, \dots, \psi^{(i)}_{n^c}}$ is an orthonormal basis of subspace $X^{(i)}$ with respectively to $\RoundBrackets{\cdot, \cdot}_{\tilde{\kappa}}$, thus $\CurlyBrackets{\xi^{(i)}_1, \xi^{(i)}_2, \dots, \xi^{(i)}_{n^c}}$ is an orthonormal basis of subspace $X^{(i)}$ with respectively to $l^{2}$ standard inner product.
          Take matrix $M^{(i)} =\sum_{k=1}^{n^c}\xi_{k}^{(i)}\RoundBrackets{\xi_{k}^{(i)}}^{\intercal}$, and we can obtain
          \[
            \begin{aligned}
               & \norm{M^{(i)}}^{2}_{F} = n^{c}, \quad \tr\RoundBrackets{{M^{(p)}}^{\intercal} M^{(q)}}= \sum_{k = 1}^{n^c}\sum_{l = 1}^{n^c}\RoundBrackets{\RoundBrackets{\xi_k^{(p)}}^\intercal \xi_l^{(q)}}^{2}, \\
               & \norm{M^{(p)} -M^{(q)}}^{2}_{F} = \norm{M^{(p)}}^{2}_{F} + \norm{M^{(q)}}^{2}_{F} - 2 \tr\RoundBrackets{(M^{(p)})^{\intercal} M^{(q)}}                                                             \\
               & \qquad = 2 \Dist^{2}_{\tilde{\kappa}}\RoundBrackets{P^{(p)}, P^{(q)}}.
            \end{aligned}
          \]
          By triangle inequality for Frobenius norm, it holds that
          \[
            \begin{aligned}
              \norm{M^{(1)} - M^{(3)}}_{F} \leq \norm{M^{(1)} - M^{(2)}}_{F} + \norm{M^{(2)} - M^{(3)}}_{F}.
            \end{aligned}
          \]
          Therefore, we can prove the triangle inequality for the distance in the loss function, i.e.,
          \[
            \Dist_{\tilde{\kappa}}\RoundBrackets{P^{(1)}, P^{(2)}}+ \Dist_{\tilde{\kappa}}\RoundBrackets
            {P^{(2)}, P^{(3)}}\geq \Dist_{\tilde{\kappa}}\RoundBrackets{P^{(1)}, P^{(3)}}.
          \]
  \end{itemize}
\end{proof}

Before we proceed to demonstrate the suitability of our proposed loss function for two-grid preconditioner, it is essential to present some theoretical foundation.
To define the performance of the two-grid method, we provide the error operator as follows.

\begin{definition}
  \label{lemma2} The error operator for two-grid method is given by
  \[
    E = I - B^{-1}A = \RoundBrackets*{I - RA}\RoundBrackets*{I - PA_{c}^{-1}P^{\intercal} A}
    \RoundBrackets*{I - AR},
  \]
  where $B^{-1}$ is the two-grid preconditioner.
\end{definition}

Next, we need to analyze the convergence rate of this error operator, which quantifies the efficiency of the two-grid method.
The proof can be found in \cite{Xu2017}.

\begin{lemma}\label{lemma3}
  The convergence rate of an exact two-grid method is given by
  \[
    \norm{E}_{A}^{2}= 1 - \frac{1}{K\RoundBrackets{P}},
  \]
  where
  \[
    K\RoundBrackets{P}= \max_{\substack{v \in \mathbb{R}^{N}}}\min_{\substack{v_c \in \mathbb{R}^{N^c}}}
    \frac{\norm{v - Pv_c}^{2}_{R^{-1}}}{\norm{v}^{2}_{A}}.
  \]
\end{lemma}

Building upon the theoretical foundation outlined above, we now establish correlation between the proposed loss function and the convergence rate, which aims to validate the efficacy of the loss function we have introduced.

\begin{theorem}
  \label{thm2} Let $E^{(1)}$ and $E^{(2)}$ be the error operators on the prolongation operators $P^{(1)}$ and $P^{(2)}$, respectively.
  Assume that $P^{(i)} = \Diag{\CurlyBrackets{P_1^{(i)}, P_2^{(i)}, \dots, P_n^{(i)}}}$ where $P_{j}^{(i)}$ is the block of $P^{(i)}$, $P^{(i)} \in \mathbb{R}^{N \times N^c}$ with $\ker\RoundBrackets{A}\subset \im\RoundBrackets{P^{(i)}}$. Moreover, suppose there exist positive constants $C_{\mathup{a}}$, $C_{\mathup{l}}$ and $C_{\mathup{u}}$, such that $\norm{\cdot}^{2}_{R^{-1}}\leq C_{\mathup{a}}\norm{\cdot}^{2}_{A}$ and $C_{\mathup{l}}\norm{\cdot}^{2}_{\tilde{\kappa}}\leq \norm{\cdot}^{2}_{R^{-1}}\leq C_{\mathup{u}}\norm{\cdot}^{2}_{{\tilde{\kappa}}}$.
  It holds that
  \[
    \abs{\norm{E^{(1)}}_{A}^{2} - \norm{E^{(2)}}_{A}^{2}}\leq C\RoundBrackets*{C_{\mathup{a}}, C_{\mathup{u}}/C_{\mathup{l}}}
    \max_{j}\CurlyBrackets*{\Dist_{\tilde{\kappa}_j}^2\RoundBrackets*{P_{j}^{(1)}, P_{j}^{(2)}}}.
  \]
\end{theorem}
\begin{proof}
  By \cref{lemma3}, we can see that
  \[
    \abs{\norm{E^{(1)}}_{A}^{2} - \norm{E^{(2)}}_{A}^{2}}= \abs{\frac{K(P^{(2)}) - K(P^{(1)})}{K(P^{(1)}) K(P^{(2)})}}
    \leq \abs{K(P^{(2)})-K(P^{(1)})},
  \]
  where we utilize the fact that $K(P^{(1)}) \geq  1$ and $K(P^{(2)}) \geq 1$.
  We are left to bound $\abs{K(P^{(2)}) - K(P^{(1)})}$.
  To improve the readability of proof, we will omit the index $j$ whenever it does not lead to ambiguity.

  \paragraph{Step1}\label{Step1}
  We first introduce the column space $X^{(i)}$ of the prolongation operator $P^{(i)}$ as follows,
  \[
    \begin{aligned}
      P^{(i)} = \SquareBrackets*{p_1^{(i)}, p_2^{(i)}, \dots, p_{N^c}^{(i)}}, \quad X^{(i)} = \Span \CurlyBrackets*{p_1^{(i)}, p_2^{(i)}, \dots, p_{N^c}^{(i)}}.
    \end{aligned}
  \]
  Under the assumption $N^{c} < N$ and $\ker\RoundBrackets{A}\subset \im\RoundBrackets{P^{(i)}}$, it holds that
  \[
    K\RoundBrackets{P^{(i)}}= \max_{\norm{v}_{A} = 1} \min_{v_{c}\in X^{(i)}}\norm{v - v_{c}}^{2}_{R^{-1}}.
  \]
  Assume that $\{\phi^{(i)}_{1}, \phi^{(i)}_{2}, \dots ,\phi^{(i)}_{N^{c}}\}$ is the orthonormal basis of subspace $X^{(i)}$ but with respectively to the inner product $\RoundBrackets{\cdot, \cdot}_{R^{-1}}$.
  Then for all $v \in \mathbb{R}^{N}$ with $\norm{v}_{A}= 1$,
  \[
    v^*_c \coloneqq \mathop{\arg\min}_{v_{c}\in X^{(i)}}\norm{v - v_{c}}^{2}_{R^{-1}} = \sum_{p=1}^{N^c}(v, \phi^{(i)}_{p})_{R^{-1}}\phi^{(i)}_{p},
  \]
  we can obtain $\norm{v - v_{c}^*}^{2}_{R^{-1}}= \norm{v}^{2}_{R^{-1}}- \sum_{p=1}^{N^c}(v, \phi^{(i)}_{p})^{2}_{R^{-1}}$.
  Under the assumption $\norm{\cdot}^{2}_{R^{-1}}\leq C_{\mathup{a}}\norm{\cdot}^{2}_{A}$, it yields $\norm{v_c^*}^{2}_{R^{-1}}\leq \norm{v}^{2}_{R^{-1}}\leq C_{\mathup{a}}\norm{v}^{2}_{A}\leq C_{\mathup{a}}$.

  \paragraph{Step2}
  \label{Step2}
  By scaling $K\RoundBrackets{P^{(1)}}$, we can obtain
  \[
    \begin{aligned}
      K\RoundBrackets{P^{(1)}} & = \max_{\norm{v}_{A} = 1}\min_{v_{c_1}\in X^{(1)}}\norm{v - v_{c_1}}^{2}_{R^{-1}}                                                                                                           \\
                               & \leq \max_{ \norm{v}_{A} = 1}\min_{v_{c_1}\in X^{(1)}}\RoundBrackets*{\norm{v - v_{c_2}}^2_{R^{-1}} + \norm{v_{c_2} - v_{c_1}}^2_{R^{-1}}}                                                  \\
                               & \leq \max_{\norm{v}_{A} = 1}\RoundBrackets*{\norm {v - v_{c_2}}^2_{R^{-1}} + \min_{v_{c_1}\in X^{(1)}} \norm{v_{c_2} - v_{c_1}}^2_{R^{-1}}}                                                 \\
                               & \leq \underbrace{\max_{ \norm{v}_{A} = 1} \norm{v - v_{c_2}}^2_{R^{-1}} }_{K_1}+ \underbrace{ \max_{\norm{v}_{A} = 1} \min_{v_{c_1}\in X^{(1)}} \norm{v_{c_2} - v_{c_1}}^2_{R^{-1}}}_{K_2}.
    \end{aligned}
  \]
  Let $v_{c_2} = \mathop{\arg\min}_{v_{c_2}\in X^{(2)}}\norm{v - v_{c}}^{2}_{R^{-1}} = \sum_{p=1}^{N^c}\RoundBrackets{v, \phi^{2}_p}_{R^{-1}}\phi^{2}_{p}$, and take it into $K_{1}$ and $K_{2}$.
  For $K_{1}$,
  \[
    \max_{\norm{v}_{A} = 1}\norm{v - \sum_{p=1}^{N^c} \RoundBrackets{v, \phi^{2}_p}_{R^{-1}} \phi^{2}_p}^{2}_{R^{-1}}= \max_{\norm {v}_{A} = 1}\min_{v_{c}\in X^{(2)}}\norm{v - v_{c}}^{2}_{R^{-1}}= K(P^{(2)}).
  \]
  For $K_{2}$, recalling that $\norm{v_{c_2}}_{R^{-1}}^2 \leq C_{\mathup{a}}$, we have
  \begin{align*}
     & \quad \max_{\norm {v}_{A} = 1}\min_{v_{c_1}\in X^{(1)}}\norm{\sum_{p=1}^{N^c} (v, \phi^{2}_p)_{R^{-1}} \phi^{2}_p - v_{c_1}}^{2}_{R^{-1}} \\
     & \leq  \max_{v \in X^{(2)}, \norm{v}_{R^{-1}} \leq \sqrt{C_\mathup{a}}}\min_{v_{c_1}\in X^{(1)}}\norm{v - v_{c_1}}^{2}_{R^{-1}}            \\
     & = C_{\mathup{a}}\max_{v \in X^{(2)}, \norm{v }_{R^{-1}} \leq 1 }\min_{v_{c_1}\in X^{(1)}}\norm{v - v_{c_1}}^{2}_{R^{-1}}.
  \end{align*}
  By the assumption $C_{\mathup{l}}\norm{\cdot}^{2}_{\tilde{\kappa}}\leq \norm{\cdot}^{2}_{R^{-1}}\leq C_{\mathup{u}}\norm{\cdot}^{2}_{{\tilde{\kappa}}}$,
  \[
    \begin{aligned}
      K\RoundBrackets{P^{(1)}} & \leq K\RoundBrackets{P^{(2)}}+ C_{\mathup{a}} \max_{v \in X^{(2)} , \norm{v}_{R^{-1}} \leq 1 }\min_{v_{c_1}\in X^{(1)}}\norm{v - v_{c_1}}^{2}_{R^{-1}}                                                        \\
                               & \leq K\RoundBrackets{P^{(2)}}+ \frac{C_{\mathup{a}} C_{\mathup{u}}}{C_{\mathup{l}}} \max_{v \in X^{(2)} , \norm{v}_{\tilde{\kappa}} \leq 1 }\min_{v_{c_1}\in X^{(1)}}\norm{v - v_{c_1}}^{2}_{\tilde{\kappa}}.
    \end{aligned}
  \]
  Due to the symmetry of $P^{(1)}$, $P^{(2)}$, we can also derive
  \[
    K\RoundBrackets{P^{(2)}}\leq K\RoundBrackets{P^{(1)}}+ \max_{v \in X^{(1)}, \norm{v}_{\tilde{\kappa}} \leq 1 }
    \min_{v_{c_2}\in X^{(2)}}\frac{C_{\mathup{a}}C_{\mathup{u}}}{C_{\mathup{l}}}
    \norm{v - v_{c_2}}^{2}_{\tilde{\kappa}}.
  \]

  \paragraph{Step3} \label{Step3}
  Take $\CurlyBrackets{\psi^{(i)}_1, \psi^{(i)}_2, \dots, \psi^{(i)}_{N^{c}}}$ as an orthonormal basis of subspace $X^{(i)}$ with respectively to $\RoundBrackets{\cdot, \cdot}_{\tilde{\kappa}}$.
  For all $v \in X^{(2)}$, introducing the expansion $v = \sum_{p=1}^{N^{c}}a_{p} \psi^{(2)}_{p}$ with $\norm{v}^{2}_{\tilde{\kappa}}= \sum_{p=1}^{N^{c}}a_{p}^{2}$, we can show that
  \[
    \max_{v \in X^{(2)}, \norm{v}_{\tilde{\kappa}} \leq 1 }\min_{v_{c_1}\in X^{(1)}}\norm{v - v_{c_1}}^{2}_{\tilde{\kappa}} = \max_{\sum_p a_p ^2 = 1}\sum_{p=1}^{N^{c}}a_{p}^{2} \RoundBrackets*{1 - \sum_{q=1}^{N^{c}} \RoundBrackets*{\psi^{(2)}_p, \psi^{(1)}_q}_{\tilde{\kappa}}^2}.
  \]
  Utilizing the non-overlapping properties of coarse elements, we can divide $a_{p}$, $\psi^{(2)}_p$ and $\psi^{(1)}_q$ into $a_{j,k}$, $\psi^{(2)}_{j,k}$ and $\psi^{(1)}_{j,l}$ respectively, where $j = 1, 2, \dots, n$ and $k = 1, 2, \dots, n^{c}$, $l = 1, 2, \dots, n^{c}$.
  Then, we can derive
  \begin{align*}
     & \quad \max_{\sum_p a_p ^2 = 1}\sum_{p=1}^{N^{c}}a_{p}^{2} \RoundBrackets*{1 - \sum_{q=1}^{N^{c}} \RoundBrackets*{\psi^{\RoundBrackets{2}}_p, \psi^{\RoundBrackets{1}}_q}_{\tilde{\kappa}}^2}                                                                                                          \\
     & = \max_{\sum_{j,k} a_{j,k} ^2 = 1}\sum_{j=1}^{n}\sum_{k=1}^{n^c}a_{j,k}^{2} \RoundBrackets*{1 - \sum_{l=1}^{n^{c}} \RoundBrackets*{\psi^{(2)}_{j,k}, \psi^{(1)}_{j,l}}_{\tilde{\kappa}_j}^2}                                                                                                          \\
     & \leq \max_{j}\CurlyBrackets*{\max_{\sum_k a_{j,k} ^2 = 1}\sum_{k=1}^{n^c} a_{j,k} ^2 \RoundBrackets*{1 - \sum_{l=1}^{n^{c}} \RoundBrackets*{\psi^{(2)}_{j,k}, \psi^{(1)}_{j,l}}_{\tilde{\kappa}_j}^2}} \leq \max_{j}\CurlyBrackets*{\Dist_{\tilde{\kappa}_j}^2\RoundBrackets*{P^{(1)}_j, P^{(2)}_j}}.
  \end{align*}
  Follow a similar procedure, we can also obtain
  \[
    \max_{v \in X_1,\norm {v}_{\tilde{\kappa}} \leq 1 }\min_{v_{c_2}\in X_2}\norm{v - v_{c_2}}
    ^{2}_{\tilde{\kappa}}\leq \max_{j}\CurlyBrackets*{\Dist_{\tilde{\kappa}_j}^2\RoundBrackets*{P^{(1)}_j, P^{(2)}_j}}.
  \]
  Combining all results above, we hence complete the proof.
\end{proof}

\subsection{Data augmentation}
In deep learning, the challenge of limited training data is pervasive.
To address this, we can artificially expand the training dataset through a series of transformations known as data augmentation.
Depending on the specific problem at hand, we recommend employing symmetry transformation and Karhunen--Lo\`{e}ve expansion as methods for implementing data augmentation.
These techniques not only diversify the dataset but also enhance the model's ability to generalize from
limited data inputs.

\subsubsection{Symmetry transformation}
In the context of the LSP, we observe inherent symmetries when the computational domain is partitioned appropriately.
In this study, we set $h_{x} =h_{y}$ to ensure that symmetric variations such as row symmetry, column symmetry, principal diagonal symmetry, and auxiliary diagonal symmetry are effectively considered.

\begin{theorem}
  \label{symmetric invariants} By taking $\tilde{\kappa}= \kappa$ in \cref{eq:spepb}, the LSP satisfies the following symmetric invariant relation:
  Let $\mathcal{T}$ be one of the symmetric variations which contains row symmetry, column symmetry, principal diagonal symmetry and auxiliary diagonal symmetry;
  If the pair $(\lambda^{k}, \Phi_{k})$ is the solution of $\cref{eq:spepb}$ based on $\kappa$, $\RoundBrackets{\lambda^k, \mathcal{T} \Phi_{k}}$ is the solution of $\cref{eq:spepb}$ based on $\mathcal{T}\kappa$.
\end{theorem}
\begin{proof}
  The symmetric variations are one-to-one such that for the right-hand side,
  \begin{equation}\label{eq:spepbb RHS}
    \lambda^k \int_{K_j}\kappa \Phi_{k}w_{h}\di \bm{x} = \lambda^k \int_{K_j}\RoundBrackets*{\mathcal{T}\kappa}\RoundBrackets*{\mathcal{T}\Phi_{k}} \RoundBrackets*{\mathcal{T} w_h} \di \bm{x}, \quad \forall w_{h} \in W_{h}(K_{i}).
  \end{equation}
  For the left-hand side, in \cref{fig:subfigures}, we can observe that the symmetric variations $\mathcal{T}$ do not alter the neighbor pairs of the edge set.
  This implies that if $\RoundBrackets{\kappa_{e,-}, \kappa_{e,+}}$ is a neighbor pair of an edge $e \in \mathcal{E}^{h}(K_{j})$, then $\RoundBrackets{\mathcal{T}\kappa_{e,-}, \mathcal{T}\kappa_{e,+}}$ remains a neighbor pair of an edge $e \in \mathcal{E}^{h}\RoundBrackets{\mathcal{T}K_j}$.
  Therefore, we can see that
  \begin{equation} \label{eq:spepbb LHS}
    \sum_{e\in \mathcal{E}^h(K_j)}\kappa_{e}\DSquareBrackets{\Phi_{k}}_{e}\DSquareBrackets{w_h}_{e} = \sum_{e\in \mathcal{E}^h(\mathcal{T}K_j)}\RoundBrackets{\mathcal{T}\kappa}_{e}\DSquareBrackets{\mathcal{T}\Phi_{k}}_{e}\DSquareBrackets{\mathcal{T}w_h}_{e} , \quad \forall w_{h} \in W_{h}(K_{i}).
  \end{equation}
  Combine \cref{eq:spepbb RHS} and \cref{eq:spepbb LHS},
  \begin{equation} \label{eq:spepbfinal1}
    \sum_{e\in \mathcal{E}^h(K_j)}\RoundBrackets{\mathcal{T}\kappa}_{e}\DSquareBrackets{\mathcal{T}\Phi_{k}}_{e}\DSquareBrackets{\mathcal{T} w_h}_{e} = \lambda^{k} \int_{K_j}\RoundBrackets*{\mathcal{T}\kappa}\RoundBrackets*{\mathcal{T}\Phi_{k}}\RoundBrackets*{\mathcal{T} w_h} \di \bm{x}, \quad \forall w_{h} \in W_{h}(K_{j}),
  \end{equation}
  which also can be written as
  \begin{equation}\label{eq:spepbfinal2}
    \sum_{e\in \mathcal{E}^h(K_j)}\RoundBrackets{\mathcal{T}\kappa}_{e}\DSquareBrackets{\mathcal{T}\Phi_{k}}_{e}\DSquareBrackets{w_h}_{e} = \lambda^{k}\int_{K_j}\RoundBrackets*{\mathcal{T}\kappa}\RoundBrackets*{\mathcal{T}\Phi_{k}} w_{h}\di \bm{x}, \quad \forall w_{h} \in W_{h}(\mathcal{T}K_{j}).
  \end{equation}

  \begin{figure}[!ht]
    \centering
    \begin{subfigure}
      [b]{0.3\textwidth}
      \centering
      \begin{tikzpicture}[scale=0.5]
        \draw[step=1cm] (0,0) grid (4,4);
        \fill[red, fill opacity=0.5] (1,3) rectangle (2,4);
        \fill[orange, fill opacity=0.5] (1,2) rectangle (2,3);
        \fill[yellow, fill opacity=0.5] (0,2) rectangle (1,3);
      \end{tikzpicture}
      \caption{original}
      \label{fig:subfig1}
    \end{subfigure}
    \begin{subfigure}
      [b]{0.3\textwidth}
      \centering
      \begin{tikzpicture}[scale=0.5]
        \draw[step=1cm] (0,0) grid (4,4);
        \fill[red, fill opacity=0.5] (1,1) rectangle (2,0);
        \fill[orange, fill opacity=0.5] (1,2) rectangle (2,1);
        \fill[yellow, fill opacity=0.5] (0,2) rectangle (1,1);
      \end{tikzpicture}
      \caption{row}
      \label{fig:subfig2}
    \end{subfigure}
    \begin{subfigure}
      [b]{0.3\textwidth}
      \centering
      \begin{tikzpicture}[scale=0.5]
        \draw[step=1cm] (0,0) grid (4,4);
        \fill[red, fill opacity=0.5] (3,3) rectangle (2,4);
        \fill[orange, fill opacity=0.5] (3,2) rectangle (2,3);
        \fill[yellow, fill opacity=0.5] (4,2) rectangle (3,3);
      \end{tikzpicture}
      \caption{column}
      \label{fig:subfig3}
    \end{subfigure}
    \begin{subfigure}
      [b]{0.3\textwidth}
      \centering
      \begin{tikzpicture}[scale=0.5]
        \draw[step=1cm] (0,0) grid (4,4);
        \fill[red, fill opacity=0.5] (0,2) rectangle (1,3);
        \fill[orange, fill opacity=0.5] (1,2) rectangle (2,3);
        \fill[yellow, fill opacity=0.5] (1,3) rectangle (2,4);
      \end{tikzpicture}
      \caption{principal diagonal}
      \label{fig:subfig4}
    \end{subfigure}
    \begin{subfigure}
      [b]{0.3\textwidth}
      \centering
      \begin{tikzpicture}[scale=0.5]
        \draw[step=1cm] (0,0) grid (4,4);
        \fill[red, fill opacity=0.5] (3,1) rectangle (4,2);
        \fill[orange, fill opacity=0.5] (2,1) rectangle (3,2);
        \fill[yellow, fill opacity=0.5] (2,0) rectangle (3,1);
      \end{tikzpicture}
      \caption{auxiliary diagonal}
      \label{fig:subfig5}
    \end{subfigure}
    \caption{Examples of symmetry transformation of $4 \times 4$ grid.}
    \label{fig:subfigures}
  \end{figure}
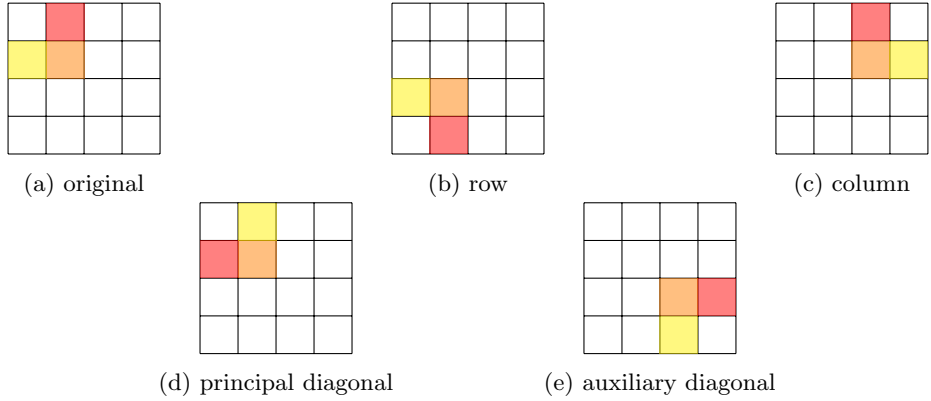
\end{proof}

\begin{remark}
  If $h_{x} \neq h_{y}$, \cref{symmetric invariants} still holds for row symmetry and column symmetry.
  It requires that $\mathcal{T}$ variation remains unchanged in the horizontal and vertical directions of the edges.
\end{remark}

Drawing on the insights gained, we have augmented the dataset by implementing transformations such as row symmetry, column symmetry, principal diagonal symmetry, and auxiliary diagonal symmetry to the data within the training set. Furthermore, the augmented data does not require separate label generation, significantly reducing the time required for training neural networks.

\subsubsection{Karhunen--Lo\`eve expansion}
The augmentation effect of symmetry transformation on data is limited.
For example, the mentioned symmetry methods can at most quadruple the original dataset.
When the data volume is still insufficient, one effective approach is to use the Karhunen--Lo\`{e}ve expansion to numerically approximate random field $\kappa\RoundBrackets{\bm{x};\omega}$ by utilizing discrete data reconstruction.
Karhunen--Lo\`{e}ve expansion seeks to represent a random field through a linear combination of a set of orthogonal basis functions.
To ensure positive permeability in the computational domain $\Omega$, we take $Z\RoundBrackets{\bm{x};\omega}= \log{\kappa\RoundBrackets{\bm{x};\omega}}$ and assume that $Z\RoundBrackets{\bm{x};\omega}$ is a Gaussian field.
To construct the Karhunen--Lo\`{e}ve expansion of $Z\RoundBrackets{\bm{x};\omega}$, we need to solve the following spectral problem
\[
  \int_{\Omega}\mathbb{C}_{Z}(\bm{x},\bm{y}) f(\bm{y}) \di \bm{y}
  = \mu f(\bm{x}),
\]
where $\mathbb{C}_{Z}(\bm{x},\bm{y})$ is the covariance kernel function and eigenpairs $\RoundBrackets{\mu^\alpha, f_\alpha}$ are obtained.
The number of total eigenpairs is equal to the degree of freedom of fine meshes, i.e., $N$ in our setting.
The Karhunen--Lo\`{e}ve expansion can be expressed as $Z(\bm{x}; \omega) = \mathbb{E}\SquareBrackets{Z}\RoundBrackets{\bm{x}}+\sum_{\alpha=1}^{N} \omega_{\alpha} \sqrt{\mu^\alpha}f_{\alpha}(\bm{x})$, where $\mathbb{E}\SquareBrackets{Z}\RoundBrackets{\bm{x}}$ is the mean value function, $\omega_{\alpha}$ are random variables with zero mean and unit variance.
For data augmentation, we take the $l$-truncated Karhunen--Lo\`{e}ve expansion
\begin{equation}\label{eq: KL truncated}
  Z(\bm{x}; \omega) \approx \mathbb{E}\SquareBrackets{Z}\RoundBrackets{\bm{x}}
  + \sum_{\alpha=1}^{l} \omega_{\alpha} \sqrt{\mu^\alpha}f_{\alpha}(\bm{x}),
\end{equation}
where $\CurlyBrackets{\mu^1, \mu^2, \dots, \mu^l}$ are the $l$ largest eigenvalues.
Consequently, the key of reconstruct random field is estimating the mean function and the covariance kernel function by some samples.
Denote $\CurlyBrackets*{Z_1\RoundBrackets{\bm{x}}, Z_2\RoundBrackets{\bm{x}}, \dots, Z_m\RoundBrackets{\bm{x}}}$ as $m$ samples, and we present the estimations for $\mathbb{E}[Z](\bm{x})$ and $\mathbb{C}_{Z}(\bm{x},\bm{y})$ as follows,
\begin{equation} \label{KL}
  \begin{aligned}
    \hat{\mathbb{E}}\SquareBrackets*{Z}\RoundBrackets{\bm{x}} & = \frac{1}{m}\sum_{s=1}^{m}Z_{s}\RoundBrackets*{{\bm{x}}},                                                                                                                                                                                                         \\
    \hat{\mathbb{C}}_{Z}\RoundBrackets*{\bm{x},\bm{y}}        & = \frac{1}{m}\sum_{s=1}^{m}\SquareBrackets*{\RoundBrackets*{Z_s\RoundBrackets*{{\bm{x}}} - \hat{\mathbb{E}}\SquareBrackets*{Z}\RoundBrackets*{\bm{x}}} \RoundBrackets*{Z_s\RoundBrackets*{{\bm{y}}} - \hat{\mathbb{E}}\SquareBrackets*{Z}\RoundBrackets*{\bm{y}}}}
  \end{aligned}
\end{equation}
which will be substituted into the Karhunen--Lo\`{e}ve expansion construction.

\section{Numerical experiments} \label{nume}
All neural network training is conducted on a single Nvidia A100-PCIE-40GB GPU, while the remaining iterative solving processes are performed on a desktop with an Intel Core i9-12900 CPU.
This arrangement is primarily due to the relatively small scale of the cases, which makes it difficult to fully utilize the GPU's high parallel computing power, meanwhile it also helps to reduce unnecessary energy consumption and improves resource efficiency.
The neural network training is implemented using PyTorch, while the iterative solver is based SciPy's sparse matrix routines.

In our simulations, we consider a heterogeneous computational domain $\Omega = \SquareBrackets{0,1}\times \SquareBrackets{0,1}$ with uncertainties.
Our square coarse grid $\mathcal{T}_{H}$ is of the size $16 \times 16$ and its refinement fine grid $\mathcal{T}_{h}$ is of $512 \times 512$ squares which means that each coarse element contains $32 \times 32$ fine elements.
Note that the smallest eigenvalue is always $0$ and the corresponding eigenvector is a constant function.
Therefore, we only need to generate $n^c-1$ local generalized multiscale basis by neural networks.
According to \cite{ye2024highly}, choosing $n^c=5$ is sufficient to achieve robustness while maintaining efficiency.
Hence, we set $n^c=5$ as the default value in our experiments.
The training ends after $30$ epochs.

\subsection{Prolongation operator learning experiments}
\subsubsection{Log-Gaussian random field}
Consider the logarithm of the random coefficients following a Gaussian random field, which has a mean function of zero and covariance kernel function is
\[
  \mathbb{C}_{Z}(\bm{x},\bm{y})  = \sigma^{2} \exp \left(-\sqrt{\frac{\left|x_1-y_1\right|^2}{\eta_1^2}+\frac{\left|x_2-y_2\right|^2}{\eta_2^2}}\right),
\]
where
$Z\RoundBrackets{\bm{x}; \omega}= \log{\kappa\RoundBrackets{\bm{x};\omega}}$, $\bm{x}= \RoundBrackets{x_1, x_2}$, $\bm{y}= \RoundBrackets{y_1, y_2}$, $\sigma^{2} = 2$, $\eta_{1}$ and $\eta_{2}$ correlation lengths in each direction of
space.

\begin{figure}[!ht]
  \centering
  \resizebox{0.9\textwidth}{!}{
    \input{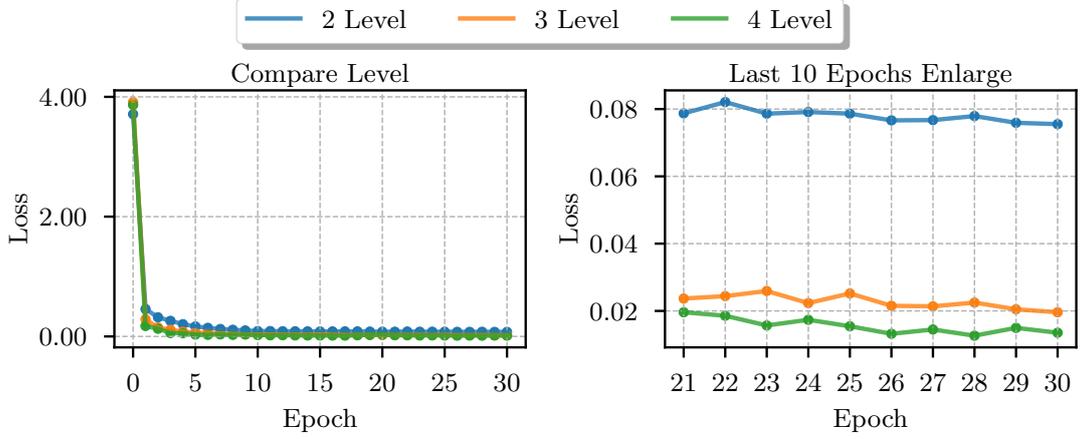}}
  \caption{Test loss for neural networks of varying depths on the same dataset is displayed in two parts: \textbf{left part} shows results across all $30$ epochs, and \textbf{right part} presents an enlarged view of the results for the last $10$ epochs.}
  \label{fig:complnn}
\end{figure}

We sample a total of $200$ instances of $\kappa(\bm{x};\omega)$ from the Log-Gaussian random field to construct our dataset.
This includes $180$ samples designated for the training dataset, with each $\kappa$ further subdivided into $16^2$ coarse elements, yielding a total of $46080$ examples.
The remaining $20$ samples are allocated to the test dataset, resulting in $5120$ examples.
As illustrated in \cref{fig:complnn}, deeper neural networks generally yield superior results.
Notably, the performance of the two-level U-Net is significantly distinct from that of the three-level and four-level U-Net.
Although the improvements between the three-level and four-level U-Net are relatively modest, the four-level U-Net is preferred due to its superior performance.

The number of training parameters and the training time required for different levels of U-Net architectures are summarized in \cref{table:unet_parameters}.
As anticipated, both the number of training parameters and the training time increase with the depth of the neural network.
However, the training time does not scale directly proportionally to the number of training parameters, which we attribute to the parallel computing capabilities and the design of the U-Net architecture.
Additionally, we highlight that by leveraging GPUs, all training are finished within $2$ minutes.

\begin{table}[!ht]
\footnotesize
  \centering
  \caption{The number of training parameters and the total training time of different levels of U-Net architectures.}
  \label{table:unet_parameters}
  \begin{tabular}{@{}ccc@{}}
    \toprule
    U-Net depth & Training parameters (M) & Training time (s) \\
    \midrule
    $2$-Level   & $0.11$                  & $83.10$           \\
    $3$-Level   & $0.52$                  & $97.49$           \\
    $4$-Level   & $2.14$                  & $106.75$          \\
    \bottomrule
  \end{tabular}
\end{table}

\begin{figure}[!ht]
  \centering
  \resizebox{0.9\textwidth}{!}{
    \input{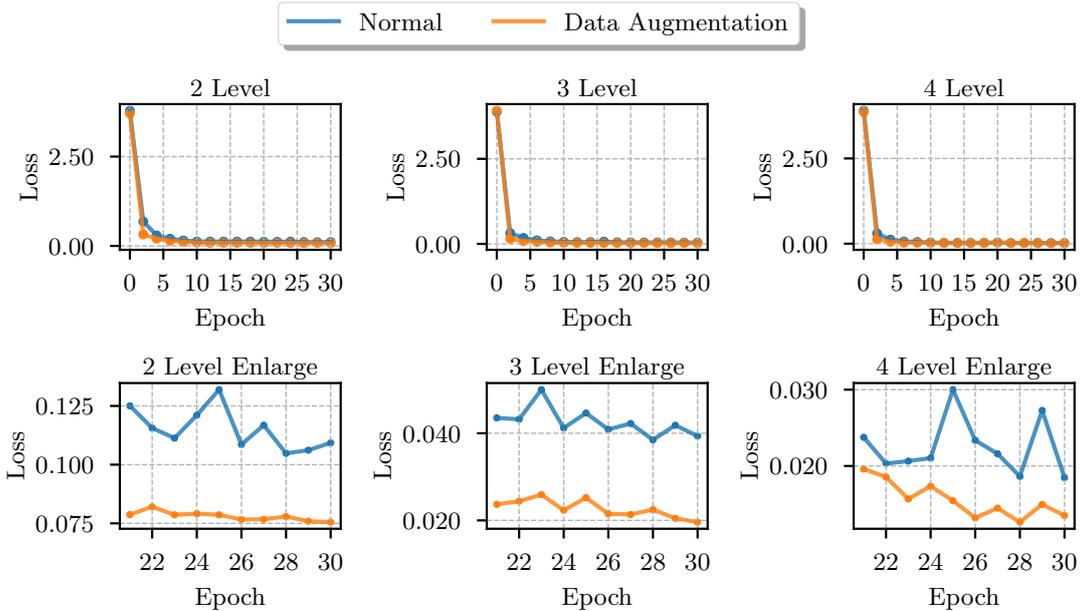}}
  \caption{
    Evaluate the impact of data augmentation on test loss for neural networks of varying depths by examining the differences before and after augmentation on the same test dataset.
    \textbf{Above part} displays the test loss trajectories for neural networks of different depths throughout all $30$
    training epochs.
    \textbf{Below part} offers an expanded view of the test loss for the final $10$ epochs.
  }
  \label{fig:DA}
\end{figure}

In situations where training data are scarce, we employ symmetry transformation as a data augmentation technique.
Specifically, we initially sample only $36$ instances of $\kappa(\bm{x};\omega)$ from Log-Gaussian random field to form our training dataset, while the test dataset remains unchanged.
We then substantially expand this dataset through symmetry transformation.
As indicated in \cref{fig:DA}, after applying symmetry transformation for data augmentation, there is a noticeable improvement in the performance of the two-level, three-level, and four-level U-Net architectures on the test set.
This enhancement substantiates the efficacy of symmetry transformation as a method for data augmentation.

\begin{figure}[!ht]
  \centering
  \resizebox{0.9\textwidth}{!}{
    \input{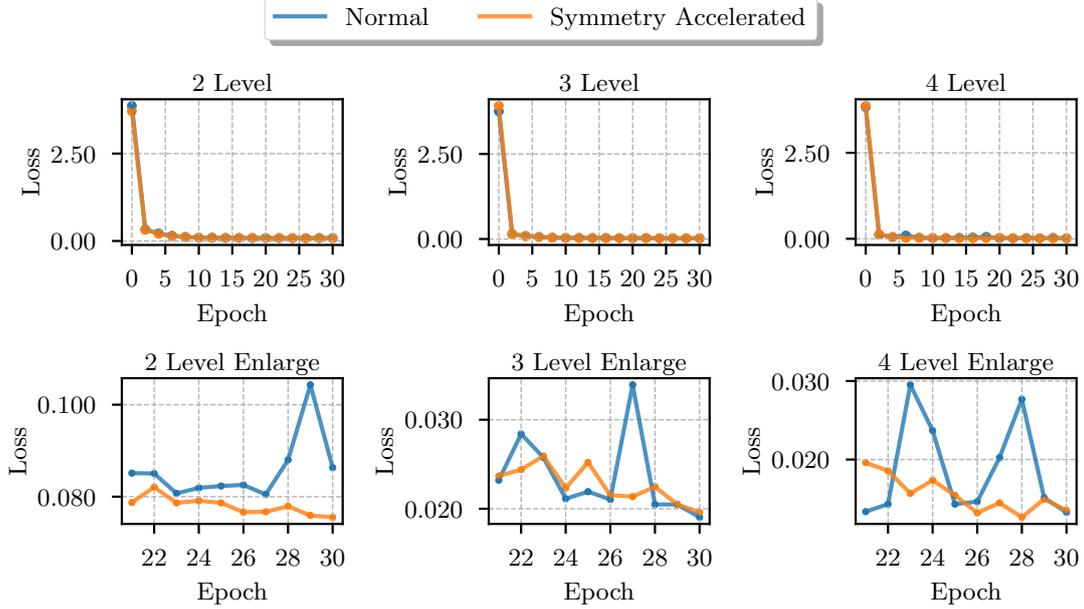}}
  \caption{Compare the test loss between training data all from random field and
    some data from augmentation for neural networks of different depths on the
    same test dataset. \textbf{Above part} displays the test loss trajectories for
    neural networks of different depths throughout all $30$ training epochs. \textbf{Below
      part} offers an expanded view of the test loss for the final $10$ epochs.}
  \label{fig:SA}
\end{figure}

As previously noted, due to inherent symmetries in the solutions of the LSP, we can avoid label generation for the augmented training data, significantly reducing the time required for generating neural networks models.
In our experiment setup, we compare two methods to generate the \emph{same amount} of training data for neural networks training: one involves sampling exclusively from Log-Gaussian random field and another one combines direct sampling from Log-Gaussian random field for $1/5$ of the data with the remaining $4/5$ generated through symmetry transformation augmentation.
The loss is calculated on the same test dataset.
According to the results shown in \cref{fig:SA}, using symmetry transformation to generate labels for the training set does not negatively impact the accuracy.
Furthermore, the training set augmented with symmetry transformations exhibits more stable performance during testing.

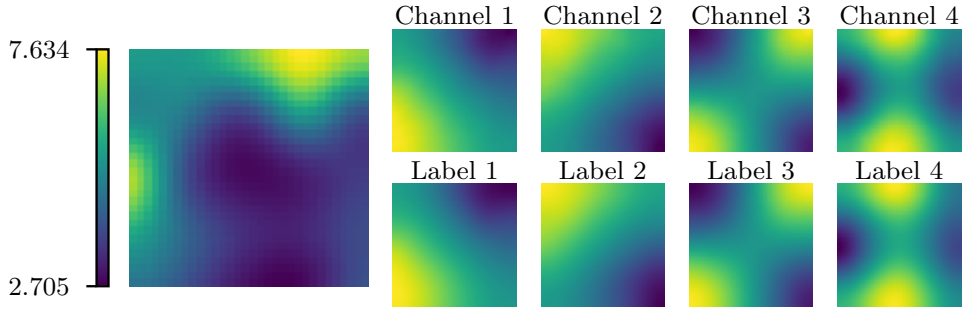
\begin{figure}[!ht]
  \centering
  \resizebox{0.9\textwidth}{!}{
    \input{RF.pgf}}
  \caption{A sample of an input coarse element from the test set is shown on the
    left.
    The labels are derived from the LSP, and the channels represent the
    projection of the corresponding labels onto the generalized multiscale
    subspace learned by the four-level U-Net.}
  \label{fig:R}
\end{figure}

Next, we opt for the four-level U-Net architecture, sampling $56$ instances of $\kappa\RoundBrackets{\bm{x};\omega}$, with $36$ samples designated for the training set, which is augmented using symmetry transformations, and $20$ samples are allocated to the test set.
A sample of the results from the four-level U-Net, depicted in \cref{fig:R}, reveals outcomes that closely resemble the generalized multiscale subspace derived from the LSP, illustrating the effectiveness of this architecture in capturing the essential characteristics of the data.

\subsubsection{Random disk inclusion coefficients}
In the subsequent experiments, we consider the presence of $15$ non-overlapping disks in the computational domain $\Omega$.
The locations and radii of these disks are randomly assigned.
The permeability value in disks, denoted $\kappa_{b}$, contrasts significantly with the rest of the domain where $\kappa$ is set to $1$, referred to as $\kappa_{r}$.
The coefficient $\kappa$ of the resolution $512 \times 512$ fine grid is determined by the following rules:
\begin{itemize}
  \item Elements entirely within or outside the disks adopt $\kappa_{b}$ and $\kappa_{r}$ respectively.
  \item For elements intersecting the boundaries of these disks, the effective permeability, $\kappa_{\text{element}}^{-1}$, is calculated as the harmonic mean of $\kappa_{b}^{-1}$ and $\kappa_{r}^{-1}$.
\end{itemize}
We test $\kappa_{b}$ chosen from $\CurlyBrackets{10^5, 10^4, 10^3, 10^2, 10^{-2}, 10^{-3}, 10^{-4}, 10^{-5}}$.
This setting introduces a significantly higher contrast in $\kappa$ values compared to those derived from Log-Gaussian random field, aiming to assess the model's capability to handle extreme contrast in permeability.

\begin{figure}[!ht]
  \centering
  \resizebox{0.9\textwidth}{!}{
    \input{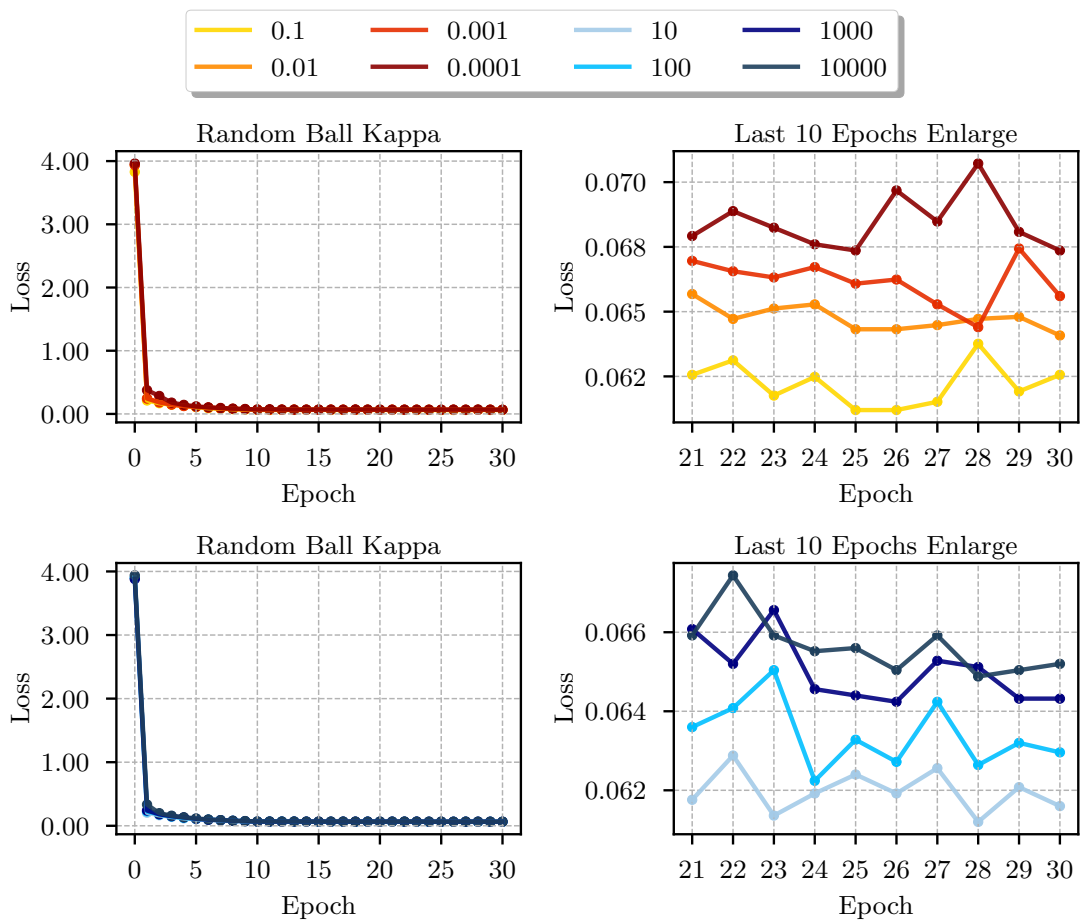}}
  \caption{Test loss for different $\kappa_{b}$ values using the four-level U-Net is displayed in two parts: \textbf{left part} shows results across all $30$ epochs, and \textbf{right part} presents an enlarged view of the results for the last $10$ epochs.}
  \label{fig:RBre}
\end{figure}

We sample $56$ instances of $\kappa\RoundBrackets{\bm{x};\omega}$ by randomly creating disks within the computational domain.
Of these, $36$ samples are designated for the training set and augmented through symmetry transformation, while the remaining $20$ samples are reserved for the test set.
The performance of the four-level U-Net across varying $\kappa_{b}$ values is depicted in \cref{fig:RBre}.
As expected, lower contrast cases yield better performance, but only minor differences are observed compared to higher contrast cases.

\begin{figure}[!ht]
  \centering
  \resizebox{0.9\textwidth}{!}{
    \input{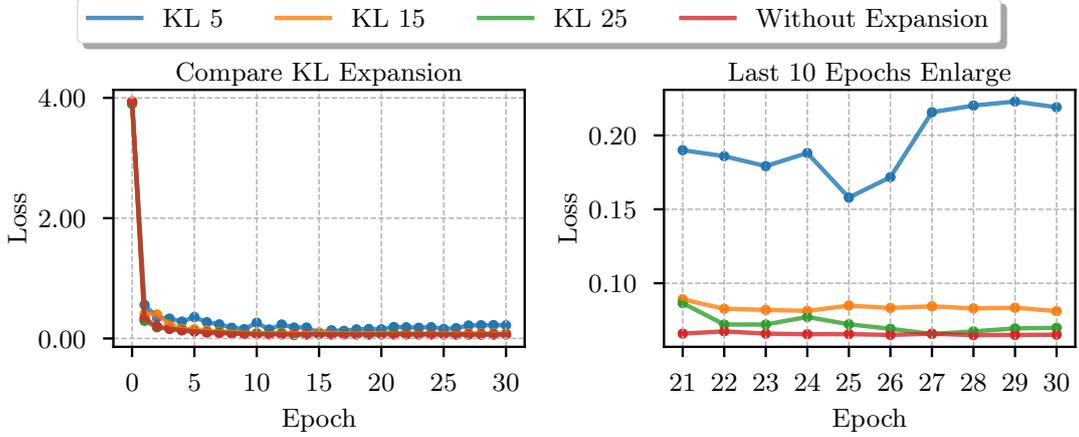}}
  \caption{The test loss for the four-level U-Net, trained with $\kappa_{b} = 10^{4}$, and data augmentation factor $M = 2$, across different truncation levels $K$ of Karhunen--Lo\`{e}ve expansion is displayed in two parts: \textbf{left part} shows results across all $30$ epochs, and \textbf{right part} presents an enlarged view of the results for the last $10$ epochs.}
  \label{fig:KLE}
\end{figure}

\begin{figure}[!ht]
  \centering
  \resizebox{0.9\textwidth}{!}{
    \input{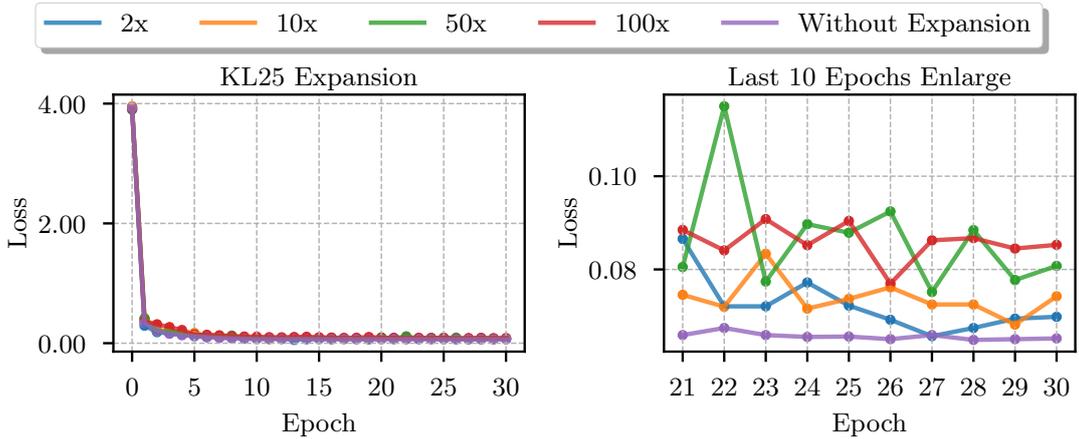}}
  \caption{The test loss for a four-level U-Net, trained with $\kappa_{b} = 10^{4}$, and truncation level $K = 25$, across varying data augmentation factors $M$ is displayed in two parts: \textbf{left part} shows results across all $30$ epochs, and \textbf{right part} presents an enlarged view of the results for the last $10$ epochs.}
  \label{fig:KLEM}
\end{figure}

We here consider apply Karhunen--Lo\`{e}ve expansion for data augmentation on this random coefficient setting.
Let $K$ represent the $K$-truncated Karhunen--Lo\`{e}ve expansion as in \cref{eq: KL truncated}, and let $M$ denote the multiplication factor required to augment the training dataset.
Specifically, we take $46080/M$ coarse elements from a training set generated by random disks, reconstruct the random field as \cref{KL} the $46080/M$ coarse elements and augment the training data to $46080$.
From \cref{fig:KLE}, we observe that increasing the truncation level $K$ of Karhunen--Lo\`{e}ve expansion enhances its
effectiveness in augmenting the dataset.
Specifically, when $K$ reaches $25$, the performance of the model closely approximates that of the training conducted without employing the Karhunen--Lo\`{e}ve expansion.
This finding suggests that a $K$-value of $25$ provides a sufficient approximation of the underlying random field, enabling substantial improvements in model training by capturing essential features within the data.
We also investigate the effect of varying the data augmentation factor $M$ on the model's performance.
In \cref{fig:KLEM}, the results for the four-level U-Net, which is trained using a dataset expanded by a factor of $2$ to $100$ through the application of a $25$-truncated Karhunen--Lo\`{e}ve expansion are reported.
Performance remains sufficiently effective at different values of $M$, although deterioration is observed when $M$ is large.
This approach highlights the robustness of the Karhunen--Lo\`{e}ve expansion in enhancing the training process under conditions of substantial data augmentation.

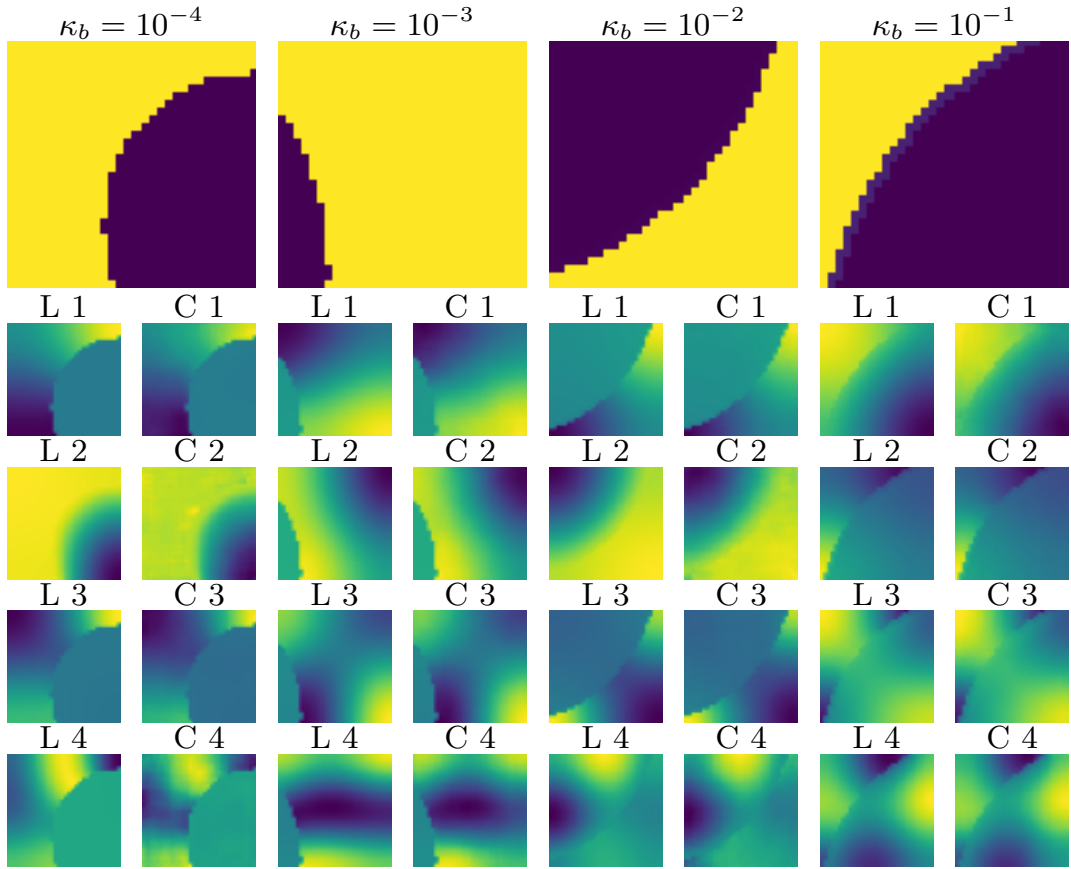
\begin{figure}[!ht]
  \centering
  \resizebox{0.9\textwidth}{!}{
    \input{_x.pgf}}
  \caption{Samples of input coarse elements ($\kappa_{b}$ values from $10^{-5}$ to $10^{-2}$) from test set shows in left part, L1 to L4 are derived from the LSP and C1 to C4 represent the projection of the corresponding labels onto the generalized multiscale subspace learned by the four-level U-Net.}
  \label{fig:x1}
\end{figure}

We also virtualize the output of neural networks in this random coefficient setting.
The dataset comprises $56$ instances of $\kappa\RoundBrackets{\bm{x};\omega}$, with $36$ samples designated for the training set, which is enhanced using symmetry transformations, and $20$ samples are allocated to the test set.
It presents samples of $\kappa_{b}$ values from $10^{-5}$ to $10^{-2}$ in \cref{fig:x1}.
These samples exhibit characteristics that are notably similar to the generalized multiscale subspace derived from the LSP.
Similarly, for $\kappa_{b}$ values between $10^{2}$ and $10^{5}$, the results are equally promising, and due to space constraints, are not presented here.

\subsection{Preconditioner comparison}
In this section, we examine the efficacy of a prolongation operator generated by the four-level U-Net in two-grid preconditioner and compare its performance to other preconditioners.
All types of preconditioners are integrated within the PCG iteration scheme, which is employed to solve equations with the no-flux boundary condition.
The source term $f$ is assigned values of $1$ in the top-left and bottom-right fine elements, $-1$ in the bottom-left and top-right fine elements, and $0$ elsewhere.
The PCG with a two-grid preconditioner is accompanied by a Block-Jacobi smoother, with the block size tailored to the coarse elements' dimensions.
The relative tolerance is set to $10^{-6}$.
The coefficient profiles are generated following the same procedure as in the previous section, and we adopt notations ``Gaussian RF'' and ``RB'' to represent the Log-Gaussian random field and random disk inclusion coefficient respectively.
Each experiment is repeated $10$ times, and the time records presented below are the averages of these repeated experiments.

We first test the two-grid preconditioner with the prolongation operator generated by the four-level U-Net.
Additionally, we compare the iteration counts (iteration time) with the prolongation operators generated by precisely solving LSPs with $1$, $3$ and $5$ bases.
As depicted in \cref{tab:comparison}, in terms of iteration counts, the prolongation operator produced by the four-level U-Net exhibits performance on par with that derived from $5$ basis functions in each LSP, and it significantly
outperforms the operators derived from either one or three basis functions in each LSP.
Importantly, we observe actual acceleration on solution time up to $50\%$ by using neural network techniques.

\begin{table}[!ht]
\footnotesize
  \centering
  \caption{Iteration counts and computing time for the two-grid preconditioner with the prolongation operators generated by solving LSPs precisely and the four-level U-Net, where ``1-LSP-MG'', ``3-LSP-MG'', ``5-LSP-MG'' and ``5-NN-MG'' represent the prolongation operators derived from LSPs with $1$, $3$, $5$ bases and the four-level U-Net respectively.}
  \label{tab:comparison}
  \begin{tabular}{@{}lcccc@{}}
    \toprule
    Coefficient profiles      & 1-LSP-MG        & 3-LSP-MG      & 5-LSP-MG               & 5-NN-MG                         \\
    \midrule
    Gaussian RF               & $67$ $(11.92)$  & $22$ $(7.49)$ & $\textbf{12}$ $(6.19)$ & $13$ ($\textbf{3.06}$)          \\
    RB $\kappa_{b} = 10^{2}$  & $67$ $(11.62)$  & $29$ $(8.85)$ & $\textbf{17}$ $(7.08)$ & $19$ ($\textbf{3.89}$)          \\
    RB $\kappa_{b} = 10^{-2}$ & $68$ $(11.70)$  & $24$ $(7.76)$ & $\textbf{16}$ $(7.13)$ & $17$ ($\textbf{3.68}$)          \\
    RB $\kappa_{b} = 10^{3}$  & $72$ $(12.88)$  & $29$ $(9.11)$ & $\textbf{19}$ $(7.40)$ & $21$ ($\textbf{4.33}$)          \\
    RB $\kappa_{b} = 10^{-3}$ & $71$ $(12.11)$  & $24$ $(8.03)$ & $\textbf{17}$ $(7.16)$ & $\textbf{17}$ ($\textbf{3.78}$) \\
    RB $\kappa_{b} = 10^{4}$  & $77$ $(13.41)$  & $30$ $(8.90)$ & $\textbf{19}$ $(7.74)$ & $21$ ($\textbf{4.39}$)          \\
    RB $\kappa_{b} = 10^{-4}$ & $82$ $(14.43)$  & $27$ $(8.19)$ & $\textbf{17}$ $(7.45)$ & $18$ ($\textbf{3.92}$)          \\
    RB $\kappa_{b} = 10^{5}$  & $136$ $(23.40)$ & $34$ $(9.76)$ & $\textbf{21}$ $(8.00)$ & $23$ ($\textbf{4.77}$)          \\
    RB $\kappa_{b} = 10^{-5}$ & $121$ $(21.06)$ & $29$ $(9.13)$ & $\textbf{20}$ $(7.85)$ & $23$ ($\textbf{4.82}$)          \\
    \bottomrule
  \end{tabular}
\end{table}

We then compare the iteration counts (iteration time) with two easily obtainable preconditioners on the Python platform, which are pyAMG, a Python implementation of algebraic multigrid methods \cite{pyamg}, and ILU, an incomplete LU decomposition provided in SciPy.
We specially take the ILU(10) preconditioner as a representative of the ILU family.
In \cref{tab:Pre_of_comparison}, we observe that, as the contrast of $\kappa$ increases, the two-grid preconditioner gradually outperforms the ILU(10) preconditioner and the AMG preconditioner in terms of iteration counts.
However, due to the time-consuming process of solving LSPs to generate the prolongation operator, it does not exhibit an advantage in terms of iteration time.
Nonetheless, by employing neural networks acceleration to generate the prolongation operator for the two-grid preconditioner, when the contrast of $\kappa$ reaches $10^{4}$, the accelerated two-grid preconditioner outperforms both the ILU(10) preconditioner and the pyAMG preconditioner in terms of both iteration counts and iteration time.
This result demonstrates the potential of proposed method in solving large linear systems with high-contrast coefficients.

\begin{table}[!ht]
\footnotesize
  \centering
  \caption{Iteration counts and computing time for ILU(10), pyAMG, the two-grid preconditioner by solving LSPs precisely, and two-grid preconditioner by the four-level U-Net.}
  \label{tab:Pre_of_comparison}
  \begin{tabular}{@{}lcccc@{}}
    \toprule
    Coefficient profiles      & ILU(10)         & pyAMG                           & 5-LSP-MG               & 5-NN-MG                \\
    \midrule
    Gaussian RF               & $17$ $(3.26)$   & $\textbf{6}$ $(\textbf{2.50})$  & $12$ $(6.19)$          & $13$ ($3.06$)          \\
    RB $\kappa_{b} = 10^{2}$  & $38$ $(3.89)$   & $\textbf{15}$ $(\textbf{2.78})$ & $17$ $(7.08)$          & $19$ ($3.89$)          \\
    RB $\kappa_{b} = 10^{-2}$ & $32$ $(3.77)$   & $\textbf{14}$ $(\textbf{2.74})$ & $16$ $(7.13)$          & $17$ ($3.68$)          \\
    RB $\kappa_{b} = 10^{3}$  & $107$ $(5.21)$  & $48$ $(\textbf{3.64})$          & $\textbf{19}$ $(7.40)$ & $21$ ($4.33$)          \\
    RB $\kappa_{b} = 10^{-3}$ & $123$ $(5.57)$  & $44$ $(\textbf{3.58})$          & $\textbf{17}$ $(7.16)$ & $\textbf{17}$ ($3.78$) \\
    RB $\kappa_{b} = 10^{4}$  & $276$ $(9.74)$  & $62$ $(4.51)$                   & $\textbf{19}$ $(7.74)$ & $21$ ($\textbf{4.39}$) \\
    RB $\kappa_{b} = 10^{-4}$ & $285$ $(9.26)$  & $59$ $(4.14)$                   & $\textbf{17}$ $(7.45)$ & $18$ ($\textbf{3.92}$) \\
    RB $\kappa_{b} = 10^{5}$  & $528$ $(13.64)$ & $101$ $(6.18)$                  & $\textbf{21}$ $(8.00)$ & $23$ ($\textbf{4.77}$) \\
    RB $\kappa_{b} = 10^{-5}$ & $527$ $(14.80)$ & $82$ $(5.79)$                   & $\textbf{20}$ $(7.85)$ & $23$ ($\textbf{4.82}$) \\
    \bottomrule
  \end{tabular}
\end{table}

In \cref{tab:ttime}, we present a comparative analysis of the average time consumption involved in obtaining the prolongation operator from the test dataset using the four-level U-Net versus the time taken to solve LSPs. 
The results indicate that the four-level U-Net architecture significantly reduces the time required to generate the prolongation operator.

\begin{table}[!ht]
\footnotesize
  \centering
  \caption{Comparison of time consumption for solving LSPs and NN across different datasets.}
  \label{tab:ttime}
  \begin{tabular}{@{}lcc@{}}
    \toprule
    Data Set                   & Time for Solving LSP (s) & Time for NN (s)   \\
    \midrule
    Log-Gaussian Random Fields & $4.2791$                 & $\textbf{0.8379}$ \\
    RB $\kappa_b = 1:10^{-5}$  & $4.2885$                 & $\textbf{0.8834}$ \\
    RB $\kappa_b = 1:10^{-4}$  & $4.2973$                 & $\textbf{0.8929}$ \\
    RB $\kappa_b = 1:10^{-3}$  & $4.1959$                 & $\textbf{0.8405}$ \\
    RB $\kappa_b = 1:10^{-2}$  & $4.2318$                 & $\textbf{0.8532}$ \\
    RB $\kappa_b = 1:10^2$     & $4.2302$                 & $\textbf{0.8381}$ \\
    RB $\kappa_b = 1:10^3$     & $4.1966$                 & $\textbf{0.8426}$ \\
    RB $\kappa_b = 1:10^4$     & $4.2201$                 & $\textbf{0.8348}$ \\
    RB $\kappa_b = 1:10^5$     & $4.2492$                 & $\textbf{0.8598}$ \\
    \midrule
    Average                    & $4.2709$                 & $\textbf{0.8537}$ \\
    \bottomrule
  \end{tabular}
\end{table}

\subsection{Generalization ability tests}
We first evaluate the generalization ability on different resolutions of the coarse element.
The neural networks we designed is trained by coarse elements of size $32 \times 32$, yet it remains applicable to smaller coarse elements.
To adapt these smaller elements to the required input size of the network, one can employ nearest neighbor interpolation.
The resulting outputs can then be refined to the target dimensions through mean downsampling.
We evaluate coarse elements of sizes of $16 \times 16$ and $8 \times 8$, which are extracted from the Gaussian RF coefficients.
The reports are presented in \cref{tab:Different Size}.
For smaller coarse elements, the two-grid preconditioner  derived from the pre-trained four-level U-Net model requires slightly more iterations to converge than that derived from precisely solved $5$ basis functions, while fewer iterations than those derived from $3$ or $1$ basis function(s) in each LSP.
While regarding iteration time, the pre-trained four-level U-Net model exhibits a notable acceleration.

\begin{table}[!ht]
\footnotesize
  \centering
  \caption{Iteration counts and computing time on different coarse element sizes for the two-grid preconditioner with the prolongation operators generated by solving LSPs precisely and the pre-trained four-level U-Net model.}
  \label{tab:Different Size}
  \begin{tabular}{@{}lcccc@{}}
    \toprule
    Coarse element sizes & 1-LSP-MG      & 3-LSP-MG      & 5-LSP-MG              & 5-NN-MG               \\
    \midrule
    $16\times 16$        & $38$ $(6.92)$ & $13$ $(6.49)$ & $\textbf{7}$ $(5.58)$ & $8$ ($\textbf{2.88}$) \\
    $8\times 8$          & $23$ $(5.96)$ & $9$ $(6.02)$  & $\textbf{5}$ $(5.59)$ & $7$ ($\textbf{2.61}$) \\
    \bottomrule
  \end{tabular}
\end{table}

We then evaluate the generalization ability of the neural network model on coefficients that fall outside the training set distribution. 
Specifically, we test the model trained on the Gaussian RF dataset using the RB test set. As shown in \cref{tab:Generalization_Ability_of_comparison }, although some performance degradation is observed, the results confirm a significant acceleration in computation time. This suggests that the neural network effectively replaces the nonlinear process of solving the LSP and does not heavily rely on $\kappa$ following a specific distribution. 
We believe this finding is highly promising and warrants further in-depth investigation, both theoretically and in practical applications.

\begin{table}[!ht]
\footnotesize
  \centering
  \caption{Iteration counts and computing time for the two-grid preconditioner with the prolongation operators generated by solving LSPs precisely and the pre-trained four-level U-Net, where ``5-NN(RF)-MG'' represents the prolongation operators derived from the four-level U-Net trained with Log-Gaussian random dataset.}
  \label{tab:Generalization_Ability_of_comparison}
  \begin{tabular}{@{}lcccc@{}}
    \toprule
    Coefficient profiles      & 1-LSP-MG        & 3-LSP-MG      & 5-LSP-MG               & 5-NN(RF)-MG            \\
    \midrule
    RB $\kappa_{b} = 10^{2}$  & $67$ $(11.62)$  & $29$ $(8.85)$ & $\textbf{17}$ $(7.08)$ & $20$ ($\textbf{3.97}$) \\
    RB $\kappa_{b} = 10^{-2}$ & $68$ $(11.70)$  & $24$ $(7.76)$ & $\textbf{16}$ $(7.13)$ & $18$ ($\textbf{3.68}$) \\
    RB $\kappa_{b} = 10^{3}$  & $72$ $(12.88)$  & $29$ $(9.11)$ & $\textbf{19}$ $(7.40)$ & $22$ ($\textbf{4.34}$) \\
    RB $\kappa_{b} = 10^{-3}$ & $71$ $(12.11)$  & $24$ $(8.03)$ & $\textbf{17}$ $(7.16)$ & $20$ ($\textbf{4.05}$) \\
    RB $\kappa_{b} = 10^{4}$  & $77$ $(13.41)$  & $30$ $(8.90)$ & $\textbf{19}$ $(7.74)$ & $25$ ($\textbf{5.08}$) \\
    RB $\kappa_{b} = 10^{-4}$ & $82$ $(14.43)$  & $27$ $(8.19)$ & $\textbf{17}$ $(7.45)$ & $22$ ($\textbf{4.52}$) \\
    RB $\kappa_{b} = 10^{5}$  & $136$ $(23.40)$ & $34$ $(9.76)$ & $\textbf{21}$ $(8.00)$ & $28$ ($\textbf{5.47}$) \\
    RB $\kappa_{b} = 10^{-5}$ & $121$ $(21.06)$ & $29$ $(9.13)$ & $\textbf{20}$ $(7.85)$ & $26$ ($\textbf{5.18}$) \\
    \bottomrule
  \end{tabular}
\end{table}

\section{Conclusion and future work} \label{con}
This research addressed the computational challenges of linear systems from Darcy flow discretization under random permeability.
Traditional generalized multiscale two-grid methods require solving a local spectral problem for each coarse element, which was time-consuming, while its robustness on the high-contrast coefficient are guaranteed.
We proposed a deep learning based surrogate model using a four-level U-Net to generate blocks of the generalized multiscale prolongation operator directly from coarse element coefficients.
Our loss function was based on the distance between subspaces defined by these coefficients.
We achieved a label generation time reduction of up to 5 folds through symmetric transformations and enhanced data efficiency via Karhunen--Lo\`{e}ve expansion.
Numerical experiments showed that our method reduced prolongation operator generation time by up to a factor of $5$ while preserving the efficiency of the two-grid preconditioner.
For the generalization ability of neural networks, our trained model is applicable to different sizes of coarse element.
Besides, the model trained on low-contrast datasets remains significantly effective when applied to high-contrast and off-distribution cases.
Recognizing the mutual influence of the smoother and the prolongation operator is essential for enhancing the effectiveness of the two-grid preconditioner.
Our future work will focus on developing an optimal smoother by deep learning to further enhance preconditioning.

\section*{Acknowledgments}
During the edition of this work, the authors utilized ChatGPT to enhance readability and language quality. 
Following the use of this tool, the authors thoroughly reviewed and as needed and accept full responsibility for the final publication.

\input{main.bbl}\end{document}

%% file: RF.pgf
%% Creator: Matplotlib, PGF backend
%%
%% To include the figure in your LaTeX document, write
%%   \input{<filename>.pgf}
%%
%% Make sure the required packages are loaded in your preamble
%%   \usepackage{pgf}
%%
%% Also ensure that all the required font packages are loaded; for instance,
%% the lmodern package is sometimes necessary when using math font.
%%   \usepackage{lmodern}
%%
%% Figures using additional raster images can only be included by \input if
%% they are in the same directory as the main LaTeX file. For loading figures
%% from other directories you can use the `import` package
%%   \usepackage{import}
%%
%% and then include the figures with
%%   \import{<path to file>}{<filename>.pgf}
%%
%% Matplotlib used the following preamble
%%   \def\mathdefault#1{#1}
%%   \everymath=\expandafter{\the\everymath\displaystyle}
%%   \usepackage{bm}%
%%   \usepackage{mismath}
%%   \newcommand{\SSSText}[1]{{\scriptscriptstyle \mathup{#1}}}%
%%   \renewcommand{\mathdefault}[1][]{}%
%%   \ifdefined\pdftexversion\else  % non-pdftex case.
%%     \usepackage{fontspec}
%%   \fi
%%   \makeatletter\@ifpackageloaded{underscore}{}{\usepackage[strings]{underscore}}\makeatother
%%
\begingroup%
\makeatletter%
\begin{pgfpicture}%
\pgfpathrectangle{\pgfpointorigin}{\pgfqpoint{5.300000in}{1.766667in}}%
\pgfusepath{use as bounding box, clip}%
\begin{pgfscope}%
\pgfsetbuttcap%
\pgfsetmiterjoin%
\definecolor{currentfill}{rgb}{1.000000,1.000000,1.000000}%
\pgfsetfillcolor{currentfill}%
\pgfsetlinewidth{0.000000pt}%
\definecolor{currentstroke}{rgb}{1.000000,1.000000,1.000000}%
\pgfsetstrokecolor{currentstroke}%
\pgfsetdash{}{0pt}%
\pgfpathmoveto{\pgfqpoint{0.000000in}{0.000000in}}%
\pgfpathlineto{\pgfqpoint{5.300000in}{0.000000in}}%
\pgfpathlineto{\pgfqpoint{5.300000in}{1.766667in}}%
\pgfpathlineto{\pgfqpoint{0.000000in}{1.766667in}}%
\pgfpathlineto{\pgfqpoint{0.000000in}{0.000000in}}%
\pgfpathclose%
\pgfusepath{fill}%
\end{pgfscope}%
\begin{pgfscope}%
\pgfsetbuttcap%
\pgfsetmiterjoin%
\definecolor{currentfill}{rgb}{1.000000,1.000000,1.000000}%
\pgfsetfillcolor{currentfill}%
\pgfsetlinewidth{0.000000pt}%
\definecolor{currentstroke}{rgb}{0.000000,0.000000,0.000000}%
\pgfsetstrokecolor{currentstroke}%
\pgfsetstrokeopacity{0.000000}%
\pgfsetdash{}{0pt}%
\pgfpathmoveto{\pgfqpoint{0.819211in}{0.307391in}}%
\pgfpathlineto{\pgfqpoint{1.953429in}{0.307391in}}%
\pgfpathlineto{\pgfqpoint{1.953429in}{1.441609in}}%
\pgfpathlineto{\pgfqpoint{0.819211in}{1.441609in}}%
\pgfpathlineto{\pgfqpoint{0.819211in}{0.307391in}}%
\pgfpathclose%
\pgfusepath{fill}%
\end{pgfscope}%
\begin{pgfscope}%
\pgfpathrectangle{\pgfqpoint{0.819211in}{0.307391in}}{\pgfqpoint{1.134218in}{1.134218in}}%
\pgfusepath{clip}%
\pgfsys@transformshift{0.819211in}{0.307391in}%
\pgftext[left,bottom]{\includegraphics[interpolate=true,width=1.140000in,height=1.140000in]{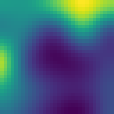}}%
\end{pgfscope}%
\begin{pgfscope}%
\pgfsetbuttcap%
\pgfsetmiterjoin%
\definecolor{currentfill}{rgb}{1.000000,1.000000,1.000000}%
\pgfsetfillcolor{currentfill}%
\pgfsetlinewidth{0.000000pt}%
\definecolor{currentstroke}{rgb}{0.000000,0.000000,0.000000}%
\pgfsetstrokecolor{currentstroke}%
\pgfsetstrokeopacity{0.000000}%
\pgfsetdash{}{0pt}%
\pgfpathmoveto{\pgfqpoint{0.662500in}{0.307391in}}%
\pgfpathlineto{\pgfqpoint{0.719211in}{0.307391in}}%
\pgfpathlineto{\pgfqpoint{0.719211in}{1.441609in}}%
\pgfpathlineto{\pgfqpoint{0.662500in}{1.441609in}}%
\pgfpathlineto{\pgfqpoint{0.662500in}{0.307391in}}%
\pgfpathclose%
\pgfusepath{fill}%
\end{pgfscope}%
\begin{pgfscope}%
\pgfsys@transformshift{0.660000in}{0.316667in}%
\pgftext[left,bottom]{\includegraphics[interpolate=true,width=0.060000in,height=1.130000in]{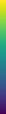}}%
\end{pgfscope}%
\begin{pgfscope}%
\pgfsetbuttcap%
\pgfsetroundjoin%
\definecolor{currentfill}{rgb}{0.000000,0.000000,0.000000}%
\pgfsetfillcolor{currentfill}%
\pgfsetlinewidth{0.803000pt}%
\definecolor{currentstroke}{rgb}{0.000000,0.000000,0.000000}%
\pgfsetstrokecolor{currentstroke}%
\pgfsetdash{}{0pt}%
\pgfsys@defobject{currentmarker}{\pgfqpoint{-0.048611in}{0.000000in}}{\pgfqpoint{-0.000000in}{0.000000in}}{%
\pgfpathmoveto{\pgfqpoint{-0.000000in}{0.000000in}}%
\pgfpathlineto{\pgfqpoint{-0.048611in}{0.000000in}}%
\pgfusepath{stroke,fill}%
}%
\begin{pgfscope}%
\pgfsys@transformshift{0.662500in}{0.307391in}%
\pgfsys@useobject{currentmarker}{}%
\end{pgfscope}%
\end{pgfscope}%
\begin{pgfscope}%
\definecolor{textcolor}{rgb}{0.000000,0.000000,0.000000}%
\pgfsetstrokecolor{textcolor}%
\pgfsetfillcolor{textcolor}%
\pgftext[x=0.243482in, y=0.263988in, left, base]{\color{textcolor}{\rmfamily\fontsize{9.000000}{10.800000}\selectfont\catcode`\^=\active\def^{\ifmmode\sp\else\^{}\fi}\catcode`\%=\active\def%{\%}$\mathdefault{2.705}$}}%
\end{pgfscope}%
\begin{pgfscope}%
\pgfsetbuttcap%
\pgfsetroundjoin%
\definecolor{currentfill}{rgb}{0.000000,0.000000,0.000000}%
\pgfsetfillcolor{currentfill}%
\pgfsetlinewidth{0.803000pt}%
\definecolor{currentstroke}{rgb}{0.000000,0.000000,0.000000}%
\pgfsetstrokecolor{currentstroke}%
\pgfsetdash{}{0pt}%
\pgfsys@defobject{currentmarker}{\pgfqpoint{-0.048611in}{0.000000in}}{\pgfqpoint{-0.000000in}{0.000000in}}{%
\pgfpathmoveto{\pgfqpoint{-0.000000in}{0.000000in}}%
\pgfpathlineto{\pgfqpoint{-0.048611in}{0.000000in}}%
\pgfusepath{stroke,fill}%
}%
\begin{pgfscope}%
\pgfsys@transformshift{0.662500in}{1.441609in}%
\pgfsys@useobject{currentmarker}{}%
\end{pgfscope}%
\end{pgfscope}%
\begin{pgfscope}%
\definecolor{textcolor}{rgb}{0.000000,0.000000,0.000000}%
\pgfsetstrokecolor{textcolor}%
\pgfsetfillcolor{textcolor}%
\pgftext[x=0.243482in, y=1.398206in, left, base]{\color{textcolor}{\rmfamily\fontsize{9.000000}{10.800000}\selectfont\catcode`\^=\active\def^{\ifmmode\sp\else\^{}\fi}\catcode`\%=\active\def%{\%}$\mathdefault{7.634}$}}%
\end{pgfscope}%
\begin{pgfscope}%
\pgfsetrectcap%
\pgfsetmiterjoin%
\pgfsetlinewidth{0.803000pt}%
\definecolor{currentstroke}{rgb}{0.000000,0.000000,0.000000}%
\pgfsetstrokecolor{currentstroke}%
\pgfsetdash{}{0pt}%
\pgfpathmoveto{\pgfqpoint{0.662500in}{0.307391in}}%
\pgfpathlineto{\pgfqpoint{0.690855in}{0.307391in}}%
\pgfpathlineto{\pgfqpoint{0.719211in}{0.307391in}}%
\pgfpathlineto{\pgfqpoint{0.719211in}{1.441609in}}%
\pgfpathlineto{\pgfqpoint{0.690855in}{1.441609in}}%
\pgfpathlineto{\pgfqpoint{0.662500in}{1.441609in}}%
\pgfpathlineto{\pgfqpoint{0.662500in}{0.307391in}}%
\pgfpathclose%
\pgfusepath{stroke}%
\end{pgfscope}%
\begin{pgfscope}%
\pgfsetbuttcap%
\pgfsetmiterjoin%
\definecolor{currentfill}{rgb}{1.000000,1.000000,1.000000}%
\pgfsetfillcolor{currentfill}%
\pgfsetlinewidth{0.000000pt}%
\definecolor{currentstroke}{rgb}{0.000000,0.000000,0.000000}%
\pgfsetstrokecolor{currentstroke}%
\pgfsetstrokeopacity{0.000000}%
\pgfsetdash{}{0pt}%
\pgfpathmoveto{\pgfqpoint{2.070786in}{0.952107in}}%
\pgfpathlineto{\pgfqpoint{2.657571in}{0.952107in}}%
\pgfpathlineto{\pgfqpoint{2.657571in}{1.538893in}}%
\pgfpathlineto{\pgfqpoint{2.070786in}{1.538893in}}%
\pgfpathlineto{\pgfqpoint{2.070786in}{0.952107in}}%
\pgfpathclose%
\pgfusepath{fill}%
\end{pgfscope}%
\begin{pgfscope}%
\pgfpathrectangle{\pgfqpoint{2.070786in}{0.952107in}}{\pgfqpoint{0.586786in}{0.586786in}}%
\pgfusepath{clip}%
\pgfsys@transformshift{2.070786in}{0.952107in}%
\pgftext[left,bottom]{\includegraphics[interpolate=true,width=0.590000in,height=0.590000in]{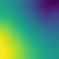}}%
\end{pgfscope}%
\begin{pgfscope}%
\definecolor{textcolor}{rgb}{0.000000,0.000000,0.000000}%
\pgfsetstrokecolor{textcolor}%
\pgfsetfillcolor{textcolor}%
\pgftext[x=2.364179in,y=1.563198in,,base]{\color{textcolor}{\rmfamily\fontsize{9.000000}{10.800000}\selectfont\catcode`\^=\active\def^{\ifmmode\sp\else\^{}\fi}\catcode`\%=\active\def%{\%}Channel 1}}%
\end{pgfscope}%
\begin{pgfscope}%
\pgfsetbuttcap%
\pgfsetmiterjoin%
\definecolor{currentfill}{rgb}{1.000000,1.000000,1.000000}%
\pgfsetfillcolor{currentfill}%
\pgfsetlinewidth{0.000000pt}%
\definecolor{currentstroke}{rgb}{0.000000,0.000000,0.000000}%
\pgfsetstrokecolor{currentstroke}%
\pgfsetstrokeopacity{0.000000}%
\pgfsetdash{}{0pt}%
\pgfpathmoveto{\pgfqpoint{2.774929in}{0.952107in}}%
\pgfpathlineto{\pgfqpoint{3.361714in}{0.952107in}}%
\pgfpathlineto{\pgfqpoint{3.361714in}{1.538893in}}%
\pgfpathlineto{\pgfqpoint{2.774929in}{1.538893in}}%
\pgfpathlineto{\pgfqpoint{2.774929in}{0.952107in}}%
\pgfpathclose%
\pgfusepath{fill}%
\end{pgfscope}%
\begin{pgfscope}%
\pgfpathrectangle{\pgfqpoint{2.774929in}{0.952107in}}{\pgfqpoint{0.586786in}{0.586786in}}%
\pgfusepath{clip}%
\pgfsys@transformshift{2.774929in}{0.952107in}%
\pgftext[left,bottom]{\includegraphics[interpolate=true,width=0.590000in,height=0.590000in]{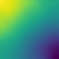}}%
\end{pgfscope}%
\begin{pgfscope}%
\definecolor{textcolor}{rgb}{0.000000,0.000000,0.000000}%
\pgfsetstrokecolor{textcolor}%
\pgfsetfillcolor{textcolor}%
\pgftext[x=3.068321in,y=1.563198in,,base]{\color{textcolor}{\rmfamily\fontsize{9.000000}{10.800000}\selectfont\catcode`\^=\active\def^{\ifmmode\sp\else\^{}\fi}\catcode`\%=\active\def%{\%}Channel 2}}%
\end{pgfscope}%
\begin{pgfscope}%
\pgfsetbuttcap%
\pgfsetmiterjoin%
\definecolor{currentfill}{rgb}{1.000000,1.000000,1.000000}%
\pgfsetfillcolor{currentfill}%
\pgfsetlinewidth{0.000000pt}%
\definecolor{currentstroke}{rgb}{0.000000,0.000000,0.000000}%
\pgfsetstrokecolor{currentstroke}%
\pgfsetstrokeopacity{0.000000}%
\pgfsetdash{}{0pt}%
\pgfpathmoveto{\pgfqpoint{3.479071in}{0.952107in}}%
\pgfpathlineto{\pgfqpoint{4.065857in}{0.952107in}}%
\pgfpathlineto{\pgfqpoint{4.065857in}{1.538893in}}%
\pgfpathlineto{\pgfqpoint{3.479071in}{1.538893in}}%
\pgfpathlineto{\pgfqpoint{3.479071in}{0.952107in}}%
\pgfpathclose%
\pgfusepath{fill}%
\end{pgfscope}%
\begin{pgfscope}%
\pgfpathrectangle{\pgfqpoint{3.479071in}{0.952107in}}{\pgfqpoint{0.586786in}{0.586786in}}%
\pgfusepath{clip}%
\pgfsys@transformshift{3.479071in}{0.952107in}%
\pgftext[left,bottom]{\includegraphics[interpolate=true,width=0.590000in,height=0.590000in]{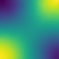}}%
\end{pgfscope}%
\begin{pgfscope}%
\definecolor{textcolor}{rgb}{0.000000,0.000000,0.000000}%
\pgfsetstrokecolor{textcolor}%
\pgfsetfillcolor{textcolor}%
\pgftext[x=3.772464in,y=1.563198in,,base]{\color{textcolor}{\rmfamily\fontsize{9.000000}{10.800000}\selectfont\catcode`\^=\active\def^{\ifmmode\sp\else\^{}\fi}\catcode`\%=\active\def%{\%}Channel 3}}%
\end{pgfscope}%
\begin{pgfscope}%
\pgfsetbuttcap%
\pgfsetmiterjoin%
\definecolor{currentfill}{rgb}{1.000000,1.000000,1.000000}%
\pgfsetfillcolor{currentfill}%
\pgfsetlinewidth{0.000000pt}%
\definecolor{currentstroke}{rgb}{0.000000,0.000000,0.000000}%
\pgfsetstrokecolor{currentstroke}%
\pgfsetstrokeopacity{0.000000}%
\pgfsetdash{}{0pt}%
\pgfpathmoveto{\pgfqpoint{4.183214in}{0.952107in}}%
\pgfpathlineto{\pgfqpoint{4.770000in}{0.952107in}}%
\pgfpathlineto{\pgfqpoint{4.770000in}{1.538893in}}%
\pgfpathlineto{\pgfqpoint{4.183214in}{1.538893in}}%
\pgfpathlineto{\pgfqpoint{4.183214in}{0.952107in}}%
\pgfpathclose%
\pgfusepath{fill}%
\end{pgfscope}%
\begin{pgfscope}%
\pgfpathrectangle{\pgfqpoint{4.183214in}{0.952107in}}{\pgfqpoint{0.586786in}{0.586786in}}%
\pgfusepath{clip}%
\pgfsys@transformshift{4.183214in}{0.952107in}%
\pgftext[left,bottom]{\includegraphics[interpolate=true,width=0.590000in,height=0.590000in]{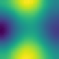}}%
\end{pgfscope}%
\begin{pgfscope}%
\definecolor{textcolor}{rgb}{0.000000,0.000000,0.000000}%
\pgfsetstrokecolor{textcolor}%
\pgfsetfillcolor{textcolor}%
\pgftext[x=4.476607in,y=1.563198in,,base]{\color{textcolor}{\rmfamily\fontsize{9.000000}{10.800000}\selectfont\catcode`\^=\active\def^{\ifmmode\sp\else\^{}\fi}\catcode`\%=\active\def%{\%}Channel 4}}%
\end{pgfscope}%
\begin{pgfscope}%
\pgfsetbuttcap%
\pgfsetmiterjoin%
\definecolor{currentfill}{rgb}{1.000000,1.000000,1.000000}%
\pgfsetfillcolor{currentfill}%
\pgfsetlinewidth{0.000000pt}%
\definecolor{currentstroke}{rgb}{0.000000,0.000000,0.000000}%
\pgfsetstrokecolor{currentstroke}%
\pgfsetstrokeopacity{0.000000}%
\pgfsetdash{}{0pt}%
\pgfpathmoveto{\pgfqpoint{2.070786in}{0.210107in}}%
\pgfpathlineto{\pgfqpoint{2.657571in}{0.210107in}}%
\pgfpathlineto{\pgfqpoint{2.657571in}{0.796893in}}%
\pgfpathlineto{\pgfqpoint{2.070786in}{0.796893in}}%
\pgfpathlineto{\pgfqpoint{2.070786in}{0.210107in}}%
\pgfpathclose%
\pgfusepath{fill}%
\end{pgfscope}%
\begin{pgfscope}%
\pgfpathrectangle{\pgfqpoint{2.070786in}{0.210107in}}{\pgfqpoint{0.586786in}{0.586786in}}%
\pgfusepath{clip}%
\pgfsys@transformshift{2.070786in}{0.210107in}%
\pgftext[left,bottom]{\includegraphics[interpolate=true,width=0.590000in,height=0.590000in]{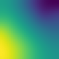}}%
\end{pgfscope}%
\begin{pgfscope}%
\definecolor{textcolor}{rgb}{0.000000,0.000000,0.000000}%
\pgfsetstrokecolor{textcolor}%
\pgfsetfillcolor{textcolor}%
\pgftext[x=2.364179in,y=0.821198in,,base]{\color{textcolor}{\rmfamily\fontsize{9.000000}{10.800000}\selectfont\catcode`\^=\active\def^{\ifmmode\sp\else\^{}\fi}\catcode`\%=\active\def%{\%}Label 1}}%
\end{pgfscope}%
\begin{pgfscope}%
\pgfsetbuttcap%
\pgfsetmiterjoin%
\definecolor{currentfill}{rgb}{1.000000,1.000000,1.000000}%
\pgfsetfillcolor{currentfill}%
\pgfsetlinewidth{0.000000pt}%
\definecolor{currentstroke}{rgb}{0.000000,0.000000,0.000000}%
\pgfsetstrokecolor{currentstroke}%
\pgfsetstrokeopacity{0.000000}%
\pgfsetdash{}{0pt}%
\pgfpathmoveto{\pgfqpoint{2.774929in}{0.210107in}}%
\pgfpathlineto{\pgfqpoint{3.361714in}{0.210107in}}%
\pgfpathlineto{\pgfqpoint{3.361714in}{0.796893in}}%
\pgfpathlineto{\pgfqpoint{2.774929in}{0.796893in}}%
\pgfpathlineto{\pgfqpoint{2.774929in}{0.210107in}}%
\pgfpathclose%
\pgfusepath{fill}%
\end{pgfscope}%
\begin{pgfscope}%
\pgfpathrectangle{\pgfqpoint{2.774929in}{0.210107in}}{\pgfqpoint{0.586786in}{0.586786in}}%
\pgfusepath{clip}%
\pgfsys@transformshift{2.774929in}{0.210107in}%
\pgftext[left,bottom]{\includegraphics[interpolate=true,width=0.590000in,height=0.590000in]{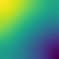}}%
\end{pgfscope}%
\begin{pgfscope}%
\definecolor{textcolor}{rgb}{0.000000,0.000000,0.000000}%
\pgfsetstrokecolor{textcolor}%
\pgfsetfillcolor{textcolor}%
\pgftext[x=3.068321in,y=0.821198in,,base]{\color{textcolor}{\rmfamily\fontsize{9.000000}{10.800000}\selectfont\catcode`\^=\active\def^{\ifmmode\sp\else\^{}\fi}\catcode`\%=\active\def%{\%}Label 2}}%
\end{pgfscope}%
\begin{pgfscope}%
\pgfsetbuttcap%
\pgfsetmiterjoin%
\definecolor{currentfill}{rgb}{1.000000,1.000000,1.000000}%
\pgfsetfillcolor{currentfill}%
\pgfsetlinewidth{0.000000pt}%
\definecolor{currentstroke}{rgb}{0.000000,0.000000,0.000000}%
\pgfsetstrokecolor{currentstroke}%
\pgfsetstrokeopacity{0.000000}%
\pgfsetdash{}{0pt}%
\pgfpathmoveto{\pgfqpoint{3.479071in}{0.210107in}}%
\pgfpathlineto{\pgfqpoint{4.065857in}{0.210107in}}%
\pgfpathlineto{\pgfqpoint{4.065857in}{0.796893in}}%
\pgfpathlineto{\pgfqpoint{3.479071in}{0.796893in}}%
\pgfpathlineto{\pgfqpoint{3.479071in}{0.210107in}}%
\pgfpathclose%
\pgfusepath{fill}%
\end{pgfscope}%
\begin{pgfscope}%
\pgfpathrectangle{\pgfqpoint{3.479071in}{0.210107in}}{\pgfqpoint{0.586786in}{0.586786in}}%
\pgfusepath{clip}%
\pgfsys@transformshift{3.479071in}{0.210107in}%
\pgftext[left,bottom]{\includegraphics[interpolate=true,width=0.590000in,height=0.590000in]{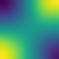}}%
\end{pgfscope}%
\begin{pgfscope}%
\definecolor{textcolor}{rgb}{0.000000,0.000000,0.000000}%
\pgfsetstrokecolor{textcolor}%
\pgfsetfillcolor{textcolor}%
\pgftext[x=3.772464in,y=0.821198in,,base]{\color{textcolor}{\rmfamily\fontsize{9.000000}{10.800000}\selectfont\catcode`\^=\active\def^{\ifmmode\sp\else\^{}\fi}\catcode`\%=\active\def%{\%}Label 3}}%
\end{pgfscope}%
\begin{pgfscope}%
\pgfsetbuttcap%
\pgfsetmiterjoin%
\definecolor{currentfill}{rgb}{1.000000,1.000000,1.000000}%
\pgfsetfillcolor{currentfill}%
\pgfsetlinewidth{0.000000pt}%
\definecolor{currentstroke}{rgb}{0.000000,0.000000,0.000000}%
\pgfsetstrokecolor{currentstroke}%
\pgfsetstrokeopacity{0.000000}%
\pgfsetdash{}{0pt}%
\pgfpathmoveto{\pgfqpoint{4.183214in}{0.210107in}}%
\pgfpathlineto{\pgfqpoint{4.770000in}{0.210107in}}%
\pgfpathlineto{\pgfqpoint{4.770000in}{0.796893in}}%
\pgfpathlineto{\pgfqpoint{4.183214in}{0.796893in}}%
\pgfpathlineto{\pgfqpoint{4.183214in}{0.210107in}}%
\pgfpathclose%
\pgfusepath{fill}%
\end{pgfscope}%
\begin{pgfscope}%
\pgfpathrectangle{\pgfqpoint{4.183214in}{0.210107in}}{\pgfqpoint{0.586786in}{0.586786in}}%
\pgfusepath{clip}%
\pgfsys@transformshift{4.183214in}{0.210107in}%
\pgftext[left,bottom]{\includegraphics[interpolate=true,width=0.590000in,height=0.590000in]{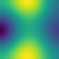}}%
\end{pgfscope}%
\begin{pgfscope}%
\definecolor{textcolor}{rgb}{0.000000,0.000000,0.000000}%
\pgfsetstrokecolor{textcolor}%
\pgfsetfillcolor{textcolor}%
\pgftext[x=4.476607in,y=0.821198in,,base]{\color{textcolor}{\rmfamily\fontsize{9.000000}{10.800000}\selectfont\catcode`\^=\active\def^{\ifmmode\sp\else\^{}\fi}\catcode`\%=\active\def%{\%}Label 4}}%
\end{pgfscope}%
\end{pgfpicture}%
\makeatother%
\endgroup%

%% file: _x.pgf
%% Creator: Matplotlib, PGF backend
%%
%% To include the figure in your LaTeX document, write
%%   \input{<filename>.pgf}
%%
%% Make sure the required packages are loaded in your preamble
%%   \usepackage{pgf}
%%
%% Also ensure that all the required font packages are loaded; for instance,
%% the lmodern package is sometimes necessary when using math font.
%%   \usepackage{lmodern}
%%
%% Figures using additional raster images can only be included by \input if
%% they are in the same directory as the main LaTeX file. For loading figures
%% from other directories you can use the `import` package
%%   \usepackage{import}
%%
%% and then include the figures with
%%   \import{<path to file>}{<filename>.pgf}
%%
%% Matplotlib used the following preamble
%%   \def\mathdefault#1{#1}
%%   \everymath=\expandafter{\the\everymath\displaystyle}
%%   \usepackage{bm}%
%%   \usepackage{mismath}
%%   \newcommand{\SSSText}[1]{{\scriptscriptstyle \mathup{#1}}}%
%%   \renewcommand{\mathdefault}[1][]{}%
%%   \ifdefined\pdftexversion\else  % non-pdftex case.
%%     \usepackage{fontspec}
%%   \fi
%%   \makeatletter\@ifpackageloaded{underscore}{}{\usepackage[strings]{underscore}}\makeatother
%%
\begingroup%
\makeatletter%
\begin{pgfpicture}%
\pgfpathrectangle{\pgfpointorigin}{\pgfqpoint{4.307500in}{3.546273in}}%
\pgfusepath{use as bounding box, clip}%
\begin{pgfscope}%
\pgfsetbuttcap%
\pgfsetmiterjoin%
\definecolor{currentfill}{rgb}{1.000000,1.000000,1.000000}%
\pgfsetfillcolor{currentfill}%
\pgfsetlinewidth{0.000000pt}%
\definecolor{currentstroke}{rgb}{1.000000,1.000000,1.000000}%
\pgfsetstrokecolor{currentstroke}%
\pgfsetdash{}{0pt}%
\pgfpathmoveto{\pgfqpoint{0.000000in}{0.000000in}}%
\pgfpathlineto{\pgfqpoint{4.307500in}{0.000000in}}%
\pgfpathlineto{\pgfqpoint{4.307500in}{3.546273in}}%
\pgfpathlineto{\pgfqpoint{0.000000in}{3.546273in}}%
\pgfpathlineto{\pgfqpoint{0.000000in}{0.000000in}}%
\pgfpathclose%
\pgfusepath{fill}%
\end{pgfscope}%
\begin{pgfscope}%
\pgfsetbuttcap%
\pgfsetmiterjoin%
\definecolor{currentfill}{rgb}{1.000000,1.000000,1.000000}%
\pgfsetfillcolor{currentfill}%
\pgfsetlinewidth{0.000000pt}%
\definecolor{currentstroke}{rgb}{0.000000,0.000000,0.000000}%
\pgfsetstrokecolor{currentstroke}%
\pgfsetstrokeopacity{0.000000}%
\pgfsetdash{}{0pt}%
\pgfpathmoveto{\pgfqpoint{0.100000in}{2.356379in}}%
\pgfpathlineto{\pgfqpoint{1.061330in}{2.356379in}}%
\pgfpathlineto{\pgfqpoint{1.061330in}{3.317709in}}%
\pgfpathlineto{\pgfqpoint{0.100000in}{3.317709in}}%
\pgfpathlineto{\pgfqpoint{0.100000in}{2.356379in}}%
\pgfpathclose%
\pgfusepath{fill}%
\end{pgfscope}%
\begin{pgfscope}%
\pgfpathrectangle{\pgfqpoint{0.100000in}{2.356379in}}{\pgfqpoint{0.961330in}{0.961330in}}%
\pgfusepath{clip}%
\pgfsys@transformshift{0.100000in}{2.356379in}%
\pgftext[left,bottom]{\includegraphics[interpolate=true,width=0.970000in,height=0.970000in]{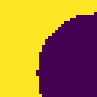}}%
\end{pgfscope}%
\begin{pgfscope}%
\definecolor{textcolor}{rgb}{0.000000,0.000000,0.000000}%
\pgfsetstrokecolor{textcolor}%
\pgfsetfillcolor{textcolor}%
\pgftext[x=0.580665in,y=3.342014in,,base]{\color{textcolor}{\rmfamily\fontsize{9.000000}{10.800000}\selectfont\catcode`\^=\active\def^{\ifmmode\sp\else\^{}\fi}\catcode`\%=\active\def%{\%}$\kappa_b = 10^{-5}$}}%
\end{pgfscope}%
\begin{pgfscope}%
\pgfsetbuttcap%
\pgfsetmiterjoin%
\definecolor{currentfill}{rgb}{1.000000,1.000000,1.000000}%
\pgfsetfillcolor{currentfill}%
\pgfsetlinewidth{0.000000pt}%
\definecolor{currentstroke}{rgb}{0.000000,0.000000,0.000000}%
\pgfsetstrokecolor{currentstroke}%
\pgfsetstrokeopacity{0.000000}%
\pgfsetdash{}{0pt}%
\pgfpathmoveto{\pgfqpoint{0.624362in}{1.779040in}}%
\pgfpathlineto{\pgfqpoint{1.061330in}{1.779040in}}%
\pgfpathlineto{\pgfqpoint{1.061330in}{2.216008in}}%
\pgfpathlineto{\pgfqpoint{0.624362in}{2.216008in}}%
\pgfpathlineto{\pgfqpoint{0.624362in}{1.779040in}}%
\pgfpathclose%
\pgfusepath{fill}%
\end{pgfscope}%
\begin{pgfscope}%
\pgfpathrectangle{\pgfqpoint{0.624362in}{1.779040in}}{\pgfqpoint{0.436968in}{0.436968in}}%
\pgfusepath{clip}%
\pgfsys@transformshift{0.624362in}{1.779040in}%
\pgftext[left,bottom]{\includegraphics[interpolate=true,width=0.440000in,height=0.440000in]{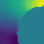}}%
\end{pgfscope}%
\begin{pgfscope}%
\definecolor{textcolor}{rgb}{0.000000,0.000000,0.000000}%
\pgfsetstrokecolor{textcolor}%
\pgfsetfillcolor{textcolor}%
\pgftext[x=0.842846in,y=2.240314in,,base]{\color{textcolor}{\rmfamily\fontsize{9.000000}{10.800000}\selectfont\catcode`\^=\active\def^{\ifmmode\sp\else\^{}\fi}\catcode`\%=\active\def%{\%}C 1}}%
\end{pgfscope}%
\begin{pgfscope}%
\pgfsetbuttcap%
\pgfsetmiterjoin%
\definecolor{currentfill}{rgb}{1.000000,1.000000,1.000000}%
\pgfsetfillcolor{currentfill}%
\pgfsetlinewidth{0.000000pt}%
\definecolor{currentstroke}{rgb}{0.000000,0.000000,0.000000}%
\pgfsetstrokecolor{currentstroke}%
\pgfsetstrokeopacity{0.000000}%
\pgfsetdash{}{0pt}%
\pgfpathmoveto{\pgfqpoint{0.624362in}{1.219360in}}%
\pgfpathlineto{\pgfqpoint{1.061330in}{1.219360in}}%
\pgfpathlineto{\pgfqpoint{1.061330in}{1.656328in}}%
\pgfpathlineto{\pgfqpoint{0.624362in}{1.656328in}}%
\pgfpathlineto{\pgfqpoint{0.624362in}{1.219360in}}%
\pgfpathclose%
\pgfusepath{fill}%
\end{pgfscope}%
\begin{pgfscope}%
\pgfpathrectangle{\pgfqpoint{0.624362in}{1.219360in}}{\pgfqpoint{0.436968in}{0.436968in}}%
\pgfusepath{clip}%
\pgfsys@transformshift{0.624362in}{1.219360in}%
\pgftext[left,bottom]{\includegraphics[interpolate=true,width=0.440000in,height=0.440000in]{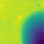}}%
\end{pgfscope}%
\begin{pgfscope}%
\definecolor{textcolor}{rgb}{0.000000,0.000000,0.000000}%
\pgfsetstrokecolor{textcolor}%
\pgfsetfillcolor{textcolor}%
\pgftext[x=0.842846in,y=1.680634in,,base]{\color{textcolor}{\rmfamily\fontsize{9.000000}{10.800000}\selectfont\catcode`\^=\active\def^{\ifmmode\sp\else\^{}\fi}\catcode`\%=\active\def%{\%}C 2}}%
\end{pgfscope}%
\begin{pgfscope}%
\pgfsetbuttcap%
\pgfsetmiterjoin%
\definecolor{currentfill}{rgb}{1.000000,1.000000,1.000000}%
\pgfsetfillcolor{currentfill}%
\pgfsetlinewidth{0.000000pt}%
\definecolor{currentstroke}{rgb}{0.000000,0.000000,0.000000}%
\pgfsetstrokecolor{currentstroke}%
\pgfsetstrokeopacity{0.000000}%
\pgfsetdash{}{0pt}%
\pgfpathmoveto{\pgfqpoint{0.624362in}{0.659680in}}%
\pgfpathlineto{\pgfqpoint{1.061330in}{0.659680in}}%
\pgfpathlineto{\pgfqpoint{1.061330in}{1.096648in}}%
\pgfpathlineto{\pgfqpoint{0.624362in}{1.096648in}}%
\pgfpathlineto{\pgfqpoint{0.624362in}{0.659680in}}%
\pgfpathclose%
\pgfusepath{fill}%
\end{pgfscope}%
\begin{pgfscope}%
\pgfpathrectangle{\pgfqpoint{0.624362in}{0.659680in}}{\pgfqpoint{0.436968in}{0.436968in}}%
\pgfusepath{clip}%
\pgfsys@transformshift{0.624362in}{0.659680in}%
\pgftext[left,bottom]{\includegraphics[interpolate=true,width=0.440000in,height=0.440000in]{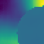}}%
\end{pgfscope}%
\begin{pgfscope}%
\definecolor{textcolor}{rgb}{0.000000,0.000000,0.000000}%
\pgfsetstrokecolor{textcolor}%
\pgfsetfillcolor{textcolor}%
\pgftext[x=0.842846in,y=1.120954in,,base]{\color{textcolor}{\rmfamily\fontsize{9.000000}{10.800000}\selectfont\catcode`\^=\active\def^{\ifmmode\sp\else\^{}\fi}\catcode`\%=\active\def%{\%}C 3}}%
\end{pgfscope}%
\begin{pgfscope}%
\pgfsetbuttcap%
\pgfsetmiterjoin%
\definecolor{currentfill}{rgb}{1.000000,1.000000,1.000000}%
\pgfsetfillcolor{currentfill}%
\pgfsetlinewidth{0.000000pt}%
\definecolor{currentstroke}{rgb}{0.000000,0.000000,0.000000}%
\pgfsetstrokecolor{currentstroke}%
\pgfsetstrokeopacity{0.000000}%
\pgfsetdash{}{0pt}%
\pgfpathmoveto{\pgfqpoint{0.624362in}{0.100000in}}%
\pgfpathlineto{\pgfqpoint{1.061330in}{0.100000in}}%
\pgfpathlineto{\pgfqpoint{1.061330in}{0.536968in}}%
\pgfpathlineto{\pgfqpoint{0.624362in}{0.536968in}}%
\pgfpathlineto{\pgfqpoint{0.624362in}{0.100000in}}%
\pgfpathclose%
\pgfusepath{fill}%
\end{pgfscope}%
\begin{pgfscope}%
\pgfpathrectangle{\pgfqpoint{0.624362in}{0.100000in}}{\pgfqpoint{0.436968in}{0.436968in}}%
\pgfusepath{clip}%
\pgfsys@transformshift{0.624362in}{0.100000in}%
\pgftext[left,bottom]{\includegraphics[interpolate=true,width=0.440000in,height=0.440000in]{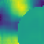}}%
\end{pgfscope}%
\begin{pgfscope}%
\definecolor{textcolor}{rgb}{0.000000,0.000000,0.000000}%
\pgfsetstrokecolor{textcolor}%
\pgfsetfillcolor{textcolor}%
\pgftext[x=0.842846in,y=0.561274in,,base]{\color{textcolor}{\rmfamily\fontsize{9.000000}{10.800000}\selectfont\catcode`\^=\active\def^{\ifmmode\sp\else\^{}\fi}\catcode`\%=\active\def%{\%}C 4}}%
\end{pgfscope}%
\begin{pgfscope}%
\pgfsetbuttcap%
\pgfsetmiterjoin%
\definecolor{currentfill}{rgb}{1.000000,1.000000,1.000000}%
\pgfsetfillcolor{currentfill}%
\pgfsetlinewidth{0.000000pt}%
\definecolor{currentstroke}{rgb}{0.000000,0.000000,0.000000}%
\pgfsetstrokecolor{currentstroke}%
\pgfsetstrokeopacity{0.000000}%
\pgfsetdash{}{0pt}%
\pgfpathmoveto{\pgfqpoint{0.100000in}{1.779040in}}%
\pgfpathlineto{\pgfqpoint{0.536968in}{1.779040in}}%
\pgfpathlineto{\pgfqpoint{0.536968in}{2.216008in}}%
\pgfpathlineto{\pgfqpoint{0.100000in}{2.216008in}}%
\pgfpathlineto{\pgfqpoint{0.100000in}{1.779040in}}%
\pgfpathclose%
\pgfusepath{fill}%
\end{pgfscope}%
\begin{pgfscope}%
\pgfpathrectangle{\pgfqpoint{0.100000in}{1.779040in}}{\pgfqpoint{0.436968in}{0.436968in}}%
\pgfusepath{clip}%
\pgfsys@transformshift{0.100000in}{1.779040in}%
\pgftext[left,bottom]{\includegraphics[interpolate=true,width=0.440000in,height=0.440000in]{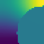}}%
\end{pgfscope}%
\begin{pgfscope}%
\definecolor{textcolor}{rgb}{0.000000,0.000000,0.000000}%
\pgfsetstrokecolor{textcolor}%
\pgfsetfillcolor{textcolor}%
\pgftext[x=0.318484in,y=2.240314in,,base]{\color{textcolor}{\rmfamily\fontsize{9.000000}{10.800000}\selectfont\catcode`\^=\active\def^{\ifmmode\sp\else\^{}\fi}\catcode`\%=\active\def%{\%}L 1}}%
\end{pgfscope}%
\begin{pgfscope}%
\pgfsetbuttcap%
\pgfsetmiterjoin%
\definecolor{currentfill}{rgb}{1.000000,1.000000,1.000000}%
\pgfsetfillcolor{currentfill}%
\pgfsetlinewidth{0.000000pt}%
\definecolor{currentstroke}{rgb}{0.000000,0.000000,0.000000}%
\pgfsetstrokecolor{currentstroke}%
\pgfsetstrokeopacity{0.000000}%
\pgfsetdash{}{0pt}%
\pgfpathmoveto{\pgfqpoint{0.100000in}{1.219360in}}%
\pgfpathlineto{\pgfqpoint{0.536968in}{1.219360in}}%
\pgfpathlineto{\pgfqpoint{0.536968in}{1.656328in}}%
\pgfpathlineto{\pgfqpoint{0.100000in}{1.656328in}}%
\pgfpathlineto{\pgfqpoint{0.100000in}{1.219360in}}%
\pgfpathclose%
\pgfusepath{fill}%
\end{pgfscope}%
\begin{pgfscope}%
\pgfpathrectangle{\pgfqpoint{0.100000in}{1.219360in}}{\pgfqpoint{0.436968in}{0.436968in}}%
\pgfusepath{clip}%
\pgfsys@transformshift{0.100000in}{1.219360in}%
\pgftext[left,bottom]{\includegraphics[interpolate=true,width=0.440000in,height=0.440000in]{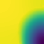}}%
\end{pgfscope}%
\begin{pgfscope}%
\definecolor{textcolor}{rgb}{0.000000,0.000000,0.000000}%
\pgfsetstrokecolor{textcolor}%
\pgfsetfillcolor{textcolor}%
\pgftext[x=0.318484in,y=1.680634in,,base]{\color{textcolor}{\rmfamily\fontsize{9.000000}{10.800000}\selectfont\catcode`\^=\active\def^{\ifmmode\sp\else\^{}\fi}\catcode`\%=\active\def%{\%}L 2}}%
\end{pgfscope}%
\begin{pgfscope}%
\pgfsetbuttcap%
\pgfsetmiterjoin%
\definecolor{currentfill}{rgb}{1.000000,1.000000,1.000000}%
\pgfsetfillcolor{currentfill}%
\pgfsetlinewidth{0.000000pt}%
\definecolor{currentstroke}{rgb}{0.000000,0.000000,0.000000}%
\pgfsetstrokecolor{currentstroke}%
\pgfsetstrokeopacity{0.000000}%
\pgfsetdash{}{0pt}%
\pgfpathmoveto{\pgfqpoint{0.100000in}{0.659680in}}%
\pgfpathlineto{\pgfqpoint{0.536968in}{0.659680in}}%
\pgfpathlineto{\pgfqpoint{0.536968in}{1.096648in}}%
\pgfpathlineto{\pgfqpoint{0.100000in}{1.096648in}}%
\pgfpathlineto{\pgfqpoint{0.100000in}{0.659680in}}%
\pgfpathclose%
\pgfusepath{fill}%
\end{pgfscope}%
\begin{pgfscope}%
\pgfpathrectangle{\pgfqpoint{0.100000in}{0.659680in}}{\pgfqpoint{0.436968in}{0.436968in}}%
\pgfusepath{clip}%
\pgfsys@transformshift{0.100000in}{0.659680in}%
\pgftext[left,bottom]{\includegraphics[interpolate=true,width=0.440000in,height=0.440000in]{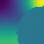}}%
\end{pgfscope}%
\begin{pgfscope}%
\definecolor{textcolor}{rgb}{0.000000,0.000000,0.000000}%
\pgfsetstrokecolor{textcolor}%
\pgfsetfillcolor{textcolor}%
\pgftext[x=0.318484in,y=1.120954in,,base]{\color{textcolor}{\rmfamily\fontsize{9.000000}{10.800000}\selectfont\catcode`\^=\active\def^{\ifmmode\sp\else\^{}\fi}\catcode`\%=\active\def%{\%}L 3}}%
\end{pgfscope}%
\begin{pgfscope}%
\pgfsetbuttcap%
\pgfsetmiterjoin%
\definecolor{currentfill}{rgb}{1.000000,1.000000,1.000000}%
\pgfsetfillcolor{currentfill}%
\pgfsetlinewidth{0.000000pt}%
\definecolor{currentstroke}{rgb}{0.000000,0.000000,0.000000}%
\pgfsetstrokecolor{currentstroke}%
\pgfsetstrokeopacity{0.000000}%
\pgfsetdash{}{0pt}%
\pgfpathmoveto{\pgfqpoint{0.100000in}{0.100000in}}%
\pgfpathlineto{\pgfqpoint{0.536968in}{0.100000in}}%
\pgfpathlineto{\pgfqpoint{0.536968in}{0.536968in}}%
\pgfpathlineto{\pgfqpoint{0.100000in}{0.536968in}}%
\pgfpathlineto{\pgfqpoint{0.100000in}{0.100000in}}%
\pgfpathclose%
\pgfusepath{fill}%
\end{pgfscope}%
\begin{pgfscope}%
\pgfpathrectangle{\pgfqpoint{0.100000in}{0.100000in}}{\pgfqpoint{0.436968in}{0.436968in}}%
\pgfusepath{clip}%
\pgfsys@transformshift{0.100000in}{0.100000in}%
\pgftext[left,bottom]{\includegraphics[interpolate=true,width=0.440000in,height=0.440000in]{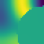}}%
\end{pgfscope}%
\begin{pgfscope}%
\definecolor{textcolor}{rgb}{0.000000,0.000000,0.000000}%
\pgfsetstrokecolor{textcolor}%
\pgfsetfillcolor{textcolor}%
\pgftext[x=0.318484in,y=0.561274in,,base]{\color{textcolor}{\rmfamily\fontsize{9.000000}{10.800000}\selectfont\catcode`\^=\active\def^{\ifmmode\sp\else\^{}\fi}\catcode`\%=\active\def%{\%}L 4}}%
\end{pgfscope}%
\begin{pgfscope}%
\pgfsetbuttcap%
\pgfsetmiterjoin%
\definecolor{currentfill}{rgb}{1.000000,1.000000,1.000000}%
\pgfsetfillcolor{currentfill}%
\pgfsetlinewidth{0.000000pt}%
\definecolor{currentstroke}{rgb}{0.000000,0.000000,0.000000}%
\pgfsetstrokecolor{currentstroke}%
\pgfsetstrokeopacity{0.000000}%
\pgfsetdash{}{0pt}%
\pgfpathmoveto{\pgfqpoint{1.148723in}{2.356379in}}%
\pgfpathlineto{\pgfqpoint{2.110053in}{2.356379in}}%
\pgfpathlineto{\pgfqpoint{2.110053in}{3.317709in}}%
\pgfpathlineto{\pgfqpoint{1.148723in}{3.317709in}}%
\pgfpathlineto{\pgfqpoint{1.148723in}{2.356379in}}%
\pgfpathclose%
\pgfusepath{fill}%
\end{pgfscope}%
\begin{pgfscope}%
\pgfpathrectangle{\pgfqpoint{1.148723in}{2.356379in}}{\pgfqpoint{0.961330in}{0.961330in}}%
\pgfusepath{clip}%
\pgfsys@transformshift{1.148723in}{2.356379in}%
\pgftext[left,bottom]{\includegraphics[interpolate=true,width=0.970000in,height=0.970000in]{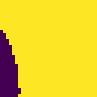}}%
\end{pgfscope}%
\begin{pgfscope}%
\definecolor{textcolor}{rgb}{0.000000,0.000000,0.000000}%
\pgfsetstrokecolor{textcolor}%
\pgfsetfillcolor{textcolor}%
\pgftext[x=1.629388in,y=3.342014in,,base]{\color{textcolor}{\rmfamily\fontsize{9.000000}{10.800000}\selectfont\catcode`\^=\active\def^{\ifmmode\sp\else\^{}\fi}\catcode`\%=\active\def%{\%}$\kappa_b = 10^{-4}$}}%
\end{pgfscope}%
\begin{pgfscope}%
\pgfsetbuttcap%
\pgfsetmiterjoin%
\definecolor{currentfill}{rgb}{1.000000,1.000000,1.000000}%
\pgfsetfillcolor{currentfill}%
\pgfsetlinewidth{0.000000pt}%
\definecolor{currentstroke}{rgb}{0.000000,0.000000,0.000000}%
\pgfsetstrokecolor{currentstroke}%
\pgfsetstrokeopacity{0.000000}%
\pgfsetdash{}{0pt}%
\pgfpathmoveto{\pgfqpoint{1.673085in}{1.779040in}}%
\pgfpathlineto{\pgfqpoint{2.110053in}{1.779040in}}%
\pgfpathlineto{\pgfqpoint{2.110053in}{2.216008in}}%
\pgfpathlineto{\pgfqpoint{1.673085in}{2.216008in}}%
\pgfpathlineto{\pgfqpoint{1.673085in}{1.779040in}}%
\pgfpathclose%
\pgfusepath{fill}%
\end{pgfscope}%
\begin{pgfscope}%
\pgfpathrectangle{\pgfqpoint{1.673085in}{1.779040in}}{\pgfqpoint{0.436968in}{0.436968in}}%
\pgfusepath{clip}%
\pgfsys@transformshift{1.673085in}{1.779040in}%
\pgftext[left,bottom]{\includegraphics[interpolate=true,width=0.440000in,height=0.440000in]{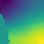}}%
\end{pgfscope}%
\begin{pgfscope}%
\definecolor{textcolor}{rgb}{0.000000,0.000000,0.000000}%
\pgfsetstrokecolor{textcolor}%
\pgfsetfillcolor{textcolor}%
\pgftext[x=1.891569in,y=2.240314in,,base]{\color{textcolor}{\rmfamily\fontsize{9.000000}{10.800000}\selectfont\catcode`\^=\active\def^{\ifmmode\sp\else\^{}\fi}\catcode`\%=\active\def%{\%}C 1}}%
\end{pgfscope}%
\begin{pgfscope}%
\pgfsetbuttcap%
\pgfsetmiterjoin%
\definecolor{currentfill}{rgb}{1.000000,1.000000,1.000000}%
\pgfsetfillcolor{currentfill}%
\pgfsetlinewidth{0.000000pt}%
\definecolor{currentstroke}{rgb}{0.000000,0.000000,0.000000}%
\pgfsetstrokecolor{currentstroke}%
\pgfsetstrokeopacity{0.000000}%
\pgfsetdash{}{0pt}%
\pgfpathmoveto{\pgfqpoint{1.673085in}{1.219360in}}%
\pgfpathlineto{\pgfqpoint{2.110053in}{1.219360in}}%
\pgfpathlineto{\pgfqpoint{2.110053in}{1.656328in}}%
\pgfpathlineto{\pgfqpoint{1.673085in}{1.656328in}}%
\pgfpathlineto{\pgfqpoint{1.673085in}{1.219360in}}%
\pgfpathclose%
\pgfusepath{fill}%
\end{pgfscope}%
\begin{pgfscope}%
\pgfpathrectangle{\pgfqpoint{1.673085in}{1.219360in}}{\pgfqpoint{0.436968in}{0.436968in}}%
\pgfusepath{clip}%
\pgfsys@transformshift{1.673085in}{1.219360in}%
\pgftext[left,bottom]{\includegraphics[interpolate=true,width=0.440000in,height=0.440000in]{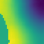}}%
\end{pgfscope}%
\begin{pgfscope}%
\definecolor{textcolor}{rgb}{0.000000,0.000000,0.000000}%
\pgfsetstrokecolor{textcolor}%
\pgfsetfillcolor{textcolor}%
\pgftext[x=1.891569in,y=1.680634in,,base]{\color{textcolor}{\rmfamily\fontsize{9.000000}{10.800000}\selectfont\catcode`\^=\active\def^{\ifmmode\sp\else\^{}\fi}\catcode`\%=\active\def%{\%}C 2}}%
\end{pgfscope}%
\begin{pgfscope}%
\pgfsetbuttcap%
\pgfsetmiterjoin%
\definecolor{currentfill}{rgb}{1.000000,1.000000,1.000000}%
\pgfsetfillcolor{currentfill}%
\pgfsetlinewidth{0.000000pt}%
\definecolor{currentstroke}{rgb}{0.000000,0.000000,0.000000}%
\pgfsetstrokecolor{currentstroke}%
\pgfsetstrokeopacity{0.000000}%
\pgfsetdash{}{0pt}%
\pgfpathmoveto{\pgfqpoint{1.673085in}{0.659680in}}%
\pgfpathlineto{\pgfqpoint{2.110053in}{0.659680in}}%
\pgfpathlineto{\pgfqpoint{2.110053in}{1.096648in}}%
\pgfpathlineto{\pgfqpoint{1.673085in}{1.096648in}}%
\pgfpathlineto{\pgfqpoint{1.673085in}{0.659680in}}%
\pgfpathclose%
\pgfusepath{fill}%
\end{pgfscope}%
\begin{pgfscope}%
\pgfpathrectangle{\pgfqpoint{1.673085in}{0.659680in}}{\pgfqpoint{0.436968in}{0.436968in}}%
\pgfusepath{clip}%
\pgfsys@transformshift{1.673085in}{0.659680in}%
\pgftext[left,bottom]{\includegraphics[interpolate=true,width=0.440000in,height=0.440000in]{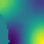}}%
\end{pgfscope}%
\begin{pgfscope}%
\definecolor{textcolor}{rgb}{0.000000,0.000000,0.000000}%
\pgfsetstrokecolor{textcolor}%
\pgfsetfillcolor{textcolor}%
\pgftext[x=1.891569in,y=1.120954in,,base]{\color{textcolor}{\rmfamily\fontsize{9.000000}{10.800000}\selectfont\catcode`\^=\active\def^{\ifmmode\sp\else\^{}\fi}\catcode`\%=\active\def%{\%}C 3}}%
\end{pgfscope}%
\begin{pgfscope}%
\pgfsetbuttcap%
\pgfsetmiterjoin%
\definecolor{currentfill}{rgb}{1.000000,1.000000,1.000000}%
\pgfsetfillcolor{currentfill}%
\pgfsetlinewidth{0.000000pt}%
\definecolor{currentstroke}{rgb}{0.000000,0.000000,0.000000}%
\pgfsetstrokecolor{currentstroke}%
\pgfsetstrokeopacity{0.000000}%
\pgfsetdash{}{0pt}%
\pgfpathmoveto{\pgfqpoint{1.673085in}{0.100000in}}%
\pgfpathlineto{\pgfqpoint{2.110053in}{0.100000in}}%
\pgfpathlineto{\pgfqpoint{2.110053in}{0.536968in}}%
\pgfpathlineto{\pgfqpoint{1.673085in}{0.536968in}}%
\pgfpathlineto{\pgfqpoint{1.673085in}{0.100000in}}%
\pgfpathclose%
\pgfusepath{fill}%
\end{pgfscope}%
\begin{pgfscope}%
\pgfpathrectangle{\pgfqpoint{1.673085in}{0.100000in}}{\pgfqpoint{0.436968in}{0.436968in}}%
\pgfusepath{clip}%
\pgfsys@transformshift{1.673085in}{0.100000in}%
\pgftext[left,bottom]{\includegraphics[interpolate=true,width=0.440000in,height=0.440000in]{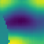}}%
\end{pgfscope}%
\begin{pgfscope}%
\definecolor{textcolor}{rgb}{0.000000,0.000000,0.000000}%
\pgfsetstrokecolor{textcolor}%
\pgfsetfillcolor{textcolor}%
\pgftext[x=1.891569in,y=0.561274in,,base]{\color{textcolor}{\rmfamily\fontsize{9.000000}{10.800000}\selectfont\catcode`\^=\active\def^{\ifmmode\sp\else\^{}\fi}\catcode`\%=\active\def%{\%}C 4}}%
\end{pgfscope}%
\begin{pgfscope}%
\pgfsetbuttcap%
\pgfsetmiterjoin%
\definecolor{currentfill}{rgb}{1.000000,1.000000,1.000000}%
\pgfsetfillcolor{currentfill}%
\pgfsetlinewidth{0.000000pt}%
\definecolor{currentstroke}{rgb}{0.000000,0.000000,0.000000}%
\pgfsetstrokecolor{currentstroke}%
\pgfsetstrokeopacity{0.000000}%
\pgfsetdash{}{0pt}%
\pgfpathmoveto{\pgfqpoint{1.148723in}{1.779040in}}%
\pgfpathlineto{\pgfqpoint{1.585691in}{1.779040in}}%
\pgfpathlineto{\pgfqpoint{1.585691in}{2.216008in}}%
\pgfpathlineto{\pgfqpoint{1.148723in}{2.216008in}}%
\pgfpathlineto{\pgfqpoint{1.148723in}{1.779040in}}%
\pgfpathclose%
\pgfusepath{fill}%
\end{pgfscope}%
\begin{pgfscope}%
\pgfpathrectangle{\pgfqpoint{1.148723in}{1.779040in}}{\pgfqpoint{0.436968in}{0.436968in}}%
\pgfusepath{clip}%
\pgfsys@transformshift{1.148723in}{1.779040in}%
\pgftext[left,bottom]{\includegraphics[interpolate=true,width=0.440000in,height=0.440000in]{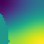}}%
\end{pgfscope}%
\begin{pgfscope}%
\definecolor{textcolor}{rgb}{0.000000,0.000000,0.000000}%
\pgfsetstrokecolor{textcolor}%
\pgfsetfillcolor{textcolor}%
\pgftext[x=1.367207in,y=2.240314in,,base]{\color{textcolor}{\rmfamily\fontsize{9.000000}{10.800000}\selectfont\catcode`\^=\active\def^{\ifmmode\sp\else\^{}\fi}\catcode`\%=\active\def%{\%}L 1}}%
\end{pgfscope}%
\begin{pgfscope}%
\pgfsetbuttcap%
\pgfsetmiterjoin%
\definecolor{currentfill}{rgb}{1.000000,1.000000,1.000000}%
\pgfsetfillcolor{currentfill}%
\pgfsetlinewidth{0.000000pt}%
\definecolor{currentstroke}{rgb}{0.000000,0.000000,0.000000}%
\pgfsetstrokecolor{currentstroke}%
\pgfsetstrokeopacity{0.000000}%
\pgfsetdash{}{0pt}%
\pgfpathmoveto{\pgfqpoint{1.148723in}{1.219360in}}%
\pgfpathlineto{\pgfqpoint{1.585691in}{1.219360in}}%
\pgfpathlineto{\pgfqpoint{1.585691in}{1.656328in}}%
\pgfpathlineto{\pgfqpoint{1.148723in}{1.656328in}}%
\pgfpathlineto{\pgfqpoint{1.148723in}{1.219360in}}%
\pgfpathclose%
\pgfusepath{fill}%
\end{pgfscope}%
\begin{pgfscope}%
\pgfpathrectangle{\pgfqpoint{1.148723in}{1.219360in}}{\pgfqpoint{0.436968in}{0.436968in}}%
\pgfusepath{clip}%
\pgfsys@transformshift{1.148723in}{1.219360in}%
\pgftext[left,bottom]{\includegraphics[interpolate=true,width=0.440000in,height=0.440000in]{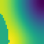}}%
\end{pgfscope}%
\begin{pgfscope}%
\definecolor{textcolor}{rgb}{0.000000,0.000000,0.000000}%
\pgfsetstrokecolor{textcolor}%
\pgfsetfillcolor{textcolor}%
\pgftext[x=1.367207in,y=1.680634in,,base]{\color{textcolor}{\rmfamily\fontsize{9.000000}{10.800000}\selectfont\catcode`\^=\active\def^{\ifmmode\sp\else\^{}\fi}\catcode`\%=\active\def%{\%}L 2}}%
\end{pgfscope}%
\begin{pgfscope}%
\pgfsetbuttcap%
\pgfsetmiterjoin%
\definecolor{currentfill}{rgb}{1.000000,1.000000,1.000000}%
\pgfsetfillcolor{currentfill}%
\pgfsetlinewidth{0.000000pt}%
\definecolor{currentstroke}{rgb}{0.000000,0.000000,0.000000}%
\pgfsetstrokecolor{currentstroke}%
\pgfsetstrokeopacity{0.000000}%
\pgfsetdash{}{0pt}%
\pgfpathmoveto{\pgfqpoint{1.148723in}{0.659680in}}%
\pgfpathlineto{\pgfqpoint{1.585691in}{0.659680in}}%
\pgfpathlineto{\pgfqpoint{1.585691in}{1.096648in}}%
\pgfpathlineto{\pgfqpoint{1.148723in}{1.096648in}}%
\pgfpathlineto{\pgfqpoint{1.148723in}{0.659680in}}%
\pgfpathclose%
\pgfusepath{fill}%
\end{pgfscope}%
\begin{pgfscope}%
\pgfpathrectangle{\pgfqpoint{1.148723in}{0.659680in}}{\pgfqpoint{0.436968in}{0.436968in}}%
\pgfusepath{clip}%
\pgfsys@transformshift{1.148723in}{0.659680in}%
\pgftext[left,bottom]{\includegraphics[interpolate=true,width=0.440000in,height=0.440000in]{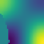}}%
\end{pgfscope}%
\begin{pgfscope}%
\definecolor{textcolor}{rgb}{0.000000,0.000000,0.000000}%
\pgfsetstrokecolor{textcolor}%
\pgfsetfillcolor{textcolor}%
\pgftext[x=1.367207in,y=1.120954in,,base]{\color{textcolor}{\rmfamily\fontsize{9.000000}{10.800000}\selectfont\catcode`\^=\active\def^{\ifmmode\sp\else\^{}\fi}\catcode`\%=\active\def%{\%}L 3}}%
\end{pgfscope}%
\begin{pgfscope}%
\pgfsetbuttcap%
\pgfsetmiterjoin%
\definecolor{currentfill}{rgb}{1.000000,1.000000,1.000000}%
\pgfsetfillcolor{currentfill}%
\pgfsetlinewidth{0.000000pt}%
\definecolor{currentstroke}{rgb}{0.000000,0.000000,0.000000}%
\pgfsetstrokecolor{currentstroke}%
\pgfsetstrokeopacity{0.000000}%
\pgfsetdash{}{0pt}%
\pgfpathmoveto{\pgfqpoint{1.148723in}{0.100000in}}%
\pgfpathlineto{\pgfqpoint{1.585691in}{0.100000in}}%
\pgfpathlineto{\pgfqpoint{1.585691in}{0.536968in}}%
\pgfpathlineto{\pgfqpoint{1.148723in}{0.536968in}}%
\pgfpathlineto{\pgfqpoint{1.148723in}{0.100000in}}%
\pgfpathclose%
\pgfusepath{fill}%
\end{pgfscope}%
\begin{pgfscope}%
\pgfpathrectangle{\pgfqpoint{1.148723in}{0.100000in}}{\pgfqpoint{0.436968in}{0.436968in}}%
\pgfusepath{clip}%
\pgfsys@transformshift{1.148723in}{0.100000in}%
\pgftext[left,bottom]{\includegraphics[interpolate=true,width=0.440000in,height=0.440000in]{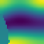}}%
\end{pgfscope}%
\begin{pgfscope}%
\definecolor{textcolor}{rgb}{0.000000,0.000000,0.000000}%
\pgfsetstrokecolor{textcolor}%
\pgfsetfillcolor{textcolor}%
\pgftext[x=1.367207in,y=0.561274in,,base]{\color{textcolor}{\rmfamily\fontsize{9.000000}{10.800000}\selectfont\catcode`\^=\active\def^{\ifmmode\sp\else\^{}\fi}\catcode`\%=\active\def%{\%}L 4}}%
\end{pgfscope}%
\begin{pgfscope}%
\pgfsetbuttcap%
\pgfsetmiterjoin%
\definecolor{currentfill}{rgb}{1.000000,1.000000,1.000000}%
\pgfsetfillcolor{currentfill}%
\pgfsetlinewidth{0.000000pt}%
\definecolor{currentstroke}{rgb}{0.000000,0.000000,0.000000}%
\pgfsetstrokecolor{currentstroke}%
\pgfsetstrokeopacity{0.000000}%
\pgfsetdash{}{0pt}%
\pgfpathmoveto{\pgfqpoint{2.197447in}{2.356379in}}%
\pgfpathlineto{\pgfqpoint{3.158777in}{2.356379in}}%
\pgfpathlineto{\pgfqpoint{3.158777in}{3.317709in}}%
\pgfpathlineto{\pgfqpoint{2.197447in}{3.317709in}}%
\pgfpathlineto{\pgfqpoint{2.197447in}{2.356379in}}%
\pgfpathclose%
\pgfusepath{fill}%
\end{pgfscope}%
\begin{pgfscope}%
\pgfpathrectangle{\pgfqpoint{2.197447in}{2.356379in}}{\pgfqpoint{0.961330in}{0.961330in}}%
\pgfusepath{clip}%
\pgfsys@transformshift{2.197447in}{2.356379in}%
\pgftext[left,bottom]{\includegraphics[interpolate=true,width=0.970000in,height=0.970000in]{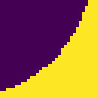}}%
\end{pgfscope}%
\begin{pgfscope}%
\definecolor{textcolor}{rgb}{0.000000,0.000000,0.000000}%
\pgfsetstrokecolor{textcolor}%
\pgfsetfillcolor{textcolor}%
\pgftext[x=2.678112in,y=3.342014in,,base]{\color{textcolor}{\rmfamily\fontsize{9.000000}{10.800000}\selectfont\catcode`\^=\active\def^{\ifmmode\sp\else\^{}\fi}\catcode`\%=\active\def%{\%}$\kappa_b = 10^{-3}$}}%
\end{pgfscope}%
\begin{pgfscope}%
\pgfsetbuttcap%
\pgfsetmiterjoin%
\definecolor{currentfill}{rgb}{1.000000,1.000000,1.000000}%
\pgfsetfillcolor{currentfill}%
\pgfsetlinewidth{0.000000pt}%
\definecolor{currentstroke}{rgb}{0.000000,0.000000,0.000000}%
\pgfsetstrokecolor{currentstroke}%
\pgfsetstrokeopacity{0.000000}%
\pgfsetdash{}{0pt}%
\pgfpathmoveto{\pgfqpoint{2.721809in}{1.779040in}}%
\pgfpathlineto{\pgfqpoint{3.158777in}{1.779040in}}%
\pgfpathlineto{\pgfqpoint{3.158777in}{2.216008in}}%
\pgfpathlineto{\pgfqpoint{2.721809in}{2.216008in}}%
\pgfpathlineto{\pgfqpoint{2.721809in}{1.779040in}}%
\pgfpathclose%
\pgfusepath{fill}%
\end{pgfscope}%
\begin{pgfscope}%
\pgfpathrectangle{\pgfqpoint{2.721809in}{1.779040in}}{\pgfqpoint{0.436968in}{0.436968in}}%
\pgfusepath{clip}%
\pgfsys@transformshift{2.721809in}{1.779040in}%
\pgftext[left,bottom]{\includegraphics[interpolate=true,width=0.440000in,height=0.440000in]{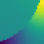}}%
\end{pgfscope}%
\begin{pgfscope}%
\definecolor{textcolor}{rgb}{0.000000,0.000000,0.000000}%
\pgfsetstrokecolor{textcolor}%
\pgfsetfillcolor{textcolor}%
\pgftext[x=2.940293in,y=2.240314in,,base]{\color{textcolor}{\rmfamily\fontsize{9.000000}{10.800000}\selectfont\catcode`\^=\active\def^{\ifmmode\sp\else\^{}\fi}\catcode`\%=\active\def%{\%}C 1}}%
\end{pgfscope}%
\begin{pgfscope}%
\pgfsetbuttcap%
\pgfsetmiterjoin%
\definecolor{currentfill}{rgb}{1.000000,1.000000,1.000000}%
\pgfsetfillcolor{currentfill}%
\pgfsetlinewidth{0.000000pt}%
\definecolor{currentstroke}{rgb}{0.000000,0.000000,0.000000}%
\pgfsetstrokecolor{currentstroke}%
\pgfsetstrokeopacity{0.000000}%
\pgfsetdash{}{0pt}%
\pgfpathmoveto{\pgfqpoint{2.721809in}{1.219360in}}%
\pgfpathlineto{\pgfqpoint{3.158777in}{1.219360in}}%
\pgfpathlineto{\pgfqpoint{3.158777in}{1.656328in}}%
\pgfpathlineto{\pgfqpoint{2.721809in}{1.656328in}}%
\pgfpathlineto{\pgfqpoint{2.721809in}{1.219360in}}%
\pgfpathclose%
\pgfusepath{fill}%
\end{pgfscope}%
\begin{pgfscope}%
\pgfpathrectangle{\pgfqpoint{2.721809in}{1.219360in}}{\pgfqpoint{0.436968in}{0.436968in}}%
\pgfusepath{clip}%
\pgfsys@transformshift{2.721809in}{1.219360in}%
\pgftext[left,bottom]{\includegraphics[interpolate=true,width=0.440000in,height=0.440000in]{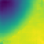}}%
\end{pgfscope}%
\begin{pgfscope}%
\definecolor{textcolor}{rgb}{0.000000,0.000000,0.000000}%
\pgfsetstrokecolor{textcolor}%
\pgfsetfillcolor{textcolor}%
\pgftext[x=2.940293in,y=1.680634in,,base]{\color{textcolor}{\rmfamily\fontsize{9.000000}{10.800000}\selectfont\catcode`\^=\active\def^{\ifmmode\sp\else\^{}\fi}\catcode`\%=\active\def%{\%}C 2}}%
\end{pgfscope}%
\begin{pgfscope}%
\pgfsetbuttcap%
\pgfsetmiterjoin%
\definecolor{currentfill}{rgb}{1.000000,1.000000,1.000000}%
\pgfsetfillcolor{currentfill}%
\pgfsetlinewidth{0.000000pt}%
\definecolor{currentstroke}{rgb}{0.000000,0.000000,0.000000}%
\pgfsetstrokecolor{currentstroke}%
\pgfsetstrokeopacity{0.000000}%
\pgfsetdash{}{0pt}%
\pgfpathmoveto{\pgfqpoint{2.721809in}{0.659680in}}%
\pgfpathlineto{\pgfqpoint{3.158777in}{0.659680in}}%
\pgfpathlineto{\pgfqpoint{3.158777in}{1.096648in}}%
\pgfpathlineto{\pgfqpoint{2.721809in}{1.096648in}}%
\pgfpathlineto{\pgfqpoint{2.721809in}{0.659680in}}%
\pgfpathclose%
\pgfusepath{fill}%
\end{pgfscope}%
\begin{pgfscope}%
\pgfpathrectangle{\pgfqpoint{2.721809in}{0.659680in}}{\pgfqpoint{0.436968in}{0.436968in}}%
\pgfusepath{clip}%
\pgfsys@transformshift{2.721809in}{0.659680in}%
\pgftext[left,bottom]{\includegraphics[interpolate=true,width=0.440000in,height=0.440000in]{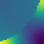}}%
\end{pgfscope}%
\begin{pgfscope}%
\definecolor{textcolor}{rgb}{0.000000,0.000000,0.000000}%
\pgfsetstrokecolor{textcolor}%
\pgfsetfillcolor{textcolor}%
\pgftext[x=2.940293in,y=1.120954in,,base]{\color{textcolor}{\rmfamily\fontsize{9.000000}{10.800000}\selectfont\catcode`\^=\active\def^{\ifmmode\sp\else\^{}\fi}\catcode`\%=\active\def%{\%}C 3}}%
\end{pgfscope}%
\begin{pgfscope}%
\pgfsetbuttcap%
\pgfsetmiterjoin%
\definecolor{currentfill}{rgb}{1.000000,1.000000,1.000000}%
\pgfsetfillcolor{currentfill}%
\pgfsetlinewidth{0.000000pt}%
\definecolor{currentstroke}{rgb}{0.000000,0.000000,0.000000}%
\pgfsetstrokecolor{currentstroke}%
\pgfsetstrokeopacity{0.000000}%
\pgfsetdash{}{0pt}%
\pgfpathmoveto{\pgfqpoint{2.721809in}{0.100000in}}%
\pgfpathlineto{\pgfqpoint{3.158777in}{0.100000in}}%
\pgfpathlineto{\pgfqpoint{3.158777in}{0.536968in}}%
\pgfpathlineto{\pgfqpoint{2.721809in}{0.536968in}}%
\pgfpathlineto{\pgfqpoint{2.721809in}{0.100000in}}%
\pgfpathclose%
\pgfusepath{fill}%
\end{pgfscope}%
\begin{pgfscope}%
\pgfpathrectangle{\pgfqpoint{2.721809in}{0.100000in}}{\pgfqpoint{0.436968in}{0.436968in}}%
\pgfusepath{clip}%
\pgfsys@transformshift{2.721809in}{0.100000in}%
\pgftext[left,bottom]{\includegraphics[interpolate=true,width=0.440000in,height=0.440000in]{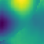}}%
\end{pgfscope}%
\begin{pgfscope}%
\definecolor{textcolor}{rgb}{0.000000,0.000000,0.000000}%
\pgfsetstrokecolor{textcolor}%
\pgfsetfillcolor{textcolor}%
\pgftext[x=2.940293in,y=0.561274in,,base]{\color{textcolor}{\rmfamily\fontsize{9.000000}{10.800000}\selectfont\catcode`\^=\active\def^{\ifmmode\sp\else\^{}\fi}\catcode`\%=\active\def%{\%}C 4}}%
\end{pgfscope}%
\begin{pgfscope}%
\pgfsetbuttcap%
\pgfsetmiterjoin%
\definecolor{currentfill}{rgb}{1.000000,1.000000,1.000000}%
\pgfsetfillcolor{currentfill}%
\pgfsetlinewidth{0.000000pt}%
\definecolor{currentstroke}{rgb}{0.000000,0.000000,0.000000}%
\pgfsetstrokecolor{currentstroke}%
\pgfsetstrokeopacity{0.000000}%
\pgfsetdash{}{0pt}%
\pgfpathmoveto{\pgfqpoint{2.197447in}{1.779040in}}%
\pgfpathlineto{\pgfqpoint{2.634415in}{1.779040in}}%
\pgfpathlineto{\pgfqpoint{2.634415in}{2.216008in}}%
\pgfpathlineto{\pgfqpoint{2.197447in}{2.216008in}}%
\pgfpathlineto{\pgfqpoint{2.197447in}{1.779040in}}%
\pgfpathclose%
\pgfusepath{fill}%
\end{pgfscope}%
\begin{pgfscope}%
\pgfpathrectangle{\pgfqpoint{2.197447in}{1.779040in}}{\pgfqpoint{0.436968in}{0.436968in}}%
\pgfusepath{clip}%
\pgfsys@transformshift{2.197447in}{1.779040in}%
\pgftext[left,bottom]{\includegraphics[interpolate=true,width=0.440000in,height=0.440000in]{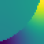}}%
\end{pgfscope}%
\begin{pgfscope}%
\definecolor{textcolor}{rgb}{0.000000,0.000000,0.000000}%
\pgfsetstrokecolor{textcolor}%
\pgfsetfillcolor{textcolor}%
\pgftext[x=2.415931in,y=2.240314in,,base]{\color{textcolor}{\rmfamily\fontsize{9.000000}{10.800000}\selectfont\catcode`\^=\active\def^{\ifmmode\sp\else\^{}\fi}\catcode`\%=\active\def%{\%}L 1}}%
\end{pgfscope}%
\begin{pgfscope}%
\pgfsetbuttcap%
\pgfsetmiterjoin%
\definecolor{currentfill}{rgb}{1.000000,1.000000,1.000000}%
\pgfsetfillcolor{currentfill}%
\pgfsetlinewidth{0.000000pt}%
\definecolor{currentstroke}{rgb}{0.000000,0.000000,0.000000}%
\pgfsetstrokecolor{currentstroke}%
\pgfsetstrokeopacity{0.000000}%
\pgfsetdash{}{0pt}%
\pgfpathmoveto{\pgfqpoint{2.197447in}{1.219360in}}%
\pgfpathlineto{\pgfqpoint{2.634415in}{1.219360in}}%
\pgfpathlineto{\pgfqpoint{2.634415in}{1.656328in}}%
\pgfpathlineto{\pgfqpoint{2.197447in}{1.656328in}}%
\pgfpathlineto{\pgfqpoint{2.197447in}{1.219360in}}%
\pgfpathclose%
\pgfusepath{fill}%
\end{pgfscope}%
\begin{pgfscope}%
\pgfpathrectangle{\pgfqpoint{2.197447in}{1.219360in}}{\pgfqpoint{0.436968in}{0.436968in}}%
\pgfusepath{clip}%
\pgfsys@transformshift{2.197447in}{1.219360in}%
\pgftext[left,bottom]{\includegraphics[interpolate=true,width=0.440000in,height=0.440000in]{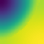}}%
\end{pgfscope}%
\begin{pgfscope}%
\definecolor{textcolor}{rgb}{0.000000,0.000000,0.000000}%
\pgfsetstrokecolor{textcolor}%
\pgfsetfillcolor{textcolor}%
\pgftext[x=2.415931in,y=1.680634in,,base]{\color{textcolor}{\rmfamily\fontsize{9.000000}{10.800000}\selectfont\catcode`\^=\active\def^{\ifmmode\sp\else\^{}\fi}\catcode`\%=\active\def%{\%}L 2}}%
\end{pgfscope}%
\begin{pgfscope}%
\pgfsetbuttcap%
\pgfsetmiterjoin%
\definecolor{currentfill}{rgb}{1.000000,1.000000,1.000000}%
\pgfsetfillcolor{currentfill}%
\pgfsetlinewidth{0.000000pt}%
\definecolor{currentstroke}{rgb}{0.000000,0.000000,0.000000}%
\pgfsetstrokecolor{currentstroke}%
\pgfsetstrokeopacity{0.000000}%
\pgfsetdash{}{0pt}%
\pgfpathmoveto{\pgfqpoint{2.197447in}{0.659680in}}%
\pgfpathlineto{\pgfqpoint{2.634415in}{0.659680in}}%
\pgfpathlineto{\pgfqpoint{2.634415in}{1.096648in}}%
\pgfpathlineto{\pgfqpoint{2.197447in}{1.096648in}}%
\pgfpathlineto{\pgfqpoint{2.197447in}{0.659680in}}%
\pgfpathclose%
\pgfusepath{fill}%
\end{pgfscope}%
\begin{pgfscope}%
\pgfpathrectangle{\pgfqpoint{2.197447in}{0.659680in}}{\pgfqpoint{0.436968in}{0.436968in}}%
\pgfusepath{clip}%
\pgfsys@transformshift{2.197447in}{0.659680in}%
\pgftext[left,bottom]{\includegraphics[interpolate=true,width=0.440000in,height=0.440000in]{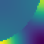}}%
\end{pgfscope}%
\begin{pgfscope}%
\definecolor{textcolor}{rgb}{0.000000,0.000000,0.000000}%
\pgfsetstrokecolor{textcolor}%
\pgfsetfillcolor{textcolor}%
\pgftext[x=2.415931in,y=1.120954in,,base]{\color{textcolor}{\rmfamily\fontsize{9.000000}{10.800000}\selectfont\catcode`\^=\active\def^{\ifmmode\sp\else\^{}\fi}\catcode`\%=\active\def%{\%}L 3}}%
\end{pgfscope}%
\begin{pgfscope}%
\pgfsetbuttcap%
\pgfsetmiterjoin%
\definecolor{currentfill}{rgb}{1.000000,1.000000,1.000000}%
\pgfsetfillcolor{currentfill}%
\pgfsetlinewidth{0.000000pt}%
\definecolor{currentstroke}{rgb}{0.000000,0.000000,0.000000}%
\pgfsetstrokecolor{currentstroke}%
\pgfsetstrokeopacity{0.000000}%
\pgfsetdash{}{0pt}%
\pgfpathmoveto{\pgfqpoint{2.197447in}{0.100000in}}%
\pgfpathlineto{\pgfqpoint{2.634415in}{0.100000in}}%
\pgfpathlineto{\pgfqpoint{2.634415in}{0.536968in}}%
\pgfpathlineto{\pgfqpoint{2.197447in}{0.536968in}}%
\pgfpathlineto{\pgfqpoint{2.197447in}{0.100000in}}%
\pgfpathclose%
\pgfusepath{fill}%
\end{pgfscope}%
\begin{pgfscope}%
\pgfpathrectangle{\pgfqpoint{2.197447in}{0.100000in}}{\pgfqpoint{0.436968in}{0.436968in}}%
\pgfusepath{clip}%
\pgfsys@transformshift{2.197447in}{0.100000in}%
\pgftext[left,bottom]{\includegraphics[interpolate=true,width=0.440000in,height=0.440000in]{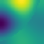}}%
\end{pgfscope}%
\begin{pgfscope}%
\definecolor{textcolor}{rgb}{0.000000,0.000000,0.000000}%
\pgfsetstrokecolor{textcolor}%
\pgfsetfillcolor{textcolor}%
\pgftext[x=2.415931in,y=0.561274in,,base]{\color{textcolor}{\rmfamily\fontsize{9.000000}{10.800000}\selectfont\catcode`\^=\active\def^{\ifmmode\sp\else\^{}\fi}\catcode`\%=\active\def%{\%}L 4}}%
\end{pgfscope}%
\begin{pgfscope}%
\pgfsetbuttcap%
\pgfsetmiterjoin%
\definecolor{currentfill}{rgb}{1.000000,1.000000,1.000000}%
\pgfsetfillcolor{currentfill}%
\pgfsetlinewidth{0.000000pt}%
\definecolor{currentstroke}{rgb}{0.000000,0.000000,0.000000}%
\pgfsetstrokecolor{currentstroke}%
\pgfsetstrokeopacity{0.000000}%
\pgfsetdash{}{0pt}%
\pgfpathmoveto{\pgfqpoint{3.246170in}{2.356379in}}%
\pgfpathlineto{\pgfqpoint{4.207500in}{2.356379in}}%
\pgfpathlineto{\pgfqpoint{4.207500in}{3.317709in}}%
\pgfpathlineto{\pgfqpoint{3.246170in}{3.317709in}}%
\pgfpathlineto{\pgfqpoint{3.246170in}{2.356379in}}%
\pgfpathclose%
\pgfusepath{fill}%
\end{pgfscope}%
\begin{pgfscope}%
\pgfpathrectangle{\pgfqpoint{3.246170in}{2.356379in}}{\pgfqpoint{0.961330in}{0.961330in}}%
\pgfusepath{clip}%
\pgfsys@transformshift{3.246170in}{2.356379in}%
\pgftext[left,bottom]{\includegraphics[interpolate=true,width=0.970000in,height=0.970000in]{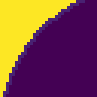}}%
\end{pgfscope}%
\begin{pgfscope}%
\definecolor{textcolor}{rgb}{0.000000,0.000000,0.000000}%
\pgfsetstrokecolor{textcolor}%
\pgfsetfillcolor{textcolor}%
\pgftext[x=3.726835in,y=3.342014in,,base]{\color{textcolor}{\rmfamily\fontsize{9.000000}{10.800000}\selectfont\catcode`\^=\active\def^{\ifmmode\sp\else\^{}\fi}\catcode`\%=\active\def%{\%}$\kappa_b = 10^{-2}$}}%
\end{pgfscope}%
\begin{pgfscope}%
\pgfsetbuttcap%
\pgfsetmiterjoin%
\definecolor{currentfill}{rgb}{1.000000,1.000000,1.000000}%
\pgfsetfillcolor{currentfill}%
\pgfsetlinewidth{0.000000pt}%
\definecolor{currentstroke}{rgb}{0.000000,0.000000,0.000000}%
\pgfsetstrokecolor{currentstroke}%
\pgfsetstrokeopacity{0.000000}%
\pgfsetdash{}{0pt}%
\pgfpathmoveto{\pgfqpoint{3.770532in}{1.779040in}}%
\pgfpathlineto{\pgfqpoint{4.207500in}{1.779040in}}%
\pgfpathlineto{\pgfqpoint{4.207500in}{2.216008in}}%
\pgfpathlineto{\pgfqpoint{3.770532in}{2.216008in}}%
\pgfpathlineto{\pgfqpoint{3.770532in}{1.779040in}}%
\pgfpathclose%
\pgfusepath{fill}%
\end{pgfscope}%
\begin{pgfscope}%
\pgfpathrectangle{\pgfqpoint{3.770532in}{1.779040in}}{\pgfqpoint{0.436968in}{0.436968in}}%
\pgfusepath{clip}%
\pgfsys@transformshift{3.770532in}{1.779040in}%
\pgftext[left,bottom]{\includegraphics[interpolate=true,width=0.440000in,height=0.440000in]{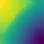}}%
\end{pgfscope}%
\begin{pgfscope}%
\definecolor{textcolor}{rgb}{0.000000,0.000000,0.000000}%
\pgfsetstrokecolor{textcolor}%
\pgfsetfillcolor{textcolor}%
\pgftext[x=3.989016in,y=2.240314in,,base]{\color{textcolor}{\rmfamily\fontsize{9.000000}{10.800000}\selectfont\catcode`\^=\active\def^{\ifmmode\sp\else\^{}\fi}\catcode`\%=\active\def%{\%}C 1}}%
\end{pgfscope}%
\begin{pgfscope}%
\pgfsetbuttcap%
\pgfsetmiterjoin%
\definecolor{currentfill}{rgb}{1.000000,1.000000,1.000000}%
\pgfsetfillcolor{currentfill}%
\pgfsetlinewidth{0.000000pt}%
\definecolor{currentstroke}{rgb}{0.000000,0.000000,0.000000}%
\pgfsetstrokecolor{currentstroke}%
\pgfsetstrokeopacity{0.000000}%
\pgfsetdash{}{0pt}%
\pgfpathmoveto{\pgfqpoint{3.770532in}{1.219360in}}%
\pgfpathlineto{\pgfqpoint{4.207500in}{1.219360in}}%
\pgfpathlineto{\pgfqpoint{4.207500in}{1.656328in}}%
\pgfpathlineto{\pgfqpoint{3.770532in}{1.656328in}}%
\pgfpathlineto{\pgfqpoint{3.770532in}{1.219360in}}%
\pgfpathclose%
\pgfusepath{fill}%
\end{pgfscope}%
\begin{pgfscope}%
\pgfpathrectangle{\pgfqpoint{3.770532in}{1.219360in}}{\pgfqpoint{0.436968in}{0.436968in}}%
\pgfusepath{clip}%
\pgfsys@transformshift{3.770532in}{1.219360in}%
\pgftext[left,bottom]{\includegraphics[interpolate=true,width=0.440000in,height=0.440000in]{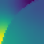}}%
\end{pgfscope}%
\begin{pgfscope}%
\definecolor{textcolor}{rgb}{0.000000,0.000000,0.000000}%
\pgfsetstrokecolor{textcolor}%
\pgfsetfillcolor{textcolor}%
\pgftext[x=3.989016in,y=1.680634in,,base]{\color{textcolor}{\rmfamily\fontsize{9.000000}{10.800000}\selectfont\catcode`\^=\active\def^{\ifmmode\sp\else\^{}\fi}\catcode`\%=\active\def%{\%}C 2}}%
\end{pgfscope}%
\begin{pgfscope}%
\pgfsetbuttcap%
\pgfsetmiterjoin%
\definecolor{currentfill}{rgb}{1.000000,1.000000,1.000000}%
\pgfsetfillcolor{currentfill}%
\pgfsetlinewidth{0.000000pt}%
\definecolor{currentstroke}{rgb}{0.000000,0.000000,0.000000}%
\pgfsetstrokecolor{currentstroke}%
\pgfsetstrokeopacity{0.000000}%
\pgfsetdash{}{0pt}%
\pgfpathmoveto{\pgfqpoint{3.770532in}{0.659680in}}%
\pgfpathlineto{\pgfqpoint{4.207500in}{0.659680in}}%
\pgfpathlineto{\pgfqpoint{4.207500in}{1.096648in}}%
\pgfpathlineto{\pgfqpoint{3.770532in}{1.096648in}}%
\pgfpathlineto{\pgfqpoint{3.770532in}{0.659680in}}%
\pgfpathclose%
\pgfusepath{fill}%
\end{pgfscope}%
\begin{pgfscope}%
\pgfpathrectangle{\pgfqpoint{3.770532in}{0.659680in}}{\pgfqpoint{0.436968in}{0.436968in}}%
\pgfusepath{clip}%
\pgfsys@transformshift{3.770532in}{0.659680in}%
\pgftext[left,bottom]{\includegraphics[interpolate=true,width=0.440000in,height=0.440000in]{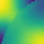}}%
\end{pgfscope}%
\begin{pgfscope}%
\definecolor{textcolor}{rgb}{0.000000,0.000000,0.000000}%
\pgfsetstrokecolor{textcolor}%
\pgfsetfillcolor{textcolor}%
\pgftext[x=3.989016in,y=1.120954in,,base]{\color{textcolor}{\rmfamily\fontsize{9.000000}{10.800000}\selectfont\catcode`\^=\active\def^{\ifmmode\sp\else\^{}\fi}\catcode`\%=\active\def%{\%}C 3}}%
\end{pgfscope}%
\begin{pgfscope}%
\pgfsetbuttcap%
\pgfsetmiterjoin%
\definecolor{currentfill}{rgb}{1.000000,1.000000,1.000000}%
\pgfsetfillcolor{currentfill}%
\pgfsetlinewidth{0.000000pt}%
\definecolor{currentstroke}{rgb}{0.000000,0.000000,0.000000}%
\pgfsetstrokecolor{currentstroke}%
\pgfsetstrokeopacity{0.000000}%
\pgfsetdash{}{0pt}%
\pgfpathmoveto{\pgfqpoint{3.770532in}{0.100000in}}%
\pgfpathlineto{\pgfqpoint{4.207500in}{0.100000in}}%
\pgfpathlineto{\pgfqpoint{4.207500in}{0.536968in}}%
\pgfpathlineto{\pgfqpoint{3.770532in}{0.536968in}}%
\pgfpathlineto{\pgfqpoint{3.770532in}{0.100000in}}%
\pgfpathclose%
\pgfusepath{fill}%
\end{pgfscope}%
\begin{pgfscope}%
\pgfpathrectangle{\pgfqpoint{3.770532in}{0.100000in}}{\pgfqpoint{0.436968in}{0.436968in}}%
\pgfusepath{clip}%
\pgfsys@transformshift{3.770532in}{0.100000in}%
\pgftext[left,bottom]{\includegraphics[interpolate=true,width=0.440000in,height=0.440000in]{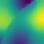}}%
\end{pgfscope}%
\begin{pgfscope}%
\definecolor{textcolor}{rgb}{0.000000,0.000000,0.000000}%
\pgfsetstrokecolor{textcolor}%
\pgfsetfillcolor{textcolor}%
\pgftext[x=3.989016in,y=0.561274in,,base]{\color{textcolor}{\rmfamily\fontsize{9.000000}{10.800000}\selectfont\catcode`\^=\active\def^{\ifmmode\sp\else\^{}\fi}\catcode`\%=\active\def%{\%}C 4}}%
\end{pgfscope}%
\begin{pgfscope}%
\pgfsetbuttcap%
\pgfsetmiterjoin%
\definecolor{currentfill}{rgb}{1.000000,1.000000,1.000000}%
\pgfsetfillcolor{currentfill}%
\pgfsetlinewidth{0.000000pt}%
\definecolor{currentstroke}{rgb}{0.000000,0.000000,0.000000}%
\pgfsetstrokecolor{currentstroke}%
\pgfsetstrokeopacity{0.000000}%
\pgfsetdash{}{0pt}%
\pgfpathmoveto{\pgfqpoint{3.246170in}{1.779040in}}%
\pgfpathlineto{\pgfqpoint{3.683138in}{1.779040in}}%
\pgfpathlineto{\pgfqpoint{3.683138in}{2.216008in}}%
\pgfpathlineto{\pgfqpoint{3.246170in}{2.216008in}}%
\pgfpathlineto{\pgfqpoint{3.246170in}{1.779040in}}%
\pgfpathclose%
\pgfusepath{fill}%
\end{pgfscope}%
\begin{pgfscope}%
\pgfpathrectangle{\pgfqpoint{3.246170in}{1.779040in}}{\pgfqpoint{0.436968in}{0.436968in}}%
\pgfusepath{clip}%
\pgfsys@transformshift{3.246170in}{1.779040in}%
\pgftext[left,bottom]{\includegraphics[interpolate=true,width=0.440000in,height=0.440000in]{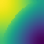}}%
\end{pgfscope}%
\begin{pgfscope}%
\definecolor{textcolor}{rgb}{0.000000,0.000000,0.000000}%
\pgfsetstrokecolor{textcolor}%
\pgfsetfillcolor{textcolor}%
\pgftext[x=3.464654in,y=2.240314in,,base]{\color{textcolor}{\rmfamily\fontsize{9.000000}{10.800000}\selectfont\catcode`\^=\active\def^{\ifmmode\sp\else\^{}\fi}\catcode`\%=\active\def%{\%}L 1}}%
\end{pgfscope}%
\begin{pgfscope}%
\pgfsetbuttcap%
\pgfsetmiterjoin%
\definecolor{currentfill}{rgb}{1.000000,1.000000,1.000000}%
\pgfsetfillcolor{currentfill}%
\pgfsetlinewidth{0.000000pt}%
\definecolor{currentstroke}{rgb}{0.000000,0.000000,0.000000}%
\pgfsetstrokecolor{currentstroke}%
\pgfsetstrokeopacity{0.000000}%
\pgfsetdash{}{0pt}%
\pgfpathmoveto{\pgfqpoint{3.246170in}{1.219360in}}%
\pgfpathlineto{\pgfqpoint{3.683138in}{1.219360in}}%
\pgfpathlineto{\pgfqpoint{3.683138in}{1.656328in}}%
\pgfpathlineto{\pgfqpoint{3.246170in}{1.656328in}}%
\pgfpathlineto{\pgfqpoint{3.246170in}{1.219360in}}%
\pgfpathclose%
\pgfusepath{fill}%
\end{pgfscope}%
\begin{pgfscope}%
\pgfpathrectangle{\pgfqpoint{3.246170in}{1.219360in}}{\pgfqpoint{0.436968in}{0.436968in}}%
\pgfusepath{clip}%
\pgfsys@transformshift{3.246170in}{1.219360in}%
\pgftext[left,bottom]{\includegraphics[interpolate=true,width=0.440000in,height=0.440000in]{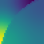}}%
\end{pgfscope}%
\begin{pgfscope}%
\definecolor{textcolor}{rgb}{0.000000,0.000000,0.000000}%
\pgfsetstrokecolor{textcolor}%
\pgfsetfillcolor{textcolor}%
\pgftext[x=3.464654in,y=1.680634in,,base]{\color{textcolor}{\rmfamily\fontsize{9.000000}{10.800000}\selectfont\catcode`\^=\active\def^{\ifmmode\sp\else\^{}\fi}\catcode`\%=\active\def%{\%}L 2}}%
\end{pgfscope}%
\begin{pgfscope}%
\pgfsetbuttcap%
\pgfsetmiterjoin%
\definecolor{currentfill}{rgb}{1.000000,1.000000,1.000000}%
\pgfsetfillcolor{currentfill}%
\pgfsetlinewidth{0.000000pt}%
\definecolor{currentstroke}{rgb}{0.000000,0.000000,0.000000}%
\pgfsetstrokecolor{currentstroke}%
\pgfsetstrokeopacity{0.000000}%
\pgfsetdash{}{0pt}%
\pgfpathmoveto{\pgfqpoint{3.246170in}{0.659680in}}%
\pgfpathlineto{\pgfqpoint{3.683138in}{0.659680in}}%
\pgfpathlineto{\pgfqpoint{3.683138in}{1.096648in}}%
\pgfpathlineto{\pgfqpoint{3.246170in}{1.096648in}}%
\pgfpathlineto{\pgfqpoint{3.246170in}{0.659680in}}%
\pgfpathclose%
\pgfusepath{fill}%
\end{pgfscope}%
\begin{pgfscope}%
\pgfpathrectangle{\pgfqpoint{3.246170in}{0.659680in}}{\pgfqpoint{0.436968in}{0.436968in}}%
\pgfusepath{clip}%
\pgfsys@transformshift{3.246170in}{0.659680in}%
\pgftext[left,bottom]{\includegraphics[interpolate=true,width=0.440000in,height=0.440000in]{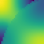}}%
\end{pgfscope}%
\begin{pgfscope}%
\definecolor{textcolor}{rgb}{0.000000,0.000000,0.000000}%
\pgfsetstrokecolor{textcolor}%
\pgfsetfillcolor{textcolor}%
\pgftext[x=3.464654in,y=1.120954in,,base]{\color{textcolor}{\rmfamily\fontsize{9.000000}{10.800000}\selectfont\catcode`\^=\active\def^{\ifmmode\sp\else\^{}\fi}\catcode`\%=\active\def%{\%}L 3}}%
\end{pgfscope}%
\begin{pgfscope}%
\pgfsetbuttcap%
\pgfsetmiterjoin%
\definecolor{currentfill}{rgb}{1.000000,1.000000,1.000000}%
\pgfsetfillcolor{currentfill}%
\pgfsetlinewidth{0.000000pt}%
\definecolor{currentstroke}{rgb}{0.000000,0.000000,0.000000}%
\pgfsetstrokecolor{currentstroke}%
\pgfsetstrokeopacity{0.000000}%
\pgfsetdash{}{0pt}%
\pgfpathmoveto{\pgfqpoint{3.246170in}{0.100000in}}%
\pgfpathlineto{\pgfqpoint{3.683138in}{0.100000in}}%
\pgfpathlineto{\pgfqpoint{3.683138in}{0.536968in}}%
\pgfpathlineto{\pgfqpoint{3.246170in}{0.536968in}}%
\pgfpathlineto{\pgfqpoint{3.246170in}{0.100000in}}%
\pgfpathclose%
\pgfusepath{fill}%
\end{pgfscope}%
\begin{pgfscope}%
\pgfpathrectangle{\pgfqpoint{3.246170in}{0.100000in}}{\pgfqpoint{0.436968in}{0.436968in}}%
\pgfusepath{clip}%
\pgfsys@transformshift{3.246170in}{0.100000in}%
\pgftext[left,bottom]{\includegraphics[interpolate=true,width=0.440000in,height=0.440000in]{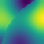}}%
\end{pgfscope}%
\begin{pgfscope}%
\definecolor{textcolor}{rgb}{0.000000,0.000000,0.000000}%
\pgfsetstrokecolor{textcolor}%
\pgfsetfillcolor{textcolor}%
\pgftext[x=3.464654in,y=0.561274in,,base]{\color{textcolor}{\rmfamily\fontsize{9.000000}{10.800000}\selectfont\catcode`\^=\active\def^{\ifmmode\sp\else\^{}\fi}\catcode`\%=\active\def%{\%}L 4}}%
\end{pgfscope}%
\end{pgfpicture}%
\makeatother%
\endgroup%